\documentclass[10pt]{amsart}
\usepackage[a4paper,top=2cm,bottom=2cm,left=2cm,right=2cm,marginparwidth=1.75cm]{geometry}
\usepackage[utf8]{inputenc} 
\usepackage[T1]{fontenc} 

\usepackage[english]{babel} 
\usepackage{booktabs,longtable,array}
\usepackage{graphicx} 
\usepackage[colorlinks]{hyperref} 
\usepackage{amsfonts, amsmath, amssymb, amsthm} 
\usepackage{stmaryrd}
\usepackage{bm}
 \usepackage{xcolor}
\usepackage{tikz}
\usepackage{mathrsfs}
\usepackage[colorinlistoftodos,bordercolor=orange,backgroundcolor=orange!20,linecolor=orange,textsize=scriptsize]{todonotes}

\newcounter{aah}
\numberwithin{equation}{section}
\counterwithin{aah}{section}
\newtheorem{thm}[aah]{Theorem}
\newtheorem{prop}[aah]{Proposition}
\newtheorem{coro}[aah]{Corollary}
\newtheorem{lem}[aah]{Lemma}

\theoremstyle{definition}

\newtheorem{defi}[aah]{Definition}
\newtheorem{propdef}[aah]{Definition-Proposition}

\theoremstyle{remark}

\newtheorem{exemple}[aah]{Example}
\newtheorem{rem}[aah]{Remark}
\title[Curvature condition on the Wasserstein space]{Curvature condition on the Wasserstein space over a closed Riemannian manifold}
\author{Aymeric Martin}
\date{November 2026}

\DeclareMathOperator{\Ric}{Ric}

\newcommand{\di}{\mathscr{D}}
\newcommand{\D}{\bar{D}}
\newcommand{\T}{\Bar{T}}
\newcommand{\grad}{\Bar{\nabla}}
\renewcommand{\P}{\mathscr{P}}
\newcommand{\SN}{\underset{i=1}{\overset{N}{\sum}}}

\newcommand{\R}{\overline{R}}
\newcommand{\C}{\mathcal{C}}

\renewcommand{\L}{\bar{\mathcal{L}}}
\newcommand{\A}{\textbf{A}}

\renewcommand{\H}{\overline{\mathcal{H}}}
\newcommand{\W}{\bm W}

\newcommand{\N}{\mathcal{N}}
\newcommand{\vol}{\mathrm{vol}}
\newcommand{\h}{\mathfrak{h}}

\newcommand{\M}{\mathcal{M}}
\newcommand{\V}{\mathscr{V}}
\renewcommand{\H}{\mathscr{H}}
\newcommand{\Inv}{\mathrm{GL}_d}
\renewcommand{\d}{\, \mathrm{d}}
\newcommand{\II}{\mathrm{I \! I}}
\newcommand{\dw}{\bar{\mathrm{d}}}
\newcommand{\dd}{\bm \d}
\newcommand{\he}{\mathrm{H}}

\newcommand{\se}{\bar{\mathrm{P}}}

\begin{document}

\begin{abstract}
	We establish regularity results for the solution flow of Stochastic Differential Equations (SDE) on the Wasserstein space $\mathscr{P}(M)$ over a closed Riemannian manifold, under suitable assumptions on the driving vector fields. We derive the SDE governing the associated tangent process and prove its well-posedness by means of the framework developped in \cite{Martin2025} for Covariant SDEs on Inverse Limit Hilbert vector bundles over the diffeomorphism group. We introduce the $\mathrm{IWC}(K)$ condition, a geometric curvature-type condition which generalizes, in the setting of SDEs on $\P(M)$, both the $\mathrm{H.C.K.}$ condition introduced in \cite{cattiaux:hal-00805299} and the ergodicity condition of \cite[Section 4]{wang2020}. Finally, combining this condition with the theory of entropic interpolations (see \cite{Leonard2014}), we establish several consequences of $\mathrm{IWC}(K)$, including control of the distance between two trajectories, gradient estimates for the associated Markov semigroup and a local Poincaré inequality. In the case $K >0$, we further obtain the existence and uniqueness of an ergodic invariant measure for the associated semigroup.
\end{abstract}

\maketitle

\tableofcontents

\section{Introduction}

	The question of sufficient conditions to ensure stability of Stochastic Differential Equations (SDEs for short) on Riemannian manifolds has long been a subject of study. In their celebrated paper \cite{BakryEmery1985}, Bakry and Emery studied general Markov processes and introduced  the so-called Curvature-Dimension condition to study the stability of such processes. This condition is highly based on analytical features such as the infinitesimal generator, the Carré du Champ or the Carré du champ itéré and has proven its worth by being very fruitful, see \cite[Chapter 5]{logsob2000} for applications to a lot of functional inequalities questions, \cite[Part III]{villani2009optimal} for a generalization to metric spaces for instance, \cite{baudoin2026} for a generalization to sub-Laplacian on Riemannian foliations for instance. In \cite{Wang1997}, Wang shows that most of the consequences of the Curvature-Dimension condition can be recovered by using a coupling by parallel displacement in the case of the Brownian motion, see also \cite{vonRenesseSturm2005} for a presentation of the equivalence of all these properties. For a generic diffusion on a Riemannian manifold, the computation of the Carré du champ itéré can turn out to be difficult and it may be more difficult to demonstrate the impact of the geometry of the manifold on the stability of the diffusions. 

	In \cite{cattiaux:hal-00805299}, Cattiaux and Guillin investigated the stability of generic diffusions on $\mathbb{R}^d$ satisfying the so-called Semi Log-Concavity property. The idea is to use the synchronous coupling instead of the parallel displacement one. In some sense, these results are more general since, to perform parallel displacement coupling we need to consider elliptic diffusions while we can dispense with this condition in the synchronous case. However, the constants obtained for the diffusions to which both criteria can be applied, for instance Brownian motions, are less sharp with the latter method, and it is not always possible to recover all the functional inequalities; this is the case, for instance, for the log-Sobolev inequality.

	In his seminal paper \cite{Otto31012001}, Otto introduces an infinite-dimensional Riemannian structure on the Wasserstein space $\P$, the counterpart of the Riemannian metric being the so-called Otto metric. This Riemannian structure has been developped in several directions, see \cite{ambrosioGradientFlowsMetric2008}, \cite{lott2007geometriccalculationswassersteinspace}, \cite{villani2009optimal}. The comparison with a standard Riemannian manifold is very limited because of the singularity of $\P$. For instance, the tangent space to a measure $\sum_{i=1}^N \alpha_i \delta_{x_i}$ is of dimension $N\cdot \mathrm{dim}(M)$, whereas the tangent space at an absolutely continuous measure $\mu = \rho \vol$ is an infinite-dimensional Hilbert space. Informally, the fact that such measures can have radically different tangent spaces while being abritrarily close to one another rules out the existence of local charts.
	
	If we restrict the study to the space $\P_\infty$ of smooth positive measures, we can get rid of these problems and use the Riemannian submersion structure of the diffeomorphism group $\mathscr{D}$ onto $\P_\infty$ to obtain a nice geometric description of this space, see \cite{EbinMarsden1970}, \cite{Martin2025}. The question of the geodesic convexity of such a space is still open and seems to be highly dependent on the geometry of the base manifold. Indeed, this is a well known fact that there exists a unique geodesic (called the displacement interpolation in Optimal Transport theory) joining two measures in $\P_\infty$, see \cite{mccann2001polar}. Under additional geometric assumptions, such as the $\mathrm{MTW}(K)$ condition, it can be shown that the associated optimal transport map is smooth, see \cite{Fig2010}, \cite{MaTrudingerWang2005}, \cite{LoeperVillani2010}. However, the fact that the geodesic remains in $\P_\infty$ is not clear, even under these assumptions. To overcome this problem, one can use the Schrödinger problem theory, see \cite{Leonard2014} for a good presentation of the problem and the main results. It is shown in \cite{GentilIvanLeonardRipaniTamini2020} that the solution to the $\varepsilon$-Schrödinger problem between two measures in $\P_\infty$, the so-called $\varepsilon$-entropic interpolation, is a nice Lipschitz approximation (in the sense of Otto geometry) of the displacement interpolation that remains in $\P_\infty$.

	The study of SDEs on the Wasserstein space $\P$ can be motived, on the one hand, by the theory of mean field games with common noise, see \cite{carmona2018probabilistic2} for instance, and, on the other hand, by the study of Stochastic Partial Differential Equations (SPDEs for short), see \cite{gassiat2022long}. In \cite{wang2020}, Wang shows well-posedness and some properties of SDEs on $\P(\mathbb{R}^d)$, and in \cite{ding2023stochasticdifferentialequationsstochastic}, Ding, Fang and Li investigate a class of SDEs on $\P(M)$, for $M$ a closed manifold. This work has been pursued in \cite{Martin2025} by generalizing the class of admissible driving vector fields and by presenting a way to study Covariant SDEs on the tangent bundle $T\P_\infty$ inspired by the infinite-dimensional calculus formalism described in \cite{Elworthy1982} and the Inverse Limit Hilbert Lie Group structure of $\mathscr{D}$ presented in \cite{EbinMarsden1970}, \cite{omori_group_1970}.

	This paper is dedicated to a generalization of the work of Cattiaux and Guillin in \cite{cattiaux:hal-00805299} in the case of the Wasserstein space $\P$ through a study of the synchronous coupling. This choice of coupling is due to, at least, two factors: the first one is that the parallel coupling is way more difficult because of the infinite-dimensional framework, for instance, the associated frame bundle is an inverse limit of Banach manifolds, this setting makes the stochastic calculus way harder than on Hilbert manifolds. The second one is that we need an infinite-dimensional counterpart of the Brownian motion. This is well known by the works of Dello-Schiavo and Sturm, see \cite{DelloSchiavo2022}, \cite{schiavo2025massiveparticlesystemswasserstein}, \cite{Sturm2024}, that there exist counterparts to the Brownian motion, but these processes are not canonical and are difficult to use in practice. \vspace{10pt}

	\noindent \textbf{Acknowledgements.} The author wishes to thank his PhD advisors, Marc Arnaudon and Michel Bonnefont, for their invaluable guidance, continuous support, and thoughtful advice throughout the course of this work. The author would also like to express his gratitude to Christian Leonard for taking the time to discuss his work. This discussion was very enriching and had a significant impact on the writing of this article.

	\subsection{Stochastic calculus and Riemannian geometry prerequisites}

	Let $(M,\langle \bullet,\bullet\rangle)$ be a complete Riemannian manifold, $(\Omega,\mathscr{F},(\mathscr{F}_t)_t, \mathbb{P})$ be a filtered probability space carrying a $\mathbb{R}^N$-valued Brownian motion $W_t =(W_t^1,\dots W_t^N)$ in the canonical basis. Thanks to the Riemannian metric $\langle \bullet,\bullet \rangle$, we will always identify $TM$ with $T^*M$ in the following. We denote by $\nabla$ both the associated Levi-Civita connection and gradient. The Riemannian Hessian of a function $f \in \C^2(M)$ is given by 
	
	$$
	\he f(A,B) = A B f - (\nabla_A B)f.
	$$
	Let $A_0,A_1, \dots,A_N \in \Gamma(TM)$. A stochastic process $(X_t)_{t \geq 0}$ is said to be solution to
	
	\begin{equation*}
\d^\nabla X_t = \SN A_i(X_t) \d W^i_t + A_0(X_t) \d t, \quad X_0 = x
	\end{equation*}
	if for all $f \in \C^2(M)$,

	$$
	f(X_t)= f(x) + \SN \int_0^t A_i f (X_s) \d W^i_s + \int_0^t A_0f(X_s) \d s + \frac{1}{2} \SN \int_0^t \he f(A_i(X_s),A_i(X_s)) \d s. 
	$$
	There always exists a unique solution to such an SDE under suitable assumptions, see for instance \cite{hsu2002stochastic} or \cite{kunita1997stochastic}. A stochastic process $(X_t)_{t \geq 0}$ is said to be solution to
	
	\begin{equation*}
	\circ \d X_t = \SN A_i(X_t) \circ \d W^i_t + A_0(X_t) \d t, \quad X_0 = x
	\end{equation*}
	if for all $f \in \C^2(M)$,

	\begin{align*}
		f(X_t)&= f(x) + \SN \int_0^t A_i f (X_s) \circ \d W^i_s + \int_0^t A_0f(X_s) \d s \\
		& = f(x) + \SN \int_0^t A_i f (X_s) \d W^i_s + \int_0^t A_0f(X_s) \d s + \frac{1}{2} \SN \int_0^t A_i^2 f(X_s) \d s . 
	\end{align*}
	The following result can be proved by using the definition of the Hessian.

	\begin{prop} \label{Conversion Ito strato dimension finie}
    The Stratonovich's SDE

$$
\circ \d X_t = \SN A_i(X_t) \circ \d W^i_t + A_0(X_t) \d t, \quad X_0 = x \in M,
$$
is equivalent to the following Itô's one:

$$
\d^\nabla X_t = \SN A_i(X_t) \d W^i_t + \left( A_0(X_t) + \frac{1}{2} \SN \nabla_{A_i} A_i(X_t) \right) \d t, \quad X_0 = x \in M.
$$
\end{prop}

	\begin{propdef}[\cite{Ito75}] \label{Prop def transport parallele stochastique M}
        Let $M$ be a Riemannian manifold with $\nabla$ the associated Levi-Civita connection, $A_0,A_1, \dots, A_N \in \Gamma(TM)$. Let $(X_t)_{t \geq 0}$ be the unique solution to:
        
            $$
    \circ \d X_t = \SN A_i(X_t) \circ \d W^i_t + A_0(X_t) \d t, \quad X_0 \in M.
    $$
    There exists a unique solution to the following stochastic covariant equation
    
    \begin{equation} \label{Transport para sto}
        \nabla_{\circ \d X_t}A_t = 0, \quad A_0 \in T_{X_0}M,
    \end{equation}
    in the following sense:
    
    $$
    \circ \d \langle A_t, B(X_t) \rangle_{X_t} = \SN \langle A_t, \nabla_{A_i(X_t)} B(X_t) \rangle_{X_t} \circ \d W^i_t + \langle A_t, \nabla_{A_0(X_t)}B(X_t) \rangle_{X_t} \d t,
    $$
    for all $B \in \Gamma(TM)$. This unique solution is called the \textit{stochastic parallel transport} along $(X_t)_{t \geq 0}$ with initial condition $A_0 \in T_{X_0}M$.
    \end{propdef}
	This is a well known fact that \eqref{Transport para sto} defines a stochastic process $(\tau_{0,t})_t$ in $O^x(M)$, where $O^x(M)$ denotes the orthonormal frame bundle, namely,
	
	$$
	O^x(M) = \sqcup_{y \in M} O^x_y(M),
	$$
	where $O^x_y(M)$ is the vector space of Euclidean isometries from $T_xM$ to $T_yM$. In particular, $\tau_{0,t} : T_xM \rightarrow T_{X_t}M$, where for any $A \in T_xM$, $(\tau_{0,t}(A))_{t}= (A_t)_t$ in \eqref{Transport para sto}. Let $J_0, J_1,\dots, J_N \in \mathrm{Hom}_{\mathrm{id}}(TM,TM)$ for $i=0,1,\dots,N$ be smooth bundle morphisms. An adapted $TM$-valued stochastic process $(U_t)_{t \geq 0}$ along $(X_t)_{t \geq 0}$ is said to be solution to the following covariant SDE:
	
	$$
	    \circ D_t U_t = \sum_{i=1}^{N} J_i(U_t) \circ \d W^i_t + J_0(U_t) \d t, \quad U_0 = U,
	$$
	if $(\tau_{0,t}^{-1} U_t)_{t \geq 0}$ is solution to the following SDE on $T_xM$:
	
	$$
	\circ \d \tau_{0,t}^{-1} U_t = \sum_{i=1}^{N} \tau_{0,t}^{-1} J_i(U_t) \circ \d W^i_t + \tau_{0,t}^{-1} J_0 (U_t) \d t, \quad U_0 = U.
	$$
 	In the same way, we say that $(U_t)_{t \geq 0}$ is solution to

	\begin{equation*}
	    D_t^\nabla U_t = \sum_{i=1}^{N} J_i(U_t) \d W^i_t + J_0(U_t) \d t, \quad U_0 = U,
	\end{equation*}
	if $(\tau_{0,t}^{-1}(U_t))_{t \geq 0}$ is solution to

	$$
	\d \tau_{0,t}^{-1} U_t = \sum_{i=1}^{N} \tau_{0,t}^{-1} J_i(U_t) \d W^i_t + \tau_{0,t}^{-1} J_0 (U_t) \d t, \quad U_0 = U.
	$$

 \subsection{Presentation of the finite-dimensional results} \label{S 1.3}
This subsection is dedicated to a presentation of, in a first time, well-known finite-dimensional results and, in a second time, to an adaptation of the results of Cattiaux and Guillin in \cite{cattiaux:hal-00805299} to the finite-dimensional Riemannian manifold case. Throughout this subsection we will consider a complete Riemannian manifold $(M, \langle \bullet, \bullet \rangle)$ and a $M$-valued solution flow $(X_t(\cdot))_t$ to an SDE. The associated Markov semigroup is denoted by $\mathrm{P}_t$ and, for a measure $\mu \in \P(M)$, the probability measure $\mu \mathrm{P}_t \in \P(M)$ is defined to be:

$$
\mu \mathrm{P}_t (\cdot) = \int_M \mathbb{P}(X_t(x) \in \cdot) \d \mu (x). 
$$
Let $f \in \mathscr{B}(M)$, we also define the measurable map $\mathrm{P}_t f(\cdot) = \mathbb{E}[ f(X_t(\cdot))] \in \mathscr{B}(M)$.

\vspace{10pt}

\noindent \textbf{Langevin dynamics on Riemannian manifolds.} A Langevin diffusion on a complete Riemannian manifold $(M, \langle \bullet, \bullet \rangle)$ is a diffusion with generator $L = \frac{1}{2} ( \Delta - \nabla \Psi)$ with $\Psi \in \C^\infty(M)$ a function such that $I(\Psi)=\int_M e^{-\Psi(x)} \d \vol(x)< +\infty$. It is well known that the process $(X_t)_{t \geq 0}$ admits a unique reversible measure $\bar{\mu}$ given by 
 
 $$
 \d \bar{\mu} = \frac{e^{-\Psi}}{I(\Psi)} \d \vol
 $$
 
 The Bakry-Emery $\mathrm{CD}(K,\infty)$ condition is the following
 
 \begin{equation} \tag{$\mathrm{CD}(K,\infty)$}
 \mathrm{Ric}(U,U) + \he \Psi(U,U) \geq K |U|^2, \quad \forall U \in TM,
 \end{equation}
 where $\mathrm{Ric}$ denotes the Ricci tensor and $\he$ denote the Hessian. Under this condition, it is well known that $(X_t)_{t \geq 0}$ satisfies
 
\begin{equation} \label{Inégalité gradient forte}
	|\nabla \mathrm{P}_t f| \leq e^{-Kt/2} \mathrm{P}_t (|\nabla f|), \quad \forall f \in \C^\infty_b(M),
\end{equation}
a local Poincaré inequality,

\begin{equation}
\mathrm{P}_t (f^2) - (\mathrm{P}_t f)^2 \leq 2 C(K,t) \mathrm{P}_t ( |\nabla f|^2),
\end{equation}
and a local Log-Sobolev inequality:

\begin{equation}
\mathrm{P}_t \left( f^2 \log f^2 \right) - \mathrm{P}_t (f^2) \log \mathrm{P}_t(f^2) \leq 4 C(K,t) \mathrm{P}_t ( |\nabla f|^2),
\end{equation}
where

\begin{equation} \label{ckt}
C(K,t) = \begin{cases}
	\frac{1-e^{-Kt}}{K} \text{ if } K \neq 0, \\
	t \text{ if } K =0.
\end{cases}
\end{equation}

\begin{figure}[t]
     \centering
     \includegraphics[width=0.7\textwidth]{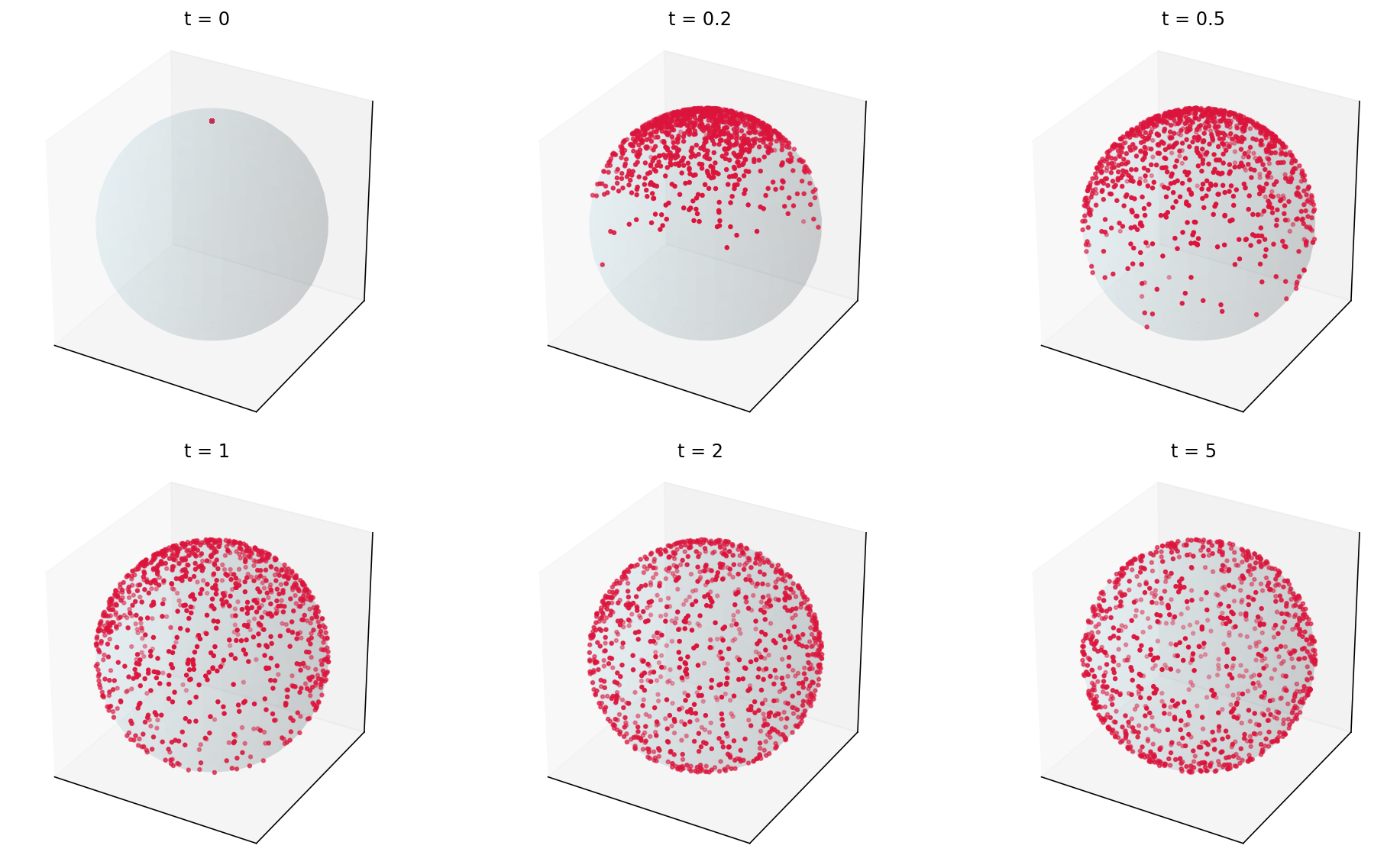}
     \caption{Illustration of the convergence of the law of the Brownian motion on the sphere $\mathbb{S}^2$ toward the volume measure.}
 \end{figure}
\vspace{10pt}
 
\noindent \textbf{Semi log-concave diffusions on Riemannian manifolds.} In \cite{cattiaux:hal-00805299}, the authors investigate the stability properties of $\mathbb{R}^d$-valued stochastic differential equations of the form
\begin{equation} \label{Cat Guil}
\d X_t = \SN A_i(X_t)  \d W^i_t + A_0(X_t)  \d t, \quad X_0 \in \mathbb{R}^d,
\end{equation}
where $\{W^i\}_{i=1}^N$ are $N$ independant Bronwnian motions, by means of the following condition:

\begin{equation} \label{E HCK condition cas euclidien} \tag{H.C.K}
\SN |A_i(x) - A_i(y)|^2 + 2 \langle A_0(x) - A_0(y), x - y \rangle \leq -K |x - y|^2, \quad \forall x,y \in \mathbb{R}^d,
\end{equation}
for some constant $K \in \mathbb{R}$. Note that, by removing the first term in the left hand side of \eqref{E HCK condition cas euclidien}, we recover the following inequality:

\begin{equation} \label{fortement dissipatif}
	2 \langle A_0(x) - A_0(y), x - y \rangle \leq -K |x - y|^2, \quad \forall x,y \in \mathbb{R}^d,
\end{equation}
In certain frameworks, when assumption \eqref{fortement dissipatif} is in force, the drift $b$ is said to be $K$-dissipative, which is equivalent to the $K$-concavity of the drift term $b$.

\begin{rem}
Let $x \in \mathbb{R}^d$ and $x_t = x +t v$ in \eqref{E HCK condition cas euclidien}. Let also

$$
f(t) = \SN |A_i(x) - A_i(x_t)|^2 + 2 \langle A_0(x) - A_0(x_t), x - x_t \rangle.
$$
Since $f(0)=0$ and $f'(0)=0$, by computing $\lim_{h\to 0} f(t)/t^2=f''(0)/2$, we obtain the infinitesimal version of \eqref{E HCK condition cas euclidien}, namely:

\begin{equation} \label{cat guil infi}
	\SN |D_x A_i(v)|^2 + 2 \langle D_x A_0(v), v \rangle \leq -K |v|^2.
\end{equation} 
Note that, in general, the infinitesimal formulation does not imply the global one.

\end{rem} Under suitable assumptions on the driving vector fields, the authors establish stability properties for the above SDE, see \cite[Theorem 1.8]{cattiaux:hal-00805299}. 
	
 	The purpose of this subsection is to extend the results of Cattiaux and Guillin in \cite{cattiaux:hal-00805299} to the framework of SDEs evolving on Riemannian manifolds. We shall see that the only difference arises from the necessity to account for the curvature of the manifold (see \eqref{E HCK sur M} below). We still consider $(M, \langle \bullet, \bullet \rangle)$ a complete Riemannian manifold. Let $A_0,A_1, \dots,A_N \in \Gamma(TM)$ and $(X_t)$ be the unique solution to
	
	\begin{equation*}
\d^\nabla X_t(x) = \SN A_i(X_t(x)) \d W^i_t + A_0(X_t(x)) \d t, \quad X_0(x) = x.
	\end{equation*}
	We will often consider its Stratonovich's formulation, namely,
	
\begin{equation} \label{EDS Strato sur M}
\circ \d X_t(x) = \SN A_i(X_t(x)) \circ \d W^i_t + A_0'(X_t(x)) \d t, \quad X_0(x) = x,
\end{equation}
	with $A_0' = A_0 - \frac{1}{2}\SN \nabla_{A_i} A_i$. Recall that, in this setting, $(X_t)_t$ admits a $\C^{0,\infty}(\mathbb{R}_+ \times M)$ and, if $f \in \C^k(M)$, $\mathrm{P}_t f \in \C^k(M)$ too, see \cite{Elworthy1982} or \cite{kunita1997stochastic}. A basic idea in order to study the long time behaviour of such a diffusion is to study the following distance for $x,y \in M$:

	\begin{equation} \label{E 4}
	\mathbb{E}\left[ d_M^2(X_t(x),X_t(y)) \right].
	\end{equation}
	Note that, by completeness, for any $x,y \in M$, there exists a (not necessarily unique) minimizing geodesic $(\gamma_a)_{0 \leq a \leq 1}$ joining $x$ to $y$ such that
	
	$$
	d_M^2(x,y) = \left(\int_0^1 |\partial_a \gamma_a|_{\gamma_a} \d a \right)^2 = \int_0^1 |\partial_a \gamma_a|_{\gamma_a}^2 \d a,
	$$
	where the geodesic can be written as $ \gamma_a = \exp_x(aV)$ for some $V \in T_xM$ and the last equality follows from the fact that the speed of a geodesic is constant. 
	Since we are interested in \eqref{E 4}, we need to obtain informations about the so-called tangent process $(TX_t(\partial_a \gamma_a))_{t \geq 0} = (\partial_a X_t(\gamma_a))_{t \geq 0}$. In the following, we will always denote by $R$ the Riemann curvature tensor.
	
	\begin{prop}[Itô's SDE of the Tangent process] \label{P Forme ito processus tangent dim finie}
		The tangent process $(TX_t)_{t \geq 0}$ is solution to the following Itô's covariant SDE:
		
		\begin{align*}
		D_t^\nabla TX_t(V) &= \SN \nabla_{TX_t(V)} A_i(X_t) \d W^i_t + \nabla_{TX_t(V)} A_0(X_t) \d t \\
		& - \frac{1}{2} \SN R(TX_t(V), A_i(X_t)) A_i(X_t), \quad TX_0(V) =  V \in T_{X_0} M.
		\end{align*}
	\end{prop}
	\begin{proof}
		There are many ways to prove this result, a proof is provided in Appendix \ref{Appendix A}, one can also prove it as in \cite[Chapter VII, Example 13D (i)]{Elworthy1982}.
	\end{proof}
	The extension of \eqref{cat guil infi} to the case of Riemannian manifolds is the following:
	
	\begin{equation} \tag{HC(\textit{K})} \label{E HCK sur M}
	\SN \left( \langle R(U, A_i)A_i, U \rangle - \langle \nabla_U A_i, \nabla_U A_i \rangle \right) - 2 \langle \nabla_U A_0, U \rangle \geq K |U|^2, \quad \forall U \in TM.
	\end{equation}

	\begin{rem}
		In the Riemannian case, it is more convenient to use such an infinitesimal formulation. Note that in that case, the global counterpart of \eqref{E HCK sur M} writes:

	$$
	\SN |P_{y,x}A_i(y)-A_i(x)|^2_x + \langle P_{y,x}A_0(y)-A_0(x), \overset{\to}{xy} \rangle + \SN \langle P_{y,x}\mathrm{Jac}_i(\overset{\to}{xy}) - A_i(x) , A_i(x) \rangle \leq -K |\overset{\to}{xy}|^2,
	$$
	where $\overset{\to}{xy} \in T_xM$ is a vector such that $\exp_x(\overset{\to}{xy}) =y$ and $\mathrm{Jac}_i(\overset{\to}{xy})=J_i(1)$ where $J_i$ is the unique Jacobi field along $(\gamma_t)_{0 \leq t \leq 1} = (\exp_x(t \overset{\to}{xy}))_{0 \leq t \leq 1}$ solution to:
		$$
\ddot{J}_i(t) = -R(J_i(t),\dot{\gamma}_t)\dot{\gamma}_t, \quad J_i(0) = A_i(x), \, \dot{J}_i(0)=0.
	$$
	As in the Euclidean case, the global formulation implies the infinitesimal one.
	\end{rem} 

	In particular, we can bound \eqref{E 4} as follows:
	
	\begin{align}
	\mathbb{E}[d_M^2(X_t(x), X_t(y))] & \leq \mathbb{E} \left[ \left( \int_0^1 |TX_t(\gamma_a)| \d a \right)^2 \right] \nonumber \\
	& \leq \mathbb{E} \left[ \int_0^1 |TX_t(\partial_a \gamma_a)|^2 \d a \right] = \int_0^1 \mathbb{E} \left[ |TX_t(\partial_a \gamma_a) |^2 \right] \d a. \label{E10}
	\end{align}
	Let 

	$$
	\P_2(M) := \left\{ \mu \in \P(M), \, \int_M d^2(x,x_0) \d \mu(x) < \infty \right\},
	$$
	where $x_0$ is a reference point. Note that the definition does not depend on the choice of $x_0$. We endow this space with the Wasserstein distance $W_2$:

\begin{equation*}
    W_2^2(\mu, \nu) = \inf \left \{ \int_{M \times M} d_M^2(x,y) \d\pi(x,y) \, : \, \pi \in \mathcal{C}(\mu,\nu) \right \},
\end{equation*}
where $\mathcal{C}(\mu,\nu)$ denotes the set of \textit{transport plans} between $\mu$ and $\nu$. It is straightforward to show that

$$
W_2^2(\delta_x, \delta_y) = d^2(x,y).
$$
	
	\begin{prop} \label{proposition décroissance exponentielle}
		Assume that the $\mathrm{HC}(K)$ condition holds. Then,
		 
		\begin{equation*}
		W_2^2(\delta_x \mathrm{P}_t, \delta_y \mathrm{P}_t) \leq \mathbb{E}[d^2(X_t(x),X_t(y))] \leq e^{-Kt} d^2_M(x,y), \quad \forall x,y \in M.
		\end{equation*}
		Consequently, for any $\vartheta_1,\vartheta_2 \in \P_2(M)$,

		$$
		W_2^2(\vartheta_1 \mathrm{P}_t, \vartheta_2 \mathrm{P}_t) \leq e^{-Kt} W_2^2(\vartheta_1,\vartheta_2).
		$$
	\end{prop}
	
	\begin{proof}
	Let $x,y \in M$ and $(\gamma_a)_{0 \leq a \leq 1}$ a minimizing geodesic from $x$ to $y$. Let also $V^a = \dot{\gamma}_a$ and $V_t^a =TX_t(\partial_a \gamma_a)$. We have by definition of the Wasserstein distance:
	
	$$
	W_2^2(\delta_x \mathrm{P}_t, \delta_y \mathrm{P}_t) \leq \mathbb{E} \left[ d_M^2(X_t(x), X_t(y)) \right].
	$$
	By \eqref{E10} we have
	
	$$
	\mathbb{E} \left[ d_M^2(X_t(x), X_t(y)) \right] \leq \int_0^1 \mathbb{E} \left[ |V_t^a |^2 \right] \d a,
	$$
	and, by applying Proposition \ref{P Forme ito processus tangent dim finie} and the $\mathrm{HC}(K)$ condition, we obtain:
	
	\begin{align}
	\d \mathbb{E}[|V_t^a|^2] &= \mathbb{E} \left[2 \langle \nabla_{V_t^a} A_0, V_t^a \rangle - \SN \langle R(V_t^a, A_i(X_t^a))A_i(X_t^a), V_t^a \rangle + |\nabla_{V_t^a}A_i(X_t^a)|^2 \right] \d t \label{E12} \\
	& \leq -K \mathbb{E}[|V_t^a|^2] \d t. \nonumber
	\end{align}
	Then, by applying Grönwall's Lemma, we obtain
	
	$$
	\mathbb{E} \left[ |V_t^a|^2 \right] \leq e^{-Kt} \mathbb{E}[ |V^a|^2].
	$$
	Thus, we obtained,
	
	$$
	\mathbb{E} \left[ d_M^2(X_t(x), X_t(y)) \right] \leq e^{-Kt} \int_0^1 \mathbb{E}[|V^a|^2] \d a = e^{-Kt} \mathbb{E} \left[ \int_0^1 |V^a|^2 \d a \right] = e^{-Kt} \mathbb{E}[ d_M^2(x,y)].
	$$
	To prove the second point, let $\mathcal{Z}_t[(x,y)] = [X_t(x),X_t(y)]$ be the synchronous coupling and denote by $\mathcal{Q}_t$ the associated Markov semigroup. Let $\pi \in \C(\vartheta_1,\vartheta_2)$ be an optimal coupling between $\vartheta_1$ and $\vartheta_2$. Then,

	\begin{align*}
		W_2^2(\vartheta_1 \mathrm{P}_t, \vartheta_2 \mathrm{P}_t) & \leq \int_{M^2} d^2(x,y) \d (\pi \mathcal{Q}_t)(x,y) \\
		& = \int_{M^2} \mathbb{E} \left[d^2(X_t(x),X_t(y)) \right] \d \pi(x,y) \\
		& \leq e^{-Kt} \int_{M^2} d^2(x,y) \d \pi(x,y) = e^{-Kt} W_2^2(\vartheta_1,\vartheta_2). 
	\end{align*}
	The result is proved.
	\end{proof}

	    \begin{coro} \label{C Mesure invariante dim finie}
    Assume that the $\mathrm{HC}(K)$ condition holds for $K > 0$. Then, there exists a unique probability law $\bar{\mu} \in \P_2(M)$ such that
    \begin{equation} \label{E 1}
		\lim_{t \rightarrow \infty } W_2^2(\vartheta \mathrm{P}_t, \bar{\mu}) =0, \quad \forall \vartheta \in \P_2(M).
	\end{equation}
    Moreover, $\bar{\mu}$ is the unique invariant measure for the semigroup $(\mathrm{P}_t)_{t \geq 0}$.
    \end{coro}
    
    \begin{proof}
    		The existence property follows from the fact that $\mathrm{P}_{t+s} = \mathrm{P}_t \mathrm{P}_s$ for $t,s \geq 0$. Indeed, by applying Theorem \ref{proposition décroissance exponentielle}, we obtain:
    		
    		$$
    		W_2^2(\vartheta\mathrm{P}_{t+s},\vartheta \mathrm{P}_t) \leq e^{-Kt} W_2^2(\vartheta \mathrm{P}_s, \vartheta).
    		$$
    	We conclude by using the fact that $\P_2(M)$ is a complete space. This limit does not depend on $\vartheta$ since for $\vartheta' \in \P(M)$ and by Proposition \ref{proposition décroissance exponentielle},
		
		$$
		W_2(\vartheta'\mathrm{P}_t, \bar{\mu}) \leq W_2(\vartheta'\mathrm{P}_t, \vartheta \mathrm{P}_t) + W_2(\vartheta \mathrm{P}_t, \bar{\mu}) \underset{t \to \infty}{\to} 0.
		$$
		Hence, \eqref{E 1} holds.
		The fact that $\bar{\mu}$ is an invariant measure comes from:
    	
    	$$
    	\bar{\mu} = \lim_{s \to \infty} \bar{\mu} \mathrm{P}_{t +s} = \lim_{s \to \infty}  (\bar{\mu} \mathrm{P}_s) \mathrm{P}_t = \bar{\mu}\mathrm{P}_t .
    	$$
    	In order to see that $\bar{\mu}$ is the unique invariant measure, suppose that there exists $\bar{\mu}'$ another invariant measure. By \eqref{E 1},

		$$
		W_2^2(\bar{\mu}',\bar{\mu}) =W_2^2(\bar{\mu}' \mathrm{P}_t, \bar{\mu}) \underset{t \to \infty}{\to} 0, 
		$$
		then, $\bar{\mu}' = \bar{\mu}$.
    \end{proof}
    
    \begin{prop}
    Assume that the $\mathrm{HC}(K)$ condition holds and let $f \in \C^1(M)$. Then, we have the following sub-commutativity property:
		
		$$
		|\nabla \mathrm{P}_t f(x)|_{x}^2 \leq e^{-Kt} \mathrm{P}_t \left( |\nabla f(x)|_x^2 \right), \quad \forall x \in M.
		$$
    \end{prop}
    
    \begin{proof}
    	Let $x \in M$, $V \in T_xM$ such that $|V|_x = 1$ and $\gamma_a = \exp(aV)$. Let $V^r = \dot{\gamma}_r$ and $V_t^r = TX_t(V^r)$. We have
    	
    	\begin{align}
    	\frac{1}{a} |\mathrm{P}_t f(x) - \mathrm{P}_t f(x_a) | & \leq \frac{1}{a} \int_0^a \mathbb{E} \left[ \langle \nabla f(X_t^r), V_t^r \rangle \right] \d r \nonumber \\
    	& \leq e^{-Kt/2} \sup_{0 \leq s \leq a} \mathbb{E}\left[ |\nabla f(X_t^s)|^2 \right]^{1/2} \frac{1}{a} \int_0^a \mathbb{E} [|V^r|^2]^{1/2} \d r, \label{E13}
    	\end{align}
    	where \eqref{E13} comes from Cauchy-Schwarz inequality and inequality \eqref{E12}. Thus, by taking the limit when $a \to 0$ and since $|V|=1$, we obtain
    	
    	$$
	\langle \nabla \mathrm{P}_t f(x), V \rangle_x \leq e^{-Kt} \mathrm{P}_t(|\nabla f(x)|^2)^{1/2}.
    	$$
    	The proof is done.
    \end{proof}
    
    \begin{rem}
    Note that this inequality is weaker than \eqref{Inégalité gradient forte} since, by Cauchy-Schwarz inequality,
    
    $$
    |\nabla \mathrm{P}_t f|^2 \leq e^{-Kt} \left( \mathrm{P}_t (|\nabla f|) \right)^2 \leq e^{-Kt} \mathrm{P}_t(|\nabla f|^2).
    $$
    \end{rem}
    
    \begin{prop}[Local Poincaré Inequality]
    		Assume that the $\mathrm{HC}(K)$ condition holds for some $K \in \mathbb{R}$ and let $f \in \C^2(M)$. Then, the following local Poincaré inequality is satisfied:
		
		\begin{equation*}
		\mathrm{P}_t (f^2) - (\mathrm{P}_t f)^2 \leq 2|A|^2_{2,\infty}  C(K,t) \mathrm{P}_t ( | \nabla f |^2), \quad \forall t \geq 0, 
		\end{equation*}
		where $|A|^2_{2,\infty} = \sup_{x \in M} \SN |A_i(x)|_x^2$ and $C(K,t)$ is defined in \eqref{ckt}.
    \end{prop}
    
    \begin{proof}
    The proof is almost the same as in \cite[Proposition 5.4.1]{logsob2000}. We just use the fact that 
    
    $$
    \Gamma(f) \leq |A|^2_{2,\infty} |\nabla f|^2,
    $$
	where $\Gamma$ denotes the associated Carré du Champ,

	$$
	\Gamma(f) = \SN (A_i f)^2.
	$$
    \end{proof}

	\noindent \textbf{The Brownian case.} Let $M$ be a closed Riemannian manifold, $\{E_i\}_{i=1}^N$ be the canonical basis of $\mathbb{R}^N$, $\phi =(\phi_1,\dots,\phi_N) : M \hookrightarrow \mathbb{R}^N$ be an isometric embedding. Let $\II : (A,B) \in\Gamma(TM) \times \Gamma(TM) \rightarrow \left(D_A B \right)^\perp \in  \Gamma(TM^\perp)$ be the so-called second fundamental form. The assumption of the isometric embedding is here to give a geometric meaning to the second fundamental form, but an isometric immersion is sufficient. We still work on $M$, $\nabla$ still denotes the Levi-Civita connection and the gradient. This is well known, see \cite[Chapter 8]{Lee2018} for instance, that $\II$ is actually a symmetric tensor and satisfies the so-called Gauss formula, namely,
		
		$$
		\II(A,B) = D_A T \phi(B) - T \phi( \nabla_A B) = \he \phi (A,B),
		$$
		where $T\phi : TM \rightarrow \mathbb{R}^N$ denotes the differential of $\phi$ and $D$ denotes the usual differential on $\mathbb{R}^N$. Let $(X_t)_{t \geq 0}$ be the unique solution to 

		$$
		\d^\nabla X_t(x) = \SN \nabla \phi_i(X_t(x)) \d W^i_t, \quad X_0(x)=x. 
		$$
		It is well known that $(X_t)_{t \geq 0}$ is a Brownian motion on $M$, see \cite[Proposition 7C]{Elw88} for instance. Note that, for all $x \in M$, $\{\nabla \phi_i(x)\}_{i=1}^N$ satisfies 

		$$
		\langle A,B \rangle = \SN \langle U, \nabla \phi_i(x)\rangle \langle V, \nabla \phi_i(x) \rangle = \SN \langle U, E_i \rangle \langle V,E_i \rangle, \quad U,V \in T_xM,
		$$
		and

		\begin{align*}
			\mathrm{H} \phi_i (U,V) & = AB \langle \phi, E_i \rangle - \nabla_A B \langle \phi, E_i \rangle \\
			& = A \langle T\phi(B), E_i \rangle - \langle T\phi(\nabla_A B), E_i \rangle \\
			& = \langle D_A T\phi(B) - T\phi(\nabla_A B), E_i \rangle \\
			&= \langle \II(A,B), E_i \rangle.
		\end{align*}
		Let $x \in M$ and $\{e_j\}_{j=1}^n$ be an orthonormal frame of $T_xM$. Then, on the one hand,

		$$
		\SN \left\langle R(U,\nabla \phi_i(x))\nabla \phi_i(x), U \right\rangle_x = \mathrm{Ric}(U,U), \quad x \in M,\, U \in T_xM.
		$$
		On the other hand, since 

		$$
		|\nabla_U \nabla \phi_i|^2 = \sum_{j=1}^n \langle \nabla_U \nabla \phi_i, e_j \rangle^2 = \sum_{j=1}^n \mathrm{H}\phi_i(U,e_j)^2,
		$$
		it holds:

		\begin{align*}
			\SN |\nabla_U \nabla \phi_i|^2 = \SN \sum_{j=1}^n \left( \mathrm{H}\phi_i(U,e_j) \right)^2 = \sum_{j=1}^n \SN \left\langle \II(U,e_j), E_i \right\rangle^2 = \sum_{j=1}^n |\II(U,e_j)|^2.
		\end{align*}
		Thus, in this case, the \eqref{E HCK sur M} condition reduces to

		$$
		\mathrm{Ric}(U,U) - \sum_{j=1}^n |\II(U,e_j)|^2 \geq K|U|^2, \quad U \in TM.
		$$

		\begin{rem}
			Alternatively, one may use the embedding and carry out all computations in $\mathbb{R}^d$. Combining the results of Appendix \ref{appendice A1} with \eqref{E HCK condition cas euclidien} yields the same result.
		\end{rem}

		\begin{exemple}[Hypersurfaces] \label{Exemple hypersurfaces}
		Let $\phi : M \hookrightarrow \mathbb{R}^{N}$ be a complete hypersurface with $\mathrm{dim}(M)=n$, $N=n+1$. The reader is refered to \cite[Chapter 8]{Lee2018} for a good introduction to Riemannian submanifolds, especially for hypersurfaces. Let us denote by $h : TM \times TM \rightarrow \mathbb{R}$ the symmetric tensor such that $\II(U,V)= h(U,V) \overset{\to}{n}$, where $\overset{\to}{n}$ denote the outward unit normal vector. Let also $\mathfrak{s} : TM \rightarrow TM$ be the unique self-adjoint operator such that $\langle \mathfrak{s}U, V \rangle = h(U,V)$. By the spectral theorem, there exists an orthonormal frame $\{e_i\}_{i=1}^n$ such that for all $i=1, \dots,n$, there exists $\kappa_i \in \mathbb{R}$ such that $\mathfrak{s} e_i = \kappa_i e_i$ (the constant $\kappa_i$ is called the $i$-th principal curvature).

		By Gauss formula (see \cite[Theorem 8.13]{Lee2018}), it follows:

		\begin{align*}
			\mathrm{Ric}(U,U)= \sum_{i=1}^n \langle R(U,e_i)e_i, U \rangle &= \sum_{i=1}^n h(U,U)h(e_i,e_i) - h(U,e_i)^2 \\
			& = \sum_{i=1}^n \langle \II(U,U), \overset{\to}{n}\rangle \langle \II(e_i,e_i), \overset{\to}{n}\rangle - \langle \II(U,e_i), \overset{\to}{n}\rangle^2 \\
			& = \sum_{i=1}^n \langle \mathfrak{s} U, U \rangle \langle \mathfrak{s}e_i, e_i \rangle - \langle \mathfrak{s} U, e_i \rangle^2  = \langle \mathfrak{s} U, U \rangle \mathrm{Tr}(\mathfrak{s}) - |\mathfrak{s}U|^2.
		\end{align*}
		Consequently, the \eqref{E HCK sur M} condition reduces to

		$$
		\langle \mathfrak{s} U, U \rangle \mathrm{Tr}(\mathfrak{s}) - 2|\mathfrak{s}U|^2 \geq K|U|^2, \quad U \in TM.
		$$
		Writing $U = \sum_{i=1}^n U_i e_i$, we obtain 

		\begin{equation}\label{Trace s}
			\langle \mathfrak{s} U, U \rangle \mathrm{Tr}(\mathfrak{s}) - 2s|\mathfrak{s}U|^2 = \sum_{j=1}^n U_j^2 \kappa_j\cdot \left(\sum_{i=1}^n \kappa_i \right) - 2 U_j^2 \kappa_j^2 = \sum_{j=1}^n \kappa_j \left( \sum_{i=1}^n \kappa_i -2 \kappa_j  \right) U_j^2.
		\end{equation}
		Therefore, \eqref{E HCK sur M} holds for $K = \min_{1 \leq j \leq n }\kappa_j \left(\sum_{i=1}^n \kappa_i - 2 \kappa_j\right)>0$ if for all $j=1,\dots,n$, there exists a constant $c >0$ such that $\kappa_j \geq c >0$ and if 

		$$
		2 \kappa_j < \sum_{i=1}^n \kappa_i.
		$$
		Note that we can also write \eqref{Trace s} in the following way:

		\begin{align*}
			\sum_{j=1}^n \kappa_j \left( \sum_{i=1}^n \kappa_i -2 \kappa_j  \right) U_j^2 &= \sum_{j=1}^n 2\left( \frac{1}{2} \mathrm{Tr}(\mathfrak{s})\kappa_j - \kappa_j^2\right)U_j^2 \\
			& = \sum_{j=1}^n \left[\frac{\mathrm{Tr}(\mathfrak{s})^2}{8} - 2 \left(\kappa_j - \frac{\mathrm{Tr}(\mathfrak{s})}{4}\right)^2 \right] U_j^2 \geq \left[ \frac{\mathrm{Tr}(\mathfrak{s})^2}{8} - 2 \max_{1 \leq j \leq n} \left( \kappa_j - \frac{1}{4} \mathrm{Tr}(\mathfrak{s}) \right)^2 \right] |U|^2.
		\end{align*}
	\end{exemple}

	\begin{exemple}[Spheres]
		Let $M= \mathbb{S}^n$ be the unit sphere, $N=n+1$ and consider the usual embedding $\phi : \mathbb{S}^n \hookrightarrow \mathbb{R}^{N}$. In this case, $\mathrm{Ric}(U,U) = (n-1) |U|^2$ and $\mathfrak{s}= -\mathrm{Id}$, thus, $\II(U,V) = -\langle U,V\rangle \overset{\to}{n}$, where $\overset{\to}{n}$ is the outward unit normal vector. Thus,
			
			$$
			\sum_{j=1}^n |\II(U,e_j)|^2 = \sum_{j=1}^n \langle U,e_j \rangle^2 = |U|^2.
			$$
			Consequently,

			$$
			\mathrm{Ric}(U,U) - \sum_{j=1}^n |\II(U,e_j)|^2 = (n-2)|U|^2.
			$$
			In particular, for $n \geq 3$, \eqref{E HCK sur M} holds for $K=n-2 >0$. Since $\kappa_j=1$ for all $j=1,\dots,n$, we could have find the same result by applying the reasoning in Example \ref{Exemple hypersurfaces}.

			If one consider $M = \mathbb{S}^{n_1} \times \dots \times \mathbb{S}^{n_k}$, with the embedding $M \hookrightarrow \mathbb{R}^{n_1+1} \times \dots \times \mathbb{R}^{n_k +1}$, one can easily compute 

			$$
			\mathrm{Ric}(U,U) - \sum_{j=1}^n |\II(U,e_j)|^2 = \sum_{l=1}^k (n_l-2)|U_l|^2, \quad U=(U_1,\dots,U_k) \in T_xM.
			$$
			Thus, if for all $l=1,\dots,k$, $n_l \geq 3$, \eqref{E HCK sur M} holds for $K =\min_{1 \leq l \leq k}n_l >0$.
	\end{exemple}

	\subsection{Structure of the paper}

	The paper is organized as follows. In Section \ref{S 2}, we introduce the infinite-dimensional geometric framework, namely, some basic definitions and properties about the Wasserstein space $\P$ and the Diffeomorphism group $\mathscr{D}$. In Section \ref{S 3}, we present some results of stochastic calculus on the Wasserstein space $\P$ and we show some regularity results about the solution flow. Section \ref{S 4} is devoted to reminders of Stochastic calculus on ILH vector bundles, to the study of the tangent process of the flow, $\T \M_t$ and to the proof of the main result, that is:

	\begin{thm}[Summary of Theorem \ref{Décroissance exp}, Corollary \ref{Coro mesu inv}, Proposition \ref{Sub-Commutativity property} and Proposition \ref{Local Poincare infini}]
		Let $\M_t$ be the solution flow of \eqref{EDS principale Ito}, denote by $\se_t$ the associated Markov semigroup and assume that \eqref{E IWCD} holds for some $K \in \mathbb{R}$. Then, for any $\vartheta_1,\vartheta_2 \in \P(\P)$, 

		$$
		\bm W_2^2(\vartheta_1 \se_t,\vartheta_2 \se_t) \leq e^{-Kt} \bm W_2^2(\vartheta_1,\vartheta_2),
		$$
		for any $\mu \in \P$ and $F \in \C^1(\P)$, it holds

		$$
		|\grad \se_t F(\mu)|_{\mu}^2 \leq e^{-Kt} \se_t \left( |\grad F(\mu)|_\mu^2 \right), \quad \forall \mu \in \P_\infty,
		$$
		and for any $F \in \C^2(\P)$, a local Poincaré inequality holds, see Proposition \ref{Local Poincare infini}. Moreover, if $K >0$, there exists a unique $\Lambda \in \P(\P)$ such that

    \begin{equation*}
		\lim_{t \rightarrow \infty } \W_2^2(\vartheta_1 \se_t, \Lambda) =0, \quad \forall \vartheta_1 \in \P(\P),
	\end{equation*}
    and $\Lambda$ is the unique invariant measure in $\P(\P)$ for the semigroup $(\se_t)_{t \geq 0}$.
	\end{thm}
	We also prove similar results on $M \times \P$ with a less sharp condition and we propose another condition for the study of lifts of Brownian motions on $M$ based on the parallel coupling on the base manifold.
	Finally, Section \ref{S 5} is devoted to the treatment of some examples, the first one being lifts of SDEs on $M$ and the second one being the Gaussian space. The reader can also find at the end of this paper a notation table.

\section{Reminders about the Wasserstein space} \label{S 2}

\subsection{Infinite-dimensional geometric setting} \label{S 2.1}
 In this subsection we present some definitions and results about the geometric structure of $\P$ (Otto formalism) and of $\mathscr{D}$ (with the ILH Lie group structure à la Ebin-Marsden). A detailed presentation of this material is done in \cite{gigli2011second}, \cite{lott2007geometriccalculationswassersteinspace}, \cite{Martin2025}, \cite{Otto31012001}. \vspace{10pt}

\noindent \textbf{The Wasserstein space.} In this article we consider a smooth closed Riemannian manifold $(M, \langle\bullet, \bullet \rangle)$. The volume measure $\vol$ is normalized, that is, $\vol(M)=1$. Let $\P$ denote the \textit{Wasserstein space} over $M$ i.e. the space of probability measures endowed with the Wasserstein distance $W_2$:

\begin{equation*}
    W_2^2(\mu, \nu) = \inf \left \{ \int_{M \times M} d_M^2(x,y) \d\pi(x,y) \, : \, \pi \in \mathcal{C}(\mu,\nu) \right \},
\end{equation*}
where $\mathcal{C}(\mu,\nu)$ denotes the set of \textit{transport plans} between $\mu$ and $\nu$. Moreover, let $\P_\infty \subset \P$ denote the \textit{smooth Wasserstein space} over $M$, namely:

$$
\P_\infty = \left\{ \mu \in \P \, : \, \frac{\d \mu}{\d \vol} = \rho \in \C^\infty(M),\, \rho >0 \right\}.
$$
equipped with the same distance and the space of purely atomic measures, namely,

$$
\P_{\mathrm{p.a.}} = \left\{ \mu \in \P : \, \mu = \sum_{i=1}^\infty \alpha_i \delta_{x_i}, \, \left\{\alpha_i\right\}_{i \in \mathbb{N}} \in \mathbb{R}_+^{\mathbb{N}}, \, \left\{x_i\right\}_{i \in \mathbb{N}} \in M^{\mathbb{N}} \right\}
$$

Additionally, let us mention the following fact which is proved using topological arguments in probability spaces (see for example \cite[Remark 6.19]{villani2009optimal}).

\begin{prop}
	If $M$ is a compact metric space, the associated Wasserstein space is compact too.
\end{prop}

The \textit{tangent space} of $\P$ in $\mu \in \P$ is:

$$
\T_\mu \P = \overline{ \left \{ \nabla \phi \, : \, \phi \in \C^\infty(M) \right \} }^{L^2(\mu)},
$$
and the \textit{tangent bundle} is $ \T \P = \underset{\mu \in \P}{\sqcup}\T_\mu \P$. Let us also introduce a dense subset of the latter,

$$
T_\mu \P = \left \{ \nabla \phi \, : \, \phi \in \C^\infty(M) \right \}, \quad \mu \in \P,
$$
which will sometime be denoted by $T_\mu \P_\infty$ when $\mu \in \P_\infty$. We denote by $T\P$ (or $T\P_\infty$) the associated tangent bundle. The \textit{Otto's metric} on $\T\P$ is defined to be

$$
\langle Z_1, Z_2 \rangle_\mu = \int_M \langle Z_1(x), Z_2(x) \rangle \d\mu,
$$
for all $Z_1, Z_2 \in \T_\mu \P$.

A \textit{vector field} on $\P$ is a map $\bar{Z} : \P \rightarrow \T \P$ such that $\bar{Z}(\mu) \in \T_\mu \P$. The notion of vector fields on $\P_\infty$ are defined analogously. 

\begin{exemple}
	Let $\phi \in \C^\infty(M)$, the associated constant vector field is defined to be $V_\phi : \mu \rightarrow V_\phi(\mu) = \nabla \phi(\cdot)$. For $\mu \in \P$ we will also denote by $V_\phi$ the tangent vector $V_\phi(\mu) = \nabla \phi \in T \P$.
	\end{exemple}

Let $\mu \in \P$, we define the \textit{Lie derivative} of a functional $F$ in the direction $V_\phi \in T_\mu \P$ to be

$$
\bar{\mathcal{L}}_{V_\phi} F(\mu) = \left. \frac{\d}{\d t} \right|_{t=0} F(\mu_t),
$$
where $(\mu_t)_{t \geq 0} =((U_t)_*\mu)_{t \geq 0}$ with

\begin{equation*}
    \partial_t U_t(x) = \nabla \phi(U_t(x)), \quad U_0(x) = x.
\end{equation*}
 In particular, $(\mu_t)_{t \geq 0}$ satisfies

\begin{equation} \label{Equation dérivée mesure}
\left. \frac{\d}{\d t} \right|_{t=0} \mu_t = V_\phi, \quad \mu_0 = \mu,
\end{equation}
in the sense of distribution, i.e. 

$$
\left. \frac{\d}{\d t} \right|_{t=0} \int_M f(x) \d\mu_t(x) = \int_M \langle \nabla f(x), \nabla \phi(x) \rangle \d\mu(x), \quad \forall f \in \C^\infty(M).
$$
We now introduce two notions of differentiability, as defined in \cite{RenWang2021}: the first is associated with the gradient flow, the second with the exponential flow. A functional $F$ is said to be \textit{intrinsically differentiable} in $\mu$ if the map 

\begin{align*}
    T_\mu \P &\longrightarrow \mathbb{R} \\
    V_\phi &\longmapsto \bar{\mathcal{L}}_{V_\phi} F(\mu)
\end{align*}
is a bounded linear functional. In this case, this map extends by density to a unique bounded linear functional:

\begin{align*}
    \T_\mu \P &\longrightarrow \mathbb{R} \\
    Z &\longmapsto \bar{\mathcal{L}}_{Z} F(\mu).
\end{align*}
The \textit{gradient} of an intrinsically differentiable functional in $\mu \in \P$ is the unique vector (tangent to $\P$) $\grad F(\mu) \in \T_\mu \P$ (the existence and uniqueness being guaranteed by Riesz representation theorem) such that for all $\bar{Z} \in \T_\mu \P$:

$$
\bar{\mathcal{L}}_{Z}F(\mu) = \langle \grad F(\mu), Z \rangle_\mu.
$$

Let $A \in L^2(M,TM;\mu)$ be a vector field and consider the measurable map $U^A := \exp(A) : x \in M \mapsto \exp_x(A(x)) \in M$. Due to \cite[Corollary 2.2]{RenWang2021}, a map $F : \P \rightarrow \mathbb{R}$ is said to be \textit{$L$-differentiable} in $\mu$ if it is intrinsically differentiable in $\mu$ with

$$
\lim_{|Z|_\mu \to 0} \frac{\left|F \left(\left(U^A \right)_*\mu\right) - F(\mu) - \bar{\mathcal{L}}_{A} F(\mu) \right|}{|A|_\mu}=0,
$$
When $F$ is intrinsically differentiable (resp. L-differentiable) in all $\mu \in \P$, we say that $F$ is intrinsically differentiable (resp. L-differentiable).
These definitions immediatly extend to functions with values in an Euclidean space $(E,\langle \bullet, \bullet\rangle_{E})$ of dimension $n'$. In particular, for $F : \P \rightarrow E$, the gradient, when it exists, may be seen as a map $\grad F : \P \times M \rightarrow E \otimes TM$ such that for all $(\mu,x) \in \P \times M$, $\grad F(\mu)(x) \in E \otimes T_xM$. Let $A \in \Gamma(E \otimes TM)$, $(e_i)_{i=1,\dots,d'}$ be a basis of $E$ and $A(x)= \sum_{i=1}^{n'} e_i \otimes A_i(x) \in E \otimes T_xM$ be the decomposition of $A$ according to this basis. We also denote by

\begin{align*}
	\langle \cdot, \cdot \rangle : (E \otimes TM) \times TM &\longrightarrow E \\
	(A,U) &\longmapsto \sum_{i=1}^{n'} e_i \langle A_i,U \rangle
\end{align*} 
and 

\begin{align*}
	\langle \cdot, \cdot \rangle : (E \otimes TM) \times (E \otimes TM) & \longrightarrow \mathbb{R} \\
	(A,v \otimes U) & \longmapsto \left \langle \left \langle A,U \right \rangle, v \right \rangle_{E}
\end{align*}

\begin{defi}
	Let $N$ be a complete Riemannian manifold, $E$ an Euclidean vector space $F : \P \rightarrow E$.
\begin{enumerate}
	\item We say that $F \in \C^{1}(\P,E)$ (resp. $\C^1_S(\P,E)$) whenever $F$ is $L$-differentiable and the map 
	
	$$
	(\mu,x) \in \P \times M \longmapsto \grad F(\mu)(x) \in E \otimes TM
	$$
	belong to $\C^0(M,\times \P,E)$ (resp. $\C^{\infty,0}(M \times \P,E))$.
\item We say that $F \in \C^{(1,1)}(\P,E)$ if $F \in \C^{1}(\P,E)$ and if for all $\mu \in \P$, the vector field $x \in M \mapsto \grad F(\mu)(x) \in E \otimes TM$ is differentiable such that

\begin{align*}
	\P \times M & \longrightarrow E \otimes TM \otimes TM \\
	(\mu,x) &\longmapsto \nabla \grad F(\mu)(x)
\end{align*}
is continuous, where $\nabla$ denotes the Levi-Civita connection (here we used the identification $TM \cong T^*M$).

\item We say that $F \in \C^{2}(\P,E)$ (resp. $\C^2_S(\P,E)$) if $F \in \C^{(1,1)}(\P,E)$ (resp. $F \in \C^1_S(\P,E)$) and $\mu \in \P \mapsto \grad F(\mu)(x) \in E \otimes T_xM$ is $L$-differentiable for all $x \in M$ such that

\begin{align*}
	\P(M) \times M \times M & \longrightarrow E \otimes TM \otimes TM \\
	(\mu,x,y) &\longmapsto \grad^2 F(\mu)(x,y) = \grad \left( \grad F(\mu)(x) \right)(y) \in E \otimes T_xM \otimes T_yM
\end{align*}
belong to $\C^0(M \times M \times \P, E \otimes TM \otimes TM)$ (resp. $\C^{\infty,\infty,0}(M \times M \times \P, E \otimes TM \otimes TM)$ ).

\item We say that $F \in \C^{2,2}(N \times \P,E)$ (resp. $\C^{\infty,2}_S(N \times \P,E)$) if $F \in \C^{0}(N \times \P)$, $F(y, \cdot) \in \C^{2}(\P,E)$ (resp. $\C^2_S(\P,E)$) for all $y \in N$, $F(\cdot,\mu) \in \C^{2}(N,E)$ (resp. $\C^\infty(N,E)$) for all $\mu \in \P$, $\mu \mapsto \nabla F(y,\mu)$ is $L$-differentiable for all $y \in N$ and the map

\begin{align*}
	\P(M) \times M \times N & \longrightarrow E \otimes TM \otimes TN \\
	(\mu,x,y) &\longmapsto \left(\grad \nabla_y F(\mu,y)\right)(x)
\end{align*}
belong to $\C^0(N \times M \times \P, E \otimes TM \otimes TN)$ (resp. $\C^{\infty,\infty,0}(N \times M \times \P, E \otimes TM \otimes TN)$). We further assume that $\grad F \in \C^{1,1,0}(N \times M \times \P,E \otimes TM)$ (resp. $\C^{\infty,\infty,0}(N \times M \times \P, E \otimes TM)$) and $\grad^2 F \in \C^0(N \times M^2 \times \P, E \otimes TM \otimes TM)$ (resp. $\C^{\infty,\infty,0}(N \times M^2 \times \P, E \otimes TM \otimes TM)$). 
\end{enumerate} 

The notions of functions in $\C^{k}(\P,E)$ (resp. $\C^{k}_S(\P,E)$) or in $\C^{k,l}(\P \times N,E)$ (resp. $\C^{k,l}_S(\P \times N,E)$) are defined analogously. 
\end{defi} 

\begin{defi}
	 A map $U: (a,x) \in [0,1] \times M \mapsto U_a(x) \in M$ is said to be a $\C^k$ path of diffeomorphisms if $U \in \C^{k,\infty}([0,1] \times M, M)$, and if for all $a \in [0,1], U_a$ is invertible with $U^{-1} \in \C^{k,\infty}([0,1] \times M,M)$. We say that $U$ is a $\C^k$ flow of diffeomorphisms if the flow property holds, $U_{a+b} = U_{a} \circ U_b$.
\end{defi}

\begin{lem} \label{Lemme fondamental L}
	Let $\mu \in \P$, $(U_a)_{0 \leq a \leq 1}$ be a path of measurable maps $U_a : M \to M$ with $U_0=\mathrm{id}$ and $U_a(x)=\exp_x(aA(x))$ for some $A \in L^2(M,TM;\mu)$, $\mu_a = (U_a)_* \mu$. Assume that for all $a \in [0,1]$, there exists $U_a^{*} : M \to M$ such that $U_a \circ U_a^* = \mathrm{id}$ $\mu_a$-almost everywhere and $U_a^* \circ U_a = \mathrm{id}$ $\mu$-almost everywhere. Then, $a \mapsto \partial_a U_a \circ U_a^* \in L^1\left([0,1],L^2(M,TM;\mu);\d a\right)$ and for any $F \in \C^1(\P)$, the map $a \mapsto F(\mu_a)$ is differentiable and

	\begin{equation} \label{E 38000}
		F(\mu_1)-F(\mu) = \int_0^1 \langle \grad F(\mu_a), \partial_a U_a \circ U_a^* \rangle_{\mu_a} \d a.
	\end{equation}
\end{lem}

\begin{proof}
	The first point follows from the fact that since $U_a(x)=\exp_x(aA(x))$,

	$$
	\int_0^1\int_M |\partial_a U_a\circ U_a^*(x)|^2_{x} \d \mu_a \d a = \int_0^1 \int_M |\partial_a U_a(x)|^2_{U_a(x)} \d \mu(x) \d a = \int_0^1 \int_M |A(x)|^2_{x} \d \mu(x) \d a = \int_M |A(x)|^2_{x} \d \mu(x).
	$$
	Moreover, by definition, $U_{a+\varepsilon}(x) = \psi_\varepsilon^a \circ U_a(x)$ where $\psi_\varepsilon^a(x) = \exp_x(\varepsilon \partial_a U_a \circ U_a^*(x))$. Thus, $\mu_{a+\varepsilon} = (\psi_\varepsilon^a)_*\mu_a$ and the differentiability of $a \mapsto F(\mu_a)$ follows from the definition of $L$-differentiability. Finally, \eqref{E 38000} follows from Cauchy-Schwarz and the fact that $\sup_{(x,\mu) \in M \times \P} \grad F(x,\mu)<\infty$.
\end{proof}
The assumptions of Lemma \ref{Lemme fondamental L} are satisfied, for instance, when $A \in \Gamma(TM)$ is such that the family of maps $(U_a)_{0 \leq a \leq 1} : x \in M \mapsto (\exp_x(aA(x)))_{0 \leq a \leq 1}$ defines a $C^\infty$ path of diffeomorphisms. Another important example is provided by the displacement interpolation between two absolutely continuous measures $\mu$ and $\nu$, see \cite{mccann2001polar} for a characterization of the corresponding transport maps. These two examples lead to the following two corollaries.

	\begin{coro} \label{Lemme fondamental L2}
		Let $(\varphi_a)_{a \geq 0}$ be a $\C^1$ path of diffeomorphisms, $\mu_a = (\varphi_a)_*\mu$ for $\mu \in \P$. Then, for any $F \in \C^1(\P)$,

		$$
		\partial_a F(\mu_a) = \langle \grad F(\mu_a), \partial_a \mu_a \rangle_{\mu_a}, \quad a \geq 0.
		$$
		In particular, $a \mapsto F(\mu_a)$ is $\C^1$.
	\end{coro}

	\begin{proof}
		Let $V_a = \left(\partial_a \varphi_a\right) \circ \varphi_a^{-1}$ and $\psi_\varepsilon(x) = \exp_x(\varepsilon V_a(x))$. Since 
		
		$$
		\sup_{x \in M} d(\varphi_{a+\varepsilon}(x),\psi_\varepsilon \circ \varphi_a(x)) = o(\varepsilon),
		$$
		we can choose $\varepsilon$ small enough to obtain the existence of a smooth flow of diffeomorphisms $(g_r^\varepsilon)_{0 \leq r \leq 1}$ given by $g_r^\varepsilon(x)=\exp_{\varphi_{a+\varepsilon}(x)}(rV^\varepsilon(\varphi_{a+\varepsilon}(x)))$ for some $V^\varepsilon \in \Gamma(TM)$ and such that $g_0^\varepsilon = \varphi_{a+\varepsilon}$ and $g^\varepsilon_1 = \psi_{\varepsilon}\circ \varphi_a$. We have, by letting $\tilde{\mu}_{a+\varepsilon} = (\psi_\varepsilon \circ \varphi_a)_*\mu$ and by applying Lemma \ref{Lemme fondamental L}:

		\begin{align*}
			\frac{F(\mu_{a+\varepsilon})-F(\mu_a)}{\varepsilon} &= \frac{F(\tilde{\mu}_{a+\varepsilon})-F(\mu_a)}{\varepsilon} + \frac{F(\mu_{a+\varepsilon})-F(\tilde{\mu}_{a+\varepsilon})}{\varepsilon}  \\
			& = \frac{F(\tilde{\mu}_{a+\varepsilon})-F(\mu_a)}{\varepsilon} + \frac{1}{\varepsilon} \int_0^1 \partial_r F((g_r^\varepsilon)_*\mu) \d r \\
			& \leq \frac{F(\tilde{\mu}_{a+\varepsilon})-F(\mu_a)}{\varepsilon} +\frac{1}{\varepsilon} \sup_{(x,\mu)}|\grad F(x,\mu)| \sup_{x \in M} d(\varphi_{a+\varepsilon}(x),\psi_\varepsilon \circ \varphi_a(x))\\
			& = \frac{F(\tilde{\mu}_{a+\varepsilon})-F(\mu_a)}{\varepsilon} +o_{\varepsilon \to 0}(1).
		\end{align*}
		The conclusion follows from the definition of $L$-differentiability.
	\end{proof}
	
	\begin{coro} \label{Proposition fonctions lipschitz}
  	Let $F \in \C^1(\P)$. Then, the functional $F$ is Lipschitz:
  	
  	$$
  	|F(\mu)- F(\nu)| \leq \sup_{(x,\mu) \in M \times\P}|\grad F(x,\mu)| W_2(\mu,\nu), \quad \forall \mu, \nu \in \P.
  	$$
  	\end{coro}
  	
  	\begin{proof}
  	Let $\mu, \nu \in \P_\infty$ and $(\mu_t)_t$ be the unique minimizing geodesic from $\mu$ to $\nu$. Applying the Benamou-Brenier formulation of the Wasserstein distance yields $W_2^2(\mu,\nu) = \int_0^1 |\dot{\mu}_t|^2 \d t$. Then, by Lemma \ref{Lemme fondamental L} and Cauchy-Schwarz inequality:
  	
  	\begin{align*}
  	    |F(\mu) -F(\nu) | &= \left| \int_0^1 \langle \grad F(\mu_t), \dot{\mu}_t \rangle_{\mu_t} \d t \right| \\
  	    &\leq \sup_{(x,\mu) \in M \times\P}|\grad F(x,\mu)| \int_0^1 |\dot{\mu}_t|_{\mu_t} \d t \\
  	    &\leq \sup_{(x,\mu) \in M \times\P}|\grad F(x,\mu)| \, \left(\int_0^1 |\dot{\mu}_t|_{\mu_t}^2 \d t \right)^{1/2} = \sup_{(x,\mu) \in M \times\P}|\grad F(x,\mu)| W_2(\mu,\nu).
  	\end{align*}
  	Since $\P_\infty$ is dense in $\P$, we obtain the result on $\P$.
  	\end{proof}

\begin{exemple} \label{Exemple fonction cyl}
	\begin{enumerate}
		\item Let $f \in \C^\infty(M)$. The functional
		
		\begin{align*}
			F_f : \P &\longrightarrow \mathbb{R} \\
			\mu & \longmapsto F_f(\mu) = \int_M f(x) \d \mu(x)
		\end{align*}
		belong to $\C^\infty_S(\P)$, with $D^\mathfrak{E}F_f(\mu)(x_0) = f(x)- F_f(\mu)$ and $\grad F_f(\mu)(x) = \nabla f(x)$. 
		\item The algebra of cylinder function

	\begin{equation*}
\mathscr{A} := \left \{ F = G (F_{f_1},F_{f_2},\dots, F_{f_k}) \, : \, k \in \mathbb{N}^*, \, G \in \C^l(\mathbb{R}^k),\, f_i \in \C^\infty(M), \, \forall i=1,2, \dots, k \right \}
\end{equation*}
belong to $\C^l_S(\P)$, with $D^\mathfrak{E} F(\mu)(x)= \sum_{i=1}^k \partial_i G(F_{f_1}(\mu),\dots F_{f_n}(\mu)) \left[ f_i(x) - F_{f_i}(\mu) \right] $ and 

$$
\grad F(\mu)(x) = \sum_{i=1}^k \partial_i G(F_{f_1}(\mu),\dots F_{f_n}(\mu)) \nabla f_i(x).
$$

In the same way, one can consider 

	\begin{equation*}
\check{\mathscr{A}} := \left \{ F = G (x,F_{f_1},F_{f_2},\dots, F_{f_k}) \, : \, k \in \mathbb{N}^*, \, G \in \C^l(M \times \mathbb{R}^k),\, f_i \in \C^\infty(M), \, \forall i=1,2, \dots, k \right \}.
\end{equation*}
This class of cylinder functions is included in $\C^{\infty,l}_S(M \times \P)$. Note that by Stone-Weierstrass theorem $\mathscr{A}$ (resp. $\check{\mathscr{A}}$) is dense in $\C^0(\P)$ (resp. $\C^0(M \times \P))$.
	\end{enumerate}
\end{exemple}

In the following, $M$ is fixed and $\P$ stands for $\P(M)$. We also define the set of \textit{$\C^k$ vector fields} on $\P$ to be:
	
	\begin{equation*}
	    \Gamma^k(T\P) = \left\{\nabla \phi(\cdot) \, : \, \phi \in \C^{\infty,k}_S(M \times \P) \right\}.
	\end{equation*}
	Throughout this article we will consider such vector fields rather than general vector fields. This choice makes simpler the introduction of essential geometric objects such as Levi-Civita connection or Lie bracket. We will also consider $\mathcal{C}^k$ vector fields on $\P_\infty$, namely,
	
	\begin{equation*}
	    \Gamma^k(T\P_\infty) = \left\{ \nabla \phi(\cdot) \, : \, \phi \in \C_S^{\infty,k}(M \times \P_\infty) \right\}.
	\end{equation*}

	\begin{exemple}
    	    A representative example of a subspace of $\Gamma^\infty(T\P)$ is the following:
    	    
    	    \begin{equation*}
    	        \left\{ \sum_{i=1}^l F_i V_{\phi_i}, \quad F_i \in \mathscr{A}, \, \phi_i \in \C^\infty(M), \, l \in \mathbb{N}^* \right\},
    	    \end{equation*}
    	     where $\mathscr{A}$ is defined in Example \ref{Exemple fonction cyl}.
    	    \end{exemple}

    \vspace{10pt}
\noindent \textbf{Levi-Civita connection.} Let $\bar{Z} : \mu \in \P \mapsto \bar{Z}(\mu) = \nabla \phi(\cdot,\mu) \in \Gamma^1(T\P)$. The \textit{Levi-Civita connection} $\grad$ on $\T \P$ is the unique connection on $\T \P$ which is metric, i.e. for $\bar{Z}_1, \bar{Z}_2 \in \Gamma^1(TP)$:

$$
\L_{v} \langle \bar{Z}_1(\mu), \bar{Z}_2(\mu) \rangle_\mu = \langle \grad_{v} \bar{Z}_1(\mu), \bar{Z}_2(\mu) \rangle_\mu + \langle  \bar{Z}_1(\mu), \grad_{v} \bar{Z}_2(\mu) \rangle_\mu, \quad \mu \in \P, \, v \in \T_\mu \P,
$$
and torsion-free, namely,

$$
[\bar{Z}_1, \bar{Z}_2]F_f(\mu) = \left(\L_{\bar{Z}_1} \L_{\bar{Z}_2} - \L_{\bar{Z}_2} \L_{\bar{Z}_1} \right) F_f(\mu) = \left\langle \grad_{\bar{Z}_1} \bar{Z}_2(\mu) -\grad_{\bar{Z}_2} \bar{Z}_1(\mu), \grad F_f(\mu) \right\rangle_\mu, \quad \mu \in \P, \, v \in \T_\mu \P
$$
The Levi-Civita connection admits the following explicit form:

$$
\grad_{V_f} \bar{Z}(\mu) = \bar{\mathcal{L}}_{V_f} \bar{Z} + \Pi_{\mu} \left( \nabla_{\nabla f} \nabla \phi(\mu,\cdot) \right) = \nabla \bar{\mathcal{L}}_{V_f} \phi(\mu,\cdot) + \Pi_{\mu} \left( \nabla_{\nabla f} \nabla \phi(\mu,\cdot) \right),
$$
where $\Pi_\mu : L^2(M,TM;\mu) \rightarrow \T_\mu \P$ denotes the orthogonal projection obtained by the Hodge decomposition.

Let $F \in \C^2(\P)$ and let $\bar{Z}_1, \bar{Z}_2 \in \Gamma^1(T\P)$, the quantity
    
    $$
    \bar{\he } F(\bar{Z}_1, \bar{Z}_2) = \L_{\bar{Z}_1} \L_{\bar{Z}_2} F -  (\grad_{\bar{Z}_1} \bar{Z}_2) F
    $$
    is well defined. We call this quantity the \textit{Hessian} of $F$. Note that, as in the finite-dimensional case, we have the following identity:
    
    $$
    \bar{\he } F(\mu)(\bar{Z}_1, \bar{Z}_2) = \left\langle \grad_{\bar{Z}_1(\mu)}\grad F(\mu), \bar{Z}_2(\mu) \right\rangle_\mu.
    $$
    Moreover, since
    
    $$
    \bar{\he } F(\bar{Z}_1, \bar{Z}_2) - \bar{\he } F(\bar{Z}_2, \bar{Z}_1) = \bar{\mathcal{L}}_{ [\bar{Z}_1, \bar{Z}_2]} F  - \bar{\mathcal{L}}_{[\bar{Z}_1,\bar{Z}_2]}F = 0,
    $$
    $\bar{\he }F$ is clearly symmetric. Furthermore, since $\bar{\he }F$ is tensorial in the first component and symmetric, we deduce that it is tensorial in both components. With the above notations,

	$$
	\bar{\he}F(\mu)(V_f,V_g) = \int_{M^2} \langle \grad^2 F(\mu)(x,y) , \nabla g(x) \otimes \nabla f(y) \rangle \d \mu^{\otimes 2}(x,y) + \int_M \langle \nabla \grad F(\mu)(x), \nabla g(x) \otimes \nabla f(x)\rangle \d \mu(x).
	$$
    Indeed,

	\begin{align*}
		\bar{\he}F(\mu)(V_f,V_g) & = \left\langle \grad_{V_f}\grad F(\mu), V_g \right\rangle_\mu \\
		& = \int_M \left\langle \bar{\mathcal{L}}_{V_f} \grad F(\mu)(x), \nabla g(x) \right\rangle \d \mu(x) + \int_M \left\langle\nabla_{\nabla f} \grad F(\mu)(x), \nabla g(x) \right\rangle \d \mu(x) \\
		& = \int_{M^2} \left\langle \grad^2 F(\mu)(x,y), \nabla g(x) \otimes \nabla f(y) \right\rangle \d \mu^{\otimes 2}(x,y) + \int_M \left\langle\nabla_{\nabla f} \grad F(\mu)(x), \nabla g(x) \right\rangle \d \mu(x).
	\end{align*}

	\begin{rem}
		Note that the above formula makes clear the difference between $\langle \grad^2 F, V_\phi \otimes \cdot \rangle$ and $\langle \grad_{V_\phi} \grad F, \cdot \rangle$. 
	\end{rem}

\noindent \textbf{The group of diffeomorphisms.} For all $s \in \mathbb{N}$, let $\mathscr{M}^s$ be a Hilbert manifold such that $\mathscr{M}^{s+1} \subset \mathscr{M}^s$ and the injection $\iota : \mathscr{M}^{s+1} \rightarrow \mathscr{M}^s$ is smooth. We consider the following Inverse Limit Hilbert (ILH for short) manifold:

\begin{equation}
	\mathscr{M} = \bigcap_{s \geq 0} \mathscr{M}^s.
\end{equation}
Let $\mathscr{N} = \cap_{s \geq 0} \mathscr{N}^s$ be another ILH manifold.
A map $f : \mathscr{M} \rightarrow \mathscr{N}$ is said to be of class $\C^k$ for the ILH topology if for all $s$ there exists $j(s) \geq s$ such that $f$ admits a $\C^k$ extension $f^s : \mathscr{M}^{j(s)} \rightarrow \mathscr{N}^s$.

Let $M$ be a closed manifold and let $\pi : E$ be a fiber bundle over $M$. 
The set $\Gamma^s(E)$ of $H^s$ sections of $E$ has a natural structure of Hilbert manifold, see \cite{Palais66} for instance.
This is the case for the bundle $M \times N$ where $N$ is another Riemannian manifold.
Note that if $E$ is a vector bundle, $\Gamma^s(E)$ can be endowed with a natural Hilbert space structure. 
Classical examples are given by $E= M \times \mathbb{R}$ or $E= TM$.
The group of diffeomorphisms $\di := \mathrm{Diff}(M,M)$ is an ILH manifold. 
More precisely, for each $s> m/2$, we consider the following Hilbert manifold:

\begin{equation}
	\di^s = \left\{ \varphi \in H^s(M,M) \text{ bijective} \, : \, \varphi^{-1} \in H^s(M,M) \right\},
\end{equation}
and the inverse limit description is given by:

\begin{equation}
	\di = \bigcap_{s >m/2} \di^s.
\end{equation}
It is well known that $\mathscr{D}$ carries an Infinite-dimensional Lie group structure with respect to the inverse limit topology, even though, for every $s >m/2$, the Hilbert manifold $\mathscr{D}^s$ fails to be a Lie group, see \cite[Section 2]{EbinMarsden1970} for instance. Moreover, the projection 

\begin{align*}
    p : \mathscr{D} &\longrightarrow \P_\infty \\
     \varphi &\longmapsto (\varphi)_* \vol,
\end{align*}
is onto and continuous. Throughout the rest of the article we will let $\mathscr{D}_\mu$ denotes the \textit{fiber} above $\mu$, i.e $\mathscr{D}_\mu = p^{-1}(\mu)$.
    
The \textit{tangent space} of $\mathscr{D}$ in $\varphi \in \mathscr{D}$ is
    
    $$
    T_{\varphi} \mathscr{D} := \Gamma(\varphi^*TM),
    $$
    and the tangent space of $\mathscr{D}^s$ in $\varphi \in \mathscr{D}^s$ is defined to be:
    
    $$
    T_{\varphi} \mathscr{D}^s := \Gamma^{(s)}(\varphi ^*TM),
    $$
    In particular, for $\varphi \in \mathscr{D}$, $T_\varphi \mathscr{D} = \cap_{s > m/2} T_\varphi \mathscr{D}^s$. Moreover, $T_{id}\mathscr{D}^s = \Gamma^{(s)}(TM)$ and $T_\varphi \mathscr{D}^s \cong T_{id} \mathscr{D}^s$ for all $\varphi \in \mathscr{D}^s$, the isomorphism being given by $\bm T R_\varphi : A \mapsto A \circ \varphi$. The set $T \mathscr{D} =  \underset{\varphi \in \mathscr{D}}{\sqcup} T_{\varphi}\mathscr{D}$ is the tangent bundle of $\mathscr{D}$ and it admits a structure of smooth ILH vector bundle over $\mathscr{D}$.
    
    \begin{defi}
    A \textit{vector field} on $\mathscr{D}$ is a map $\bm A : \mathscr{D} \rightarrow T \mathscr{D}$ such that $\bm A(\varphi) = A(\varphi) \circ \varphi$ with $A(\varphi) \in \Gamma(TM)$ for all $\varphi \in \mathscr{D}$.
    \end{defi}  
    
    We can build \textit{right invariant vector fields} on $\mathscr{D}$ similarly as on Lie groups by setting $ \check{A}(\varphi) = A \circ \varphi$ for $A \in T_{\mathrm{id}}\mathscr{D}$.
    
    \begin{defi}\begin{enumerate}
        \item A vector field $\bm A$ on $\mathscr{D}$ is said to be of class $\C^k$ if $\bm A : \mathscr{D} \rightarrow T\mathscr{D}$ is $\C^k$ as a ILH manifold map. We denote by $\Gamma^k(T\mathscr{D})$ the set of $\C^k$ vector fields on $\mathscr{D}$.
        \item A vector field $\bm A \in \Gamma^k(T\mathscr{D})$ on $\mathscr{D}$ is said to be a $\C^k$ stratified vector field if there exists $s_0$ such that for all $s > s_0$ there exists a $\C^k$ extension $A^s : \mathscr{D}^s \rightarrow T\mathscr{D}^s$. The set of $\C^k$ stratified vector fields is denoted by $\Gamma^k_S(T \mathscr{D})$. 
    \end{enumerate}
        
    \end{defi}
    In the following we will only consider stratified vector fields. Let $s > m/2$, we consider the following $L^2$ \textit{Riemannian metric} on $T\mathscr{D}^s$:
    
    \begin{equation} \label{E Métrique Riemanienne}
        \langle A \circ \varphi, B \circ \varphi \rangle_{\varphi} = \int_M \langle A(\varphi(x)), B(\varphi(x)) \rangle \d\vol(x), \quad A \circ \varphi, B \circ \varphi \in T_\varphi \mathscr{D}^s.
    \end{equation}
    We also endow $T\mathscr{D}$ with this Riemannian metric. Let us define
    
    $$
    G := \left \{ g \in \mathscr{D} \, : \, p(g) = \vol \right \} = \mathscr{D}_{\vol}.
    $$
    Note that the group $\mathscr{D}$ has a structure of $G$-principal bundle over $\P_\infty$. Note also that the Riemannian metric \eqref{E Métrique Riemanienne} on $T\mathscr{D}$ is $G$-right invariant.
    
	For $\bm F \in \C^1(\mathscr{D})$, the Lie derivative in the direction $A \circ \varphi \in T_\varphi \mathscr{D}$ is denoted by

	$$
	\mathcal{L}^{\mathscr{D}}_{A \circ \varphi} \bm F(\varphi).
	$$

    Let $k \geq 1$. The \textit{Lie bracket} between two vector fields $\bm A, \bm B \in \Gamma^k_S(T\mathscr{D})$ is defined to be the unique vector field satisfying:
    
    $$
    [\bm A, \bm B] \bm F(\varphi) = \mathcal{L}^\mathscr{D}_{\bm A} \mathcal{L}^\mathscr{D}_{\bm B} \bm F(\varphi) - \mathcal{L}^\mathscr{D}_{\bm B} \mathcal{L}^\mathscr{D}_{\bm A} \bm F(\varphi), \quad \forall \bm F \in \C^{2}(\mathscr{D}).
    $$
    Note that $[\bm A,\bm B] \in \Gamma^{k-1}_S(T\mathscr{D})$. 

	The Levi-Civita connection $\bm \nabla$ on $T\mathscr{D}$ has been proved (see \cite[Theorem 9.1]{EbinMarsden1970}) to be the unique torsion free linear connection, namely,
    
    \begin{equation} \label{Torsion free D}
        \bm \nabla_{\bm A_1} \bm A_2 - \bm \nabla_{\bm A_2} \bm A_1 = [\bm A_1, \bm A_2],
    \end{equation}
    which is compatible with the Riemannian metric, i.e.
    
    \begin{equation} \label{Compatible with metric D}
        \mathcal{L}^{\mathscr{D}}_{\bm C \circ \varphi} \langle \bm A_1(\varphi), \bm A_2 (\varphi) \rangle_\varphi = \langle \bm \nabla_{C \circ \varphi} \bm A_1 (\varphi), \bm A_2 (\varphi) \rangle_\varphi + \langle \bm A_1(\varphi), \bm \nabla_{C \circ \varphi} \bm A_2(\varphi) \rangle_\varphi,
    \end{equation}
	for all $\bm A_1, \bm A_2 \in \Gamma^{1}_S(T\mathscr{D})$.

\begin{figure}[h]
     \centering
     \includegraphics[width=0.8\textwidth]{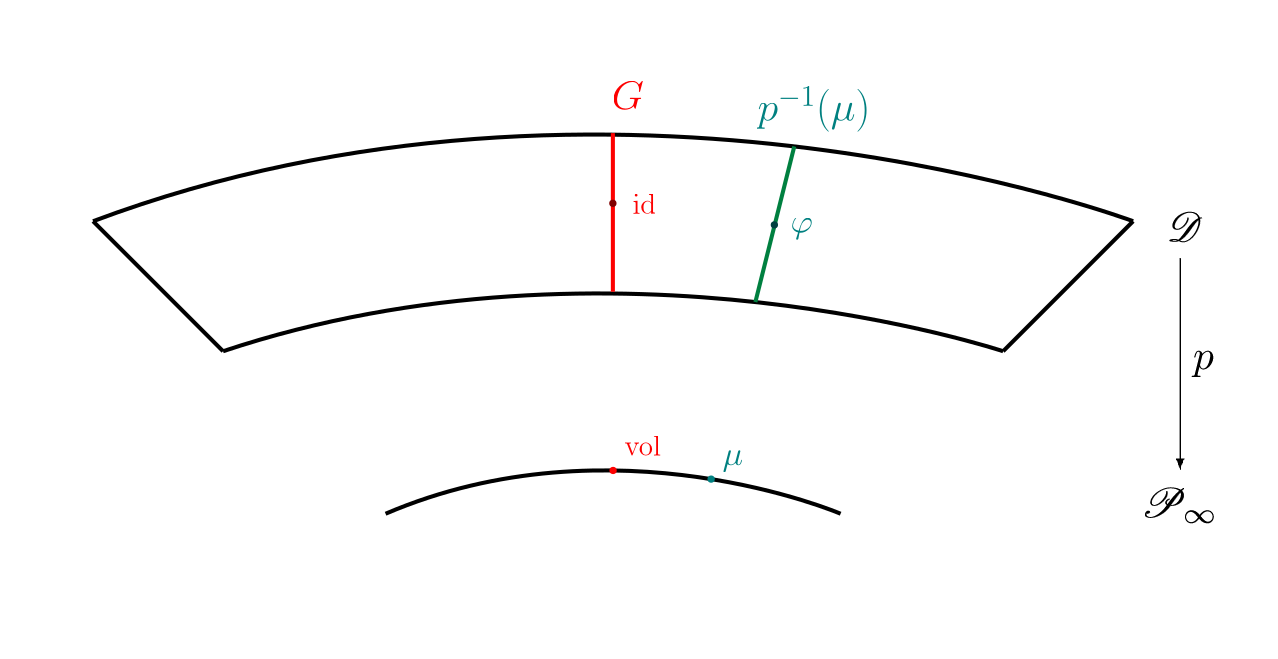}
     \caption{Illustration of the $G$-principal bundle structure of $\mathscr{D}$ over $\P_\infty$.}
 \end{figure}

	\vspace{10pt}
	\noindent \textbf{The Riemannian submersion structure.} The projection $p : \mathscr{D} \rightarrow \P_\infty$ carries a natural Riemannian submersion structure, with a the following splitting of the tangent bundle:
	
	$$
	T\mathscr{D} = \V \mathscr{D} \oplus \H \mathscr{D},
	$$
	with
	$$
    \H_{\varphi} \mathscr{D} = \left \{ \nabla \phi \circ \varphi \, : \, \phi \in \C^\infty(M) \right \}
    $$
	and
	$$
    \V_\varphi \mathscr{D} = \left \{ Y \in \Gamma(\varphi^*TM), \, \text{div}_{p(\varphi)}(Y) = 0 \right \},
    $$
	where $\mathrm{div}_{p(\varphi)}$ is the divergence associated to the measure $p(\varphi)$. It also writes

	$$
	\mathrm{div}_{p(\varphi)} = \mathrm{div} + \nabla \log (\rho),
	$$
	where $\mathrm{div}$ is the divergence with respect to the volume measure and $\rho = (\d p(\varphi))/ \vol$. The right-inverse of $\bm T p|_{\H \mathscr{D}} : \H \mathscr{D} \rightarrow T \P_\infty$ is given by $\mathfrak{h} : T\P_\infty \rightarrow \H\mathscr{D}$ where for all $\varphi \in \mathscr{D}$,

	$$
	\mathfrak{h}_\varphi(V_\phi) = \nabla \phi \circ \varphi.
	$$
	It is shown in \cite{Martin2025} that $\H \di$ and $\V \di$ are smooth subbundles of $T \di$ and that the projections $P_{\H} : T \di \rightarrow \H \di$ and $P_{\V}\di : T \di \rightarrow \V \di$ induced by the Hodge decomposition are smooth vector bundle morphisms, in the sense of definition \ref{Defi vector bundle morphism}. Note that the notion of $L$-differentiability is compatible with the notion of differentiability in $\mathscr{D}$ since the charts are defined by the exponential map (see for instance \cite{EbinMarsden1970}). Furthermore, the following result holds:

	\begin{prop}
		If $F \in \C^1_S(\P)$, then $F \circ p \in \C^1(\mathscr{D})$ and $\bm \nabla (F \circ p) \in \Gamma^0_S(T\mathscr{D})$. Moreover, if $\bar{Z} \in \Gamma^k(T\P)$, then $\bm Z = \mathfrak{h}(\bar{Z}) \in \Gamma^k_{S}(T\mathscr{D})$.
	\end{prop}

	\begin{proof}
		Let $F \in \C^1_S(\P)$, one can check that $\bm \nabla (F \circ p) (\varphi) = \grad F(p(\varphi)) \circ \varphi$. Since $p : \mathscr{D} \to \P_\infty$ is continuous, $F \in \C^1(\mathscr{D})$ and $\bm \nabla (F \circ p) \in \Gamma^0_S(T\mathscr{D})$. In the same way, if $\bar{Z} \in \Gamma^1(T\P)$ with $\bar{Z} : \mu \in \P \mapsto \nabla \phi(\mu,\cdot) \in T_\mu \P$, one can easily check that $\bm Z(\varphi) := \mathfrak{h}(\bar{Z})(\varphi) = \nabla \phi(p(\varphi)) \circ \varphi$ for any $\varphi \in \mathscr{D}^s$ for $s >m/2$. Let $C \circ \varphi \in T_\varphi \mathscr{D}^s$, its covariant derivative is given by

	\begin{align*}
		& \quad \bm \nabla_{C \circ \varphi} \bm Z(\varphi) \\
		& = \left( \nabla_C \nabla \phi(\varphi) \right) \circ \varphi + \mathcal{L}^\mathscr{D}_{C \circ \varphi} \left(\nabla  \phi(\varphi) \right) \circ \varphi \\
		& = \int_M \langle \nabla^2 \phi (p(\varphi),\varphi(x)), \Pi_{p(\varphi)}(C)(\varphi(x)) \otimes \cdot \rangle \d \vol(x) + \int_M \langle \grad \nabla \phi(p(\varphi),\varphi(x))(\varphi(y)), \Pi_{p(\varphi)}(C)(\varphi(y)) \otimes \cdot \rangle \d \vol(y).
	\end{align*}
	Thus, since $\varphi \in \mathscr{D}^s$, $p : \varphi \in \mathscr{D}^s \mapsto (\varphi)_* \vol \in \P$ is continuous, $(x,\mu) \mapsto \nabla \phi(\mu, x) \in \C^{\infty,1}(M \times \P, TM \otimes TM)$, $(x,y,\mu) \mapsto \grad \nabla \phi(\mu,x,y) \in \C^{\infty,\infty,0}(M \times M \times \P, TM \otimes TM)$ and since the map $(\varphi,A) \in \mathscr{D}^s \times \Gamma(TM) \mapsto \A \circ \varphi $ is smooth (see \cite{EbinMarsden1970} for example), we obtained that $\bm Z \in \Gamma^1_S(T\mathscr{D})$. The same reasoning can be performed to show that for any $\bar{Z} \in \Gamma^k(\T \P)$, $\bm Z = \mathfrak{h} (\bar{Z}) \in \Gamma^k_S(T\mathscr{D})$. 
	\end{proof}

\subsection{Riemann curvature tensors} \label{S 2.2} The \textit{normal tensor} is the antisymmetric linear map defined as:
	
	\begin{align*}
		\mathcal{N}_{\varphi} : T_\varphi \mathscr{D} \times T_\varphi \mathscr{D}  &\longrightarrow \V_\varphi \mathscr{D}. \\
		(U \circ \varphi, A \circ \varphi) &\longmapsto \mathcal{N}_{\varphi} (U \circ \varphi, A \circ \varphi) = P_{\V_{\varphi}}\left(  ((\nabla U)^T \cdot A) \circ \varphi  \right),
	\end{align*}
    Let $\bm Z \in \Gamma^1_S(\H \mathscr{D})$, then, for any $\varphi \in \mathscr{D}$ and any $A \circ \varphi \in T_\varphi \mathscr{D}$, the following holds

	$$
	P_{\V_{\varphi}} \left( \bm \nabla_{A \circ \varphi} \bm Z(\varphi) \right) = \mathcal{N}_\varphi(\bm Z(\varphi), A \circ \varphi).
	$$
	In particular, if we denote by $\mathcal{O}_{U \circ \varphi}$ the linear counterpart of the normal tensor, namely,

	$$
		\mathcal{O}_{U \circ \varphi}(A\circ \varphi) = \mathcal{N}_{\varphi} (U \circ \varphi,A \circ \varphi), \quad \forall A \circ \varphi \in T_\varphi\mathscr{D}^s,
	$$
	then, for $\bm Z \in \Gamma^k_S(\H \mathscr{D})$,

	$$
	\mathcal{O}_{\bm Z} = P_\V \circ \bm \nabla \bm Z,
	$$
	and consequently belongs to $\mathrm{Hom}_{\mathrm{id}}^{k-1}(T\mathscr{D}, T\mathscr{D})$, see Definition \ref{Defi vector bundle morphism}. Moreover, note that for a $G$-right invariant vector field $\bm Z \in \Gamma^1(\H \mathscr{D}^s)$, the following property holds:
	
	$$
	\mathcal{O}^\top_{\bm Z} = \bm \nabla \bm Z \circ P_\V,
	$$
	where $\mathcal{O}^{\top}_{\bm Z}$ is the adjoint for the Riemannian metric on $T\mathscr{D}$.
	
	\noindent \textbf{Riemann curvature operator on $\di$.}
	
	\begin{defi}
		The Riemann curvature operator on $\di$ is defined to be
		
		$$
		\bm R(\bm A_1,\bm A_2)\bm A_3(\varphi) = \bm \nabla_{\bm A_1} \bm \nabla_{\bm A_2} \bm A_3 (\varphi) - \bm \nabla_{\bm A_2} \bm \nabla_{\bm A_1} \bm A_3 - \bm\nabla_{[\bm A_1,\bm A_2]} \bm A_3,
		$$
		for all $\bm A_1, \bm A_2, \bm A_3 \in \Gamma^2(T\mathscr{D}^s)$.
	\end{defi}
	
	\begin{prop} \label{Forme sympa tenseur riemann D}
		Let $\bm A_1, \bm A_2, \bm A_3,\bm A_4 \in \Gamma^2(T \mathscr{D})$, we have
		
		\begin{equation} \label{Tenseur Riemann difféos}
		\left \langle (\bm R(\bm A_1,\bm A_2)\bm A_3)(\varphi), \bm A_4(\varphi) \right \rangle_\varphi = \int_M \left \langle R(A_1(\varphi)(\varphi(x)),A_2(\varphi)(\varphi(x)))A_3(\varphi)(\varphi(x)), A_4(\varphi)(\varphi(x)) \right \rangle d \vol(x).
		\end{equation}
	\end{prop}
	
	\begin{prop}
		The curvature operator given by \eqref{Tenseur Riemann difféos} inherits the symmetries of the Riemann curvature operator on $M$ and is tensor.
	\end{prop}
	
	\begin{proof}
	The fact that this operator inherits the symmetries of the Riemann curvature operator on $M$ is a direct consequence of the above proposition. Now, it is clear that the operator given by \eqref{Tenseur Riemann difféos} is tensor in the fourth coordinate. The conclusion follows from the symmetries.
	\end{proof}
	
	\noindent \textbf{Riemann curvature tensor on $\P$.}
	
	\begin{defi}
	The Riemann curvature operator on $\P$ is defined to be
	
	$$
	\bar{R}(\bar{Z}_1,\bar{Z}_2)\bar{Z}_3 = \grad_{\bar{Z}_1} \grad_{\bar{Z}_2} \bar{Z}_3 - \grad_{\bar{Z}_2} \grad_{\bar{Z}_1} \bar{Z}_3 - \grad_{[\bar{Z}_1, \bar{Z}_2]} \bar{Z}_3,
	$$
	where $\bar{Z}_1, \bar{Z}_2, \bar{Z}_3 \in \Gamma^2(\T\P)$.
	\end{defi}
	
	\begin{prop}[\cite{gigli2011second}, Proposition 6.4]
	Let $\mu \in \P_\infty$, $\bar{Z}_1,\bar{Z}_2,\bar{Z}_3,\bar{Z}_4 \in \Gamma^2_S(T\mathscr{D})$ and $\bm Z_1, \bm Z_2, \bm Z_3, \bm Z_4 \in \Gamma^2(T\mathscr{D})$ their horizontal lifts, we have
	
	\begin{align*}
	\left \langle \bar{R}(\bar{Z}_1,\bar{Z}_2)\bar{Z}_3, \bar{Z}_4 \right \rangle_{\mu} & = \left \langle \bm R (\bm Z_1, \bm Z_2) \bm Z_3, \bm Z_4 \right \rangle_\varphi  + \left \langle \mathcal{N}_{\varphi}(\bm Z_3, \bm Z_2), \mathcal{N}_\varphi(\bm Z_4, \bm Z_1 ) \right \rangle_\varphi \\
	& - \left \langle \mathcal{N}_\varphi(\bm Z_3, \bm Z_1), \mathcal{N}_\varphi( \bm Z_4, \bm Z_2) \right \rangle_\varphi
 + 2\left \langle \mathcal{N}_{\varphi}(\bm Z_3, \bm Z_4), \mathcal{N}_\varphi(\bm Z_2, \bm Z_1) \right \rangle_\varphi,
	\end{align*}
	for all $\varphi \in \di_\mu$. In particular, $\bar{R}$ is a tensor and inherits the same symmetries as $R$.
	
	\end{prop}
	
	\begin{coro} \label{Trace tenseur Riemann sur P}
	Let $\mu \in \P_\infty$, $Z_1, Z_2, U \in T_\mu \P_\infty$. Then, the following holds:
	
	$$
	\left \langle \bar{R}(U, Z_1) Z_1, Z_2 \right \rangle_\mu = \left \langle \bm R (\mathfrak{h}_\varphi(U), \mathfrak{h}_\varphi(Z_1)) \mathfrak{h}_\varphi(Z_1), \mathfrak{h}_\varphi(Z_2) \right \rangle_\varphi + 3 \langle \mathcal{N}_{\varphi}(\mathfrak{h}_\varphi(Z_1), \mathfrak{h}_\varphi(U)), \mathcal{N}_{\varphi}(\mathfrak{h}_\varphi(Z_1), \mathfrak{h}_\varphi(Z_2)) \rangle_\varphi,
	$$
	for all $\varphi \in \mathscr{D}$. Equivalently,

	\begin{equation}
		 \bm R (\mathfrak{h}_\varphi(U), \mathfrak{h}_\varphi(Z_1)) \mathfrak{h}_\varphi(Z_1) = \mathfrak{h}_{\varphi}\left(\bar{R}(U, Z_1) Z_1\right) - 3 \mathcal{O}^*_{\mathfrak{h}_\varphi(Z_1)} \left( \mathcal{O}_{\mathfrak{h}_\varphi(Z_1)} \left( \mathfrak{h}_\varphi(U) \right) \right).
	\end{equation}
	\end{coro}
	
	\begin{rem}
		Interpreting $\di$ as a Riemannian submersion onto $\P_\infty$, this result is not surprising at all. Indeed, let $p : N \rightarrow M$ be a Riemannian submersion, it is just the infinite-dimensional equivalent of \cite[Corollary 1]{oneill1966fundamental} which states that for $X,Y \in \Gamma(N)$ horizontal vector fields and $p$-related to $\bar{X}, \bar{Y} \in \Gamma(M)$,
		
		$$
		\bar{K}(\bar{X},\bar{Y}) = K(X,Y) + \frac{3}{4} \frac{\|[X,Y]^V  \|^2}{\|X\|^2\|Y\|^2 - \langle X,Y \rangle^2},
		$$
		where $K, \bar{K}$ denote the sectional curvature of $N, M$.
		
	\end{rem}

		\begin{exemple} \label{Exemple tore}
	The Wasserstein curvature of $\P_\infty(\mathbb{S}^1)$ is non-negative. Moreover, we show in this example that the normal part of a vector field can be explicitly computed on $\mathbb{S}^1$ which is not generally the case.
	
	 The fact that $\mathbb{S}^1$ is parallelizable allows us to see vector fields $V$ as functions from $\mathbb{S}^1$ to $\mathbb{R}$. Moreover, a vector field $v$ will be a gradient of function if and only if $\int_{\mathbb{S}^1} v(x) dx = 0$. Let $\mu \in \P_\infty(\mathbb{S}^1)$, $\rho = \d \mu / \d \vol$, $v \in \Gamma(\mathbb{S}^1)$ and $\partial_x \phi = \Pi_\mu(v)$. Let $f \in \C^\infty(\mathbb{S}^1)$, we have:
	
	$$
	\int_M \partial_x f(x) v(x) \rho \d \vol(x) = \int_M \partial_x f(x) \partial_x \phi(x) \rho \d \vol(x).
	$$
	By using integration by parts we get
	
	$$
	\int_M f(x) \partial_x( v(x)\rho(x)) \d \vol(x) = \int_M f(x) \partial_x (\partial_x \phi(x) \rho(x)) \d \vol(x).
	$$
	Thus, since the set of divergence-free vector fields is the set of constant vector fields, we obtain the following equality in $L^2(\vol )$:
	
	$$
	v(x)\rho(x) = \partial_x \phi(x) \rho(x) + K.
	$$
	Dividing by $\rho$ on both sides and integrating we get
	
	$$
	K \int_{\mathbb{S}^1} \frac{1}{\rho} \d \vol(x) = \int_{\mathbb{S}^1} v(x) \d \vol(x),
	$$
	thus, we get the following decomposition in $\T_\mu \oplus \T_\mu^\perp$ :
	
	$$
	v = \partial_x \phi + \hat{\rho} \int_{\mathbb{S}^1} v(x) \d \vol(x)
	$$
	where $\hat{\rho} = \frac{1}{\rho} \cdot \frac{1}{\int_{\mathbb{S}^1} \frac{1}{\rho(x)}dx}$. In particular, we obtain:
	
	$$
	\int_{\mathbb{S}^1} | \Pi_{\mu}^\perp(v)(x) |^2 \d \mu(x) = \left( \int_{\mathbb{S}^1} v(x) \d \vol(x) \right)^2 \frac{1}{\int_{\mathbb{S}^1} \frac{\d \vol(x)}{\rho(x)} }.
	$$	
	As a consequence, for $\phi, \psi \in \C^\infty(M)$, Corollary \ref{Trace tenseur Riemann sur P} yields:
	
	\begin{align*}
	\langle \bar{R}(V_{\phi}, V_{\psi}) V_{\psi}, V_{\phi} \rangle  &= \int_{\mathbb{S}^1} \langle R(\nabla \phi(x), \nabla \psi(x)) \nabla \psi(x), \nabla \phi(x) \rangle \d \mu(x) \\
	& + \left( \int_{\mathbb{S}^1} \partial_x \phi(x) \partial_x^2 \psi(x) \d \vol(x) \right)^2 \frac{3}{\int_{\mathbb{S}^1}\frac{1}{\rho(x)}dx} \\
	& = \left( \int_{\mathbb{S}^1} \partial_x \phi(x) \partial_x^2 \psi(x) \d \vol(x) \right)^2 \frac{3}{\int_{\mathbb{S}^1}\frac{1}{\rho(x)}dx}
	\end{align*}
	since $R \equiv 0$ on $\mathbb{S}^1$. A consequence of this computation is the following: the sectional curvature of $\P_\infty(\mathbb{S}^1)$ is not bounded since it is possible to find gradient of functions $\partial_x \phi$, $\partial_x \psi$ such that for any $C \in \mathbb{R}_+$,
	
	$$
	\int_{\mathbb{S}^1} |\partial_x \phi(x)|^2 \d \vol(x) = \int_{\mathbb{S}^1} |\partial_x \psi(x)|^2 \d \vol(x) = 1 \, \, \text{ and } \, \, \int_{\mathbb{S}^1} \partial_x \phi(x) \partial_x^2 \psi(x) \d \vol(x) \geq C.
	$$
	A very simple example is given by $\partial_x \phi(x) = \cos(nx)/\sqrt{\pi}$, $\partial_x \psi(x) = \sin(nx)/ \sqrt{\pi}$. In particular, it is impossible to construct regular geodesic balls on $\P_\infty(\mathbb{S}^1)$ in the sense of \cite[Definition 1.6]{Kendall1990}.
\end{exemple}

\subsection{Regular paths on \texorpdfstring{$\P$}{P}}	
	
		We start this subsection by defining what is a regular path on $\P$ and we introduce some examples.

	\begin{defi}
	A continuous path $(\mu_a)_{0 \leq a \leq 1}$ on $\P_\infty$ is said to be regular if there exists a $\C^1$ path of diffeomorphisms $(U_a)_{0 \leq a \leq 1}$ such that $\mu_a= (U_a)_*\mu$. Let $\mu, \nu \in \P_\infty$, the set of regular paths $(\mu_a)_{0 \leq a \leq 1}$ with $\mu_0 =\mu$ and $\mu_1 = \nu$ is denoted by $\mathscr{R}(\mu,\nu)$.

	\end{defi}

	\begin{exemple}[Linear interpolation]
		Let $\mu_0, \mu_1 \in \P_\infty$. The linear interpolation $(\mu_a)_{0 \leq a \leq 1}$ between $\mu_0$ and $\mu_1$ defined as
		
		$$
		\mu_a = (1-a)\mu_0 + a \mu_1, \quad 0 \leq a \leq 1
		$$
		is a regular path on $\P_\infty$. Indeed, for all $0 \leq a \leq 1$, let $\rho_a$ denotes the density $ \d\mu_a = \rho_a \d \vol$. It is easy to see that $\rho_a = (1-a) \rho_0 + a \rho_1$, thus, $(\mu_a)_{0 \leq a \leq 1}$ lies in $\P_\infty$. Note that
	
$$
\frac{\d}{\d a}  \int_M f(x) \d\mu_a(x) = \int_M f(x) (\rho_1(x) - \rho_0(x)) \d \vol(x),
$$	
Since $\rho_1 - \rho_0$ is a smooth function with $\int_M \rho_1 - \rho_0 \d \vol = 0$, there exists a unique smooth solution to
	
	$$
	\mathrm{div}(X) = \rho_0 - \rho_1,
	$$
	given by $X = \nabla \Delta^{-1} (\rho_0 - \rho_1)$. In particular, $(a,x) \mapsto X_a(x)= X(x)/\rho_a(x)$ is a smooth family of smooth vector fields. Thus, the unique solution to
	
	$$
	\partial_a U_a = X_a(U_a), \quad U_0= \varphi \in \mathscr{D}_{\mu_0}
	$$
	is a smooth path of diffeomorphisms. Let $\nu_a = (U_a)_* \vol$, then,

	$$
	\frac{\d}{\d a} \int_M f(x) \d \nu_a(x) = -\int_M f(x) \mathrm{div}(\rho_a X_a) \d \vol(x) = \int_M f(x) (\rho_1 -\rho_0)(x) \d \vol(x).
	$$
	As a consequence, $(\mu_a)_{0 \leq a \leq 1} = (\nu_a)_{0 \leq a \leq 1}$. The reader should keep in mind that $(\mu_a)_{0 \leq a \leq 1}$ does not minimize the Wasserstein distance and is not adapted to many problems of optimal transport.
	\end{exemple}

	\begin{figure}[h] \label{interpolation lineaire gaussiennes}
		\centering
		\includegraphics[width=0.7\textwidth]{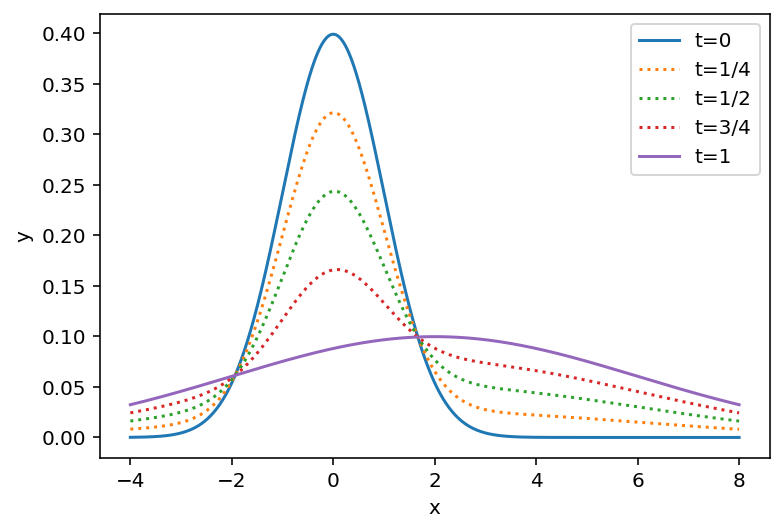}
     \caption{The linear interpolation between two Gaussians $N(0,1)$ and $N(2,4)$ at time $t\in \left\{1,\frac{1}{4},\frac{1}{2},\frac{3}{4},1\right\}$. In particular, it appears clearly that the linear interpolation between two Gaussians is not a Gaussian (Compare with Figure \ref{Transport optimal gaussiennes} below).}
	\end{figure}

	The question of whether the displacement interpolation can be represented as the pushforward of a smooth path of diffeomorphisms remains open. The seminal works of Ma, Trudinger, and Wang, followed by those of Loeper and Villani (see \cite{Loep2006}, \cite{LoeperVillani2010}, \cite{MaTrudingerWang2005}, see also \cite{Fig2010} for a survey of these results), provide sufficient geometric conditions on the underlying manifold to ensure the smoothness of the optimal transport map $T=\exp(\nabla\phi)$. Even under these assumptions, the maps $T_a=\exp(a\nabla\phi)$ need not be diffeomorphisms for all $a\in[0,1]$, as their differentials may become singular for some values of $a$. The following result nevertheless remains valid and enables us to work with the displacement interpolation at least in an infinitesimal way.

	\begin{prop}
		Let $\mu \in \P_\infty$, $A \in \Gamma(TM)$. Then, there exists $b>0$ small enough such that the path $(\mu_{ab})_{0 \leq a \leq 1}$ defined by
		
		$$
	\mu_{ab} = (\exp(ab A))_*\mu, \quad 0 \leq a \leq 1
	$$
	belongs to $\mathscr{R}(\mu,\mu_b)$ with $\mu_b \in \P_\infty$.

	\end{prop}
	\begin{proof}
		This proposition is just a consequence of the fact that the exponential is a diffeomorphism on a sufficiently small neighborhood of any diffeomorphism (take for instance $b$ such that $b|A|_\infty < \frac{1}{2}\mathrm{inj}(M)$).
	\end{proof}
	We also mention the following result, that will not be used in the remainder of the article since we are interested in positively curved manifolds and the torus is flat, but which is the first positive answer to the question of the smoothness of the displacement interpolation in the compact case.

	A map $f : \mathbb{R}^d \rightarrow \mathbb{R}^d$ is said to be additive if for all $(x,n) \in \mathbb{R}^d \times \mathbb{Z}^d, f(x+n)=f(x)+n$. A smooth additive map descends to a smooth map $\tilde{f} : \mathbb{T}^d \rightarrow \mathbb{T}^d$.

	 \begin{prop} \label{P GEOD CONV S}
		Let $\mu,\nu \in \P_\infty(\mathbb{T}^d)$ and let $T$ be the optimal transport map. Then, the flow $T_t$ defined by $T_t = (1-t) \mathrm{id} + tT$ for $t \in [0,1]$ is a smooth flow of diffeomorphisms. As a consequence, the displacement interpolation $(\mu_t)_{t \geq 0} = (T_t)_*\mu$ belongs to $\mathscr{R}(\mu,\nu)$.
	\end{prop}

	\begin{proof}

		Let $\mu,\nu \in \P_\infty(\mathbb{T}^d)$. By Theorem \cite[Theorem 1]{CorderoErausquin1999}, there exists a unique smooth diffeomorphism which is also an additive map $T: \mathbb{R}^d \rightarrow \mathbb{R}^d$ such that $T(x)= x + \partial \phi(x)$ for some $\phi$ and $\tilde{T} : \mathbb{T}^d \rightarrow \mathbb{T}^d$ is the unique optimal transport map from $\mu$ to $\nu$. Moreover, $T$ is an Aleksandrov solution to a Monge-Ampere equation. In particular (see \cite[Sections 2, 3, Theorem 5.5]{DePhilippis2013}), there exist $c_1,c_2>0$ such that:

		$$
		c_1 \mathrm{Id} \leq D T \leq c_2 \mathrm{Id}.
		$$
		Let $T_t(x) =(1-t)x + t T(x)$ for $t \in [0,1]$ (note that $T_t$ is also additive), the following holds:

		$$
		\partial T_t(x) = t D T(x) + (1-t) \mathrm{Id} \geq 1 - t(1-c_1) \mathrm{Id}, \quad \forall (x,t) \in \mathbb{S}^1 \times (0,1).
		$$
		Thus, $T_t$ is a smooth immersion for all $t\in [0,1]$. Moreover, the Mather's shortening principle (see \cite[Chapter 8]{villani2009optimal}) together with the fact that $T=T_1$ is a smooth diffeomorphism ensure that $T_t$ is injective for all $t \in [0,1]$. In particular, $\tilde{T}_t$ is a one-to-one smooth immersion from a compact manifold to itself, and thus a smooth diffeomorphism. Consequently, $(T_t)_{0 \leq t \leq 1}$ defines a smooth flow of smooth diffeomorphisms, the conclusion follows.
	\end{proof}	

	However, the displacement interpolation can still be approximated in a suitable manner by a regular path, known as the entropic interpolation.

	\noindent $\varepsilon$-\textbf{Entropic interpolation.} Let $\mathscr{X}$ be a measured space. Let $B$ be a positive measure on $\mathscr{X}$ and $P$ be a probability measure on $\mathscr{X}$. The relative entropy is defined to be

	$$
	H(P|B) = \begin{cases}
		\int_{\mathscr{X}} \log\left( \frac{\d P}{\d B} \right) \d P \quad \text{whenever } P << B, \\
		\infty \quad  \text{elsewhere}.
		\end{cases}
	$$

		For a good introduction to the Schrödinger problem and its links with optimal transport, the reader is referred to \cite{Leonard2014}. Let $B^\varepsilon$ denote the law of the reversible $\varepsilon$-slowed down Brownian motion $(B_{a\varepsilon})_{0 \leq a \leq 1}$ on $\C([0,1],M)$. Given $\mu_0,\mu_1 \in \P(M)$, the associated dynamical Schrödinger problem is about minimizing

		\begin{equation} \label{Sdyn} \tag{$\mathrm{S}_{\mathrm{dyn}}^\varepsilon$}
			\inf \left\{ H(P|B^\varepsilon), \quad P \in \P(\C([0,1],M)), \, P_0 = \mu_0, \,P_1=\mu_1\right\}.
		\end{equation}
		
		\begin{prop}
			Assume that $\mu_0,\mu_1 \in \P_\infty(M)$. Then, there exists a unique solution $(P_a)_{0 \leq a \leq 1}$ to \eqref{Sdyn} and two positive maps $f^\varepsilon, g^\varepsilon \in \C^\infty([0,1] \times M, \mathbb{R})$ such that

		$$
		\d P_t = f^\varepsilon_t g^\varepsilon_t \d \vol.
		$$
		\end{prop}
		\begin{proof}
			The positivity property can be found in \cite[Theorem 2.8]{Leonard2014} and the regularity property in \cite[Proposition 2.7]{GentilIvanLeonardRipaniTamini2020}
		\end{proof}
		Let $\varphi^\varepsilon, \, \psi^\varepsilon$ be the so-called $\varepsilon$-entropic potentials defined as:

		$$
		\varphi^\varepsilon = \varepsilon \log f^\varepsilon \quad \text{and} \quad \psi^\varepsilon = \varepsilon \log g^\varepsilon. 
		$$
		These potentials obviously belong to $\C^\infty([0,1] \times M, \mathbb{R})$ and are respectively solutions to the following forward (resp. backward) Hamilton-Jacobi-Bellman equations:

		$$
		\partial_a \varphi_a^\varepsilon = \frac{1}{2} \left( |\nabla \varphi^\varepsilon_a|^2 + \varepsilon \Delta \varphi^\varepsilon_a \right), \quad -\partial_a \psi_a^\varepsilon = \frac{1}{2} \left( |\nabla \psi^\varepsilon_a|^2 + \varepsilon \Delta \psi^\varepsilon_a \right).
		$$
		Moreover, by letting $\upsilon^\varepsilon_a = (\psi^\varepsilon_a - \varphi^\varepsilon_a)/2$, the entropic interpolation is solution to the following continuity equation

		$$
		\partial_a \mu^\varepsilon_a + \mathrm{div}(\nabla \upsilon^\varepsilon_a \mu^\varepsilon_a).
		$$
		Consequently, the solution to
		$$
		\partial_a U_a^\varepsilon = \nabla \upsilon^\varepsilon_a(U_a^\varepsilon), \quad U_0^\varepsilon = \mathrm{id}
		$$
		is a smooth path of diffeomorphisms. In particular, the $\varepsilon$-entropic interpolation is a regular path on $\P_\infty$ and is given by $\mu_a^\varepsilon= (U_a^\varepsilon)_* \mu_0$. Note that, by \cite[Lemma 3.3]{GentilIvanLeonardRipaniTamini2020},

		\begin{equation*} 
			\lim_{\varepsilon \to 0} \int_0^1 |\nabla \upsilon^\varepsilon_a|_{\mu^\varepsilon_a}^2 \d a = W_2^2(\mu_0,\mu_1).
		\end{equation*}
		We gather the results discussed above in the following proposition:

	\begin{prop}\label{Convergence entropic interpolation}
		The $\varepsilon$-entropic interpolation $(\mu^\varepsilon_a)_{0 \leq a \leq 1}$ belong to $\mathscr{R}(\mu,\nu)$ and converges toward the displacement interpolation $(\mu_a)_{0 \leq a \leq 1}$ for the $\C^0([0,1],\P)$ topology. Moreover, 

		\begin{equation*} 
			\lim_{\varepsilon \to 0} \int_0^1 |\partial_a \mu^\varepsilon_a|_{\mu^\varepsilon_a}^2 \d a = W_2^2(\mu_0,\mu_1).
		\end{equation*}
	\end{prop}

	 \noindent \textbf{Stochastic regular paths.}  Let $\mathscr{G} \subset \mathscr{F}$ be a sub-$\sigma$-algebra of $\mathscr{F}$. A mapping $U: (\omega,a,x) \in \Omega \times [0,1] \times M \mapsto U_a(x,\omega) \in M$ is said to be a $\mathscr{G}$-measurable $\C^k$ stochastic path of diffeomorphisms if \begin{enumerate}
		\item For all $(a,x) \in [0,1] \times M$, $\omega \in \Omega \mapsto U_a(x,\omega) \in M$ is $\mathscr{G}$-measurable.
	 \item There exists $\Omega_1 \subset \Omega$ with $\mathbb{P}(\Omega_1)=1$ such that for all $\omega \in \Omega_1$, $U(\cdot,\omega)$ belong to $\C^{k,\infty}([0,1] \times M,M)$, and if for all $a \in [0,1], U_a(\omega,\cdot)$ is invertible with $U^{-1}(\omega,\cdot) \in \C^{k,\infty}([0,1] \times M,M)$.
	 \end{enumerate} 

	\begin{defi}
		Let $\mu, \nu$ be two $\mathscr{G}$-measurable random variables on $\P_\infty$. A $\mathscr{G}$-measurable random variable on $\C([0,1],\P_\infty)$ is said to be a \textit{stochastic regular path from $\mu$ to $\nu$} if there exists a $\C^1$ stochastic path of diffeomorphisms $(\varphi_a)_{0 \leq a \leq 1}$ such that for almost all $\omega \in \Omega$, $(\mu_a)_{0 \leq a \leq 1}(\omega)=((\varphi_a(\omega))_*\mu)_{0 \leq a \leq 1}$.
		The set of stochastic regular paths from $\mu$ to $\nu$ is denoted by $\mathscr{R}_{\mathrm{s}}(\mu,\nu)$.
	\end{defi}

	\begin{prop}
	Let $\mu,\nu$ be two random variables on $\P$. Then,

	$$
	\mathbb{E}[W_2^2(\mu,\nu)] \leq \inf_{(\mu_a)_{0 \leq a \leq 1} \in \mathscr{R}_s (\mu,\nu)} \mathbb{E} \left[ \int_0^1 |\partial_a \mu_a|^2_{\mu_a} \d a \right],
	$$
	where the right-hand side quantity is $+\infty$ when $\mathscr{R}_s(\mu,\nu)= \emptyset$.
	\end{prop}

	\begin{proof}
		Let $(\mu_a)_{0 \leq a \leq 1} \in \mathscr{R}_{s}(\mu,\nu)$. Then, for almost all $\omega \in \Omega$, the Benamou-Brenier formulation yields
		
		$$
		W_2^2(\mu(\omega),\nu(\omega)) \leq \inf_{(m_a) \in \mathscr{R}(\mu(\omega),\nu(\omega))} \int_0^1 |\partial_a m_a|^2 \d a \leq \int_0^1 |\partial_a \mu_a(\omega)|^2_{\mu_a(\omega)} \d a.
		$$
		The conclusion follows by taking the expectation on both sides and by taking the infimum on the right-hand side.
	\end{proof}

	 \section{Infinite-dimensional stochastic calculus}\label{S 3}

	 	Let $\bar{Z}_0, \bar{Z}_1, \dots, \bar{Z}_N \in \Gamma^2(T\P)$. We consider the following Itô's SDE:
	
	\begin{equation} \label{EDS principale Ito}
	 \dw^{\grad} \M_t(\mu) = \SN \bar{Z}_i(\M_t(\mu)) \d W^i_t + \bar{Z}_0(\M_t(\mu)) \d t, \quad \M_0(\mu) = \mu \in \P_\infty,
	\end{equation}
	which is equivalent to the following Stratonovich's SDE (see \cite[Section 4.2]{Martin2025}):
	
	\begin{equation} \label{EDS principale Strato}
	\circ \dw \M_t(\mu) = \SN \bar{Z}_i(\M_t(\mu)) \circ \d W^i_t + \bar{Z}_0'(\M_t(\mu)) \d t, \quad \M_0(\mu) = \mu \in \P_\infty,
	\end{equation}
	where $\bar{Z}_0' = \bar{Z}_0 - \frac{1}{2} \SN \grad_{\bar{Z}_i} \bar{Z}_i$.

	In this section we will recall some results about SDEs on $\P$ that can be found in \cite{Martin2025}, \cite{wang2020}. After that, we will investigate the regularity properties of the solution flow.

	\subsection{SDEs on \texorpdfstring{$\P$}{P}} \label{S 3.1}
	
	\noindent \textbf{ID solutions.} A stochastic process $(\mathcal{M}_t(\mu))_{t \geq 0}$ on $\P$ is said to be a Image-dependent (ID for short) solution to Equation \eqref{EDS principale Ito} if for all $\mu \in \P$ there exists a family of adapted processes $\{ x \in M \mapsto (\varphi_t(x,\mu))_{t \geq 0} \}$ such that:
	
	\begin{enumerate}
	\item The map $t \mapsto \varphi_t(x,\mu)$ is continuous for all $x \in M$ and $\mu \in \P$.
	\item The map $x \mapsto \varphi_t(x,\mu)$ is measurable for all $t \in \mathbb{R}_+$ and $\mu \in \P$ and $\M_t(\mu) = (\varphi_t(\cdot,\mu))_* \mu$.
	\item For all $f \in \C^\infty(M)$, $x\in M$ and $\mu \in \P$,
	\begin{align*}
	f(\varphi_t(x,\mu))-f(x) &= \SN \int_0^t Z_i(\cdot, \M_s(\mu)) f(\varphi_s(x,\mu)) \d W^i_s + \int_0^t Z_0(\cdot,\M_s(\mu)) f(\varphi_s(x,\mu)) \d s \\
	& + \frac{1}{2} \SN \int_0^t \he f(Z_i(\varphi_s(x,\mu),\M_s(\mu)),Z_i(\varphi_s(x,\mu),\M_s(\mu))) \d s, \quad \forall t \geq 0.
	\end{align*}
	\end{enumerate}
	As a consequence, $(\M_t(\mu))_{t \geq 0}$ satisfies an Itô formula since we obtained that for all $f \in \C^\infty(M)$ and for all $\mu \in \P$ we have:
	
	\begin{align*}
	F_f(\M_t(\mu)) -F_f(\mu) &= \SN \int_0^t \bar{\mathcal{L}}_{\bar{Z}_i}F_f(\M_s(\mu)) \d W^i_s + \int_0^t \bar{\mathcal{L}}_{\bar{Z}_0} F_f(\M_s(\mu)) \d s \\
	& + \frac{1}{2} \SN \int_0^t \bar{\he} F_f(\bar{Z}_i(\M_s(\mu)), \bar{Z}_i(\M_s(\mu))) \d s, \quad \forall t \geq 0.
	\end{align*}
	We shall refer to the flow $(\varphi_t)_{t \geq 0}$ as the ID flow.
	
	\noindent \textbf{Well-posedness.} 

	\begin{thm}\label{Existence Unicité}
		Let $A_0, A_1, \dots, A_N \in \C^1_c(\mathbb{R}^d \times \P(M), \mathbb{R}^d)$. Then, there exists a unique solution to the following Image-dependent SDE:

		\[
		\begin{cases}
	        & dX_t(x) = \SN A_i(X_t(x), \M_t(\mu)) \d W^i_t + A_0(X_t(x), \M_t(\mu)) \d t, \\
	        & X_0(x) =x \in \mathbb{R}^d, \, \M_t(\mu) = (X_t|_M)_* \mu. \nonumber
	    \end{cases} \]
	\end{thm}

	\begin{proof}
		The proof can be performed in a similar way than \cite[Theorem 4.11]{Martin2025}. See also \cite{wang2020} for the original results on $\mathbb{R}^d \times \P_2(\mathbb{R}^d)$.
	\end{proof}

	\begin{rem}
		The hypothesis of the theorem can be weakened in the following way: $A_0, A_1, \dots, A_N \in \C^0_c(\mathbb{R}^d \times \P(M), \mathbb{R}^d)$ and there exists $K >0$ such that:

		$$
		\SN |A_i(x,\mu)-A_i(y,\nu)|^2 + 2\left\langle A_0(x,\mu)-A_0(y,\nu), x-y \right\rangle \leq K \left[ |x-y|^2 + W_2^2(\mu,\nu) \right].
		$$
	\end{rem}
	
\begin{coro}[\cite{Martin2025}, Theorem 4.11] 
	Let $\bar{Z}_0, \bar{Z}_1, \dots, \bar{Z}_N \in \Gamma^1(T\P)$. The SDE
	
	\begin{equation}\label{Brave EDS sur Wasserstein}
	\dw^{\grad} \M_t(\mu) = \SN \bar{Z}_i(\M_t(\mu)) \d W^i_t + \bar{Z}_0(\M_t(\mu)) \d t, \quad \M_0(\mu) = \mu \in \P,
	\end{equation}
	admits a unique ID solution. Moreover, if $\mu \in \P_\infty$, then $\varphi_t(\cdot,\mu) \in \di$ and $(\M_t(\mu))_{t \geq 0}$ belongs to $\P_\infty$. Furthermore, for all $p \geq 1$ there exists an increasing function $C^1_p : \mathbb{R}_+ \rightarrow \mathbb{R}_+$ such that:

	\begin{equation} \label{Estimee lip sup flot}
		\mathbb{E}\left[\sup_{t \in [0,T]}\left(d^{2p}(\varphi_t(x,\mu), \varphi_t(y,\nu))+W_2^{2p}(\M_t(\mu),\M_t(\nu)) \right) \right] \leq C^1_p(T) \left( d^{2p}(x,y) + W_2^{2p}(\mu,\nu) \right),
	\end{equation}
	for all $T >0$.

\end{coro}

Throughout the remaining of the article we consider $\P(\P)$ and $\P(\P_\infty)$ endowed with the $\bm W_2$ distance, namely,

	$$
	\bm W_2^2(\vartheta_1, \vartheta_2) = \inf \left\{ \int_{\P \times \P} W_2^2(\mu,\nu) \d \Pi(\mu,\nu), \, \Pi \in \C(\vartheta_1,\vartheta_2) \right\}, \quad \forall \vartheta_1, \vartheta_2 \in \P(\P).
	$$

\noindent \textbf{Transition semigroup.} The \textit{transition semigroup} associated to $(\M_t)_{t \geq 0}$ is the family of linear maps $(\se_t)_{t \geq 0}$ with $\se_t : L^\infty(\P) \rightarrow L^\infty(\P)$ for all $t \geq 0$ such that:

\begin{equation*}
	\se_t F(\cdot) = \mathbb{E}[ F(\M_t(\cdot))], \quad \forall F \in L^\infty(\P).
\end{equation*}
Note that $\se_t(L^\infty(\P)) \subset L^\infty(\P)$. We define in the same way $\check{\mathrm{P}}_t$, the transition semigroup associated to $(\varphi_t,\M_t)$, namely,

$$
\se_t \check{F}(x,\mu) = \mathbb{E} \left[ \check{F}(\varphi_t(x,\mu),\M_t(\mu)) \right], \quad x \in M, \, \mu \in \P, \, F \in L^\infty(M \times \P).
$$

\begin{prop} \label{Semigroup property}
	The process $(\M_t)_{t \geq 0}$ (resp. $(\varphi_t,\M_t))_{t \geq 0}$ is a strong Markov process. As a consequence, the map $\se_t : \vartheta \in \P(\P) \mapsto \vartheta \se_t \in \P(\P)$ (resp. $\check{\mathrm{P}}_t : \check{\vartheta} \in \P(M \times \P) \mapsto \check{\vartheta}\check{\mathrm{P}}_t$) is continuous for the $\bm W_2$ convergence.
\end{prop}
\begin{proof}
	We prove the proposition only for $\se_t$ since the proof is exactly the same for $\check{\mathrm{P}}_t$. The fact that $(\M_t)_{t \geq 0}$ is Markov follows from the well-posedness of the associated Image dependent SDE. 
	
	Let $F\in \C^0(\P)$, we show that $\se_t F \in \C^0(\P)$ for all $t \geq 0$. Let $\varepsilon >0$. By uniform continuity of $F$, there exists $\delta >0$ such that if $W_2^2(\mu',\nu') < \delta$, $|F(\mu')-F(\nu')| < \varepsilon$. Moreover, by Markov inequality and \eqref{Estimee lip sup flot},

	$$
	\mathbb{P}(W_2^2(\M_t(\mu),\M_t(\nu)) \geq \delta) \leq \frac{ C_t }{\delta}W_2^2(\mu,\nu),
	$$
	where $C_t$ is the constant in \eqref{Estimee lip sup flot}. Let $A = \{ \omega \, : \, W_2^2(\M_t(\mu),\M_t(\nu)) < \delta \}$ and let
	
	$$
	\delta' = \frac{\varepsilon \delta }{2 |F|_\infty C_t}
	$$
	Then, for $\mu,\nu \in \P_\infty$ such that $W_2^2(\mu,\nu)< \delta'$,

	\begin{align} \label{Inegal markov fort}
		|\se_t F(\mu) - \se_t F(\nu)| &\leq \mathbb{E} \left[ |F(\M_t(\mu))-F(\M_t(\nu))| \right] \\
		& = \mathbb{E} \left[ |F(\M_t(\mu))-F(\M_t(\nu)) | \bm 1_A \right] + \mathbb{E} \left[ |F(\M_t(\mu))-F(\M_t(\nu)) | \bm 1_{\Omega \setminus A} \right] \nonumber \\
		& < \varepsilon + 2 |F|_\infty \mathbb{P}(W_2^2(\M_t(\mu),\M_t(\nu)) \geq \delta) \nonumber \\
		& < 2\varepsilon. \nonumber
	\end{align}
	Thus, $\se_t F \in \C^0(\P)$.
	
	We now show that $\se_t F \underset{t \to 0}{\overset{\C^0(\P)}{\to}} F$. Let $A' = \{ \omega \, : \, W_2^2(\M_t(\mu),\mu) < \delta \}$, then,

	\begin{align*}
		|\se_t F(\mu) - F(\mu)| &\leq \mathbb{E} \left[ |F(\M_t(\mu)) -F(\mu)|\right] \\
		& \leq \mathbb{E} \left[ |F(\M_t(\mu))-F(\mu) | \bm 1_{A'} \right] + \mathbb{E} \left[ |F(\M_t(\mu))-F(\mu) | \bm 1_{\Omega \setminus A'} \right] \\
		& \leq \varepsilon + 2 |F|_\infty \mathbb{P}(W_2^2(\M_t(\mu),\mu) \geq \delta),
	\end{align*}
	Moreover, there exists a constant $C>0$ such that

	$$
	\mathbb{E} \left[ W_2^2(\M_t(\mu),\mu) \right] \leq \mathbb{E} \left[ \int_M |\varphi_t(x,\mu) - x|^2 \d \mu(x) \right] \leq C(t+t^2).
	$$ 
	We conclude by using once again Markov inequality.
	
	Let $F \in \C^0(\P)$ and $(\vartheta_n)_n$ a sequence in $\P(\P)$ converging toward $\vartheta \in \P(\P)$ for the $\bm W_2$ convergence. Because of the compactness of $\P$, the $\bm W_2$ convergence is equivalent to the weak convergence in $\P(\P)$, thus, since $\se_t F \in \C^0(\P)$,

	\begin{align*}
		(\vartheta_n \se_t)F &= \int_{\P} \se_t F(\mu) \d \vartheta_n(\mu) \underset{n \to \infty}{\rightarrow} \int_{\P} \se_t F \d \vartheta(\mu) = (\vartheta \se_t) F.
	\end{align*}
	\end{proof} 

\noindent \textbf{General Itô formulae.} The proof of the following result can be found in Appendix \ref{Appendix A}.

\begin{prop} \label{Ito formula}
		Let $F \in \C^2(\P)$, then the following holds:

		\begin{align*}
	F(\M_t(\mu)) -F(\mu) = \SN \int_0^t \bar{\mathcal{L}}_{\bar{Z}_i}F(\M_s(\mu)) \d W^i_s + \int_0^t \underbrace{\left[ \bar{\mathcal{L}}_{\bar{Z}_0} F(\M_s(\mu))  +  \frac{1}{2} \SN  \bar{\he} F(\bar{Z}_i(\M_s(\mu)), \bar{Z}_i(\M_s(\mu))) \right]}_{\bar{\mathrm{L}}F(\M_s(\mu))} \d s. \nonumber
	\end{align*}
	\end{prop}

	\noindent \textbf{Infinitesimal generator.} The \textit{infinitesimal generator} of $(\se_t)_{t \geq 0}$ is defined to be:

 \begin{equation*}
 \bar{\mathrm{L}} F = \left. \frac{\d}{\d t} \right|_{t=0} \se_t F,
 \end{equation*}
 for all $F$ such that this limit exists in $\C^0(\P)$. We denote the set of all the functions such that this limit exists by $\mathrm{D}(\bar{\mathrm{L}})$.
 \begin{prop} \label{Proposition dérivée semi groupe}
	Let $ F \in \mathrm{D}(\bar{\mathrm{L}})$. Then, for all $t \geq 0$, $\se_t F \in \mathrm{D}(\bar{\mathrm{L}})$ and:
	
	\begin{equation} \label{E Dérivée semi groupe}
	\frac{\d}{\d t} \se_t F = \se_t \left(\bar{\mathrm{L}} F \right) = \bar{\mathrm{L}}\left(\se_t F \right).
	\end{equation}
\end{prop}
\begin{proof}
	We have by using the semigroup property:
	
	$$
	\sup_{\mu \in \P} \left| \frac{\se_{t+\varepsilon}F - \se_t F}{\varepsilon} - \se_t \bar{\mathrm{L}}F \right| = \sup_{\mu \in \P} \left| \se_t\left( \frac{\se_\varepsilon F - F }{\varepsilon} - \bar{\mathrm{L}} F \right) \right| \leq \sup_{\mu \in \P} \left| \frac{\se_\varepsilon F - F }{\varepsilon} - \bar{\mathrm{L}} F \right|.
	$$
	Taking the limit as $\varepsilon \to 0$ yields
	
	$$
	\frac{\d}{\d t}\se_t F = \se_t \left(\bar{\mathrm{L}} F \right).
	$$
\end{proof}

 Note that, by Proposition \ref{Ito formula}, $\C^2(\P) \subset \mathrm{D}(\bar{L})$, and that 

	\begin{equation} \label{Generateur P}
		\bar{\mathrm{L}} F(\mu) = \bar{\mathcal{L}}_{\bar{Z}_0} F(\mu) + \frac{1}{2} \SN \bar{\he} F(\bar{Z}_i(\mu),\bar{Z}_i(\mu)), \quad \forall F \in \C^2(\P), \, \mu \in \P.
	\end{equation}
	In particular, $\bar{\mathrm{L}}(\C^2(\P)) \subset \C^0(\P)$.

\begin{prop} 
	Let $(x,\mu) \in M \times \P$, denote by $(X_t,\mu_t)=(\varphi_t(x,\mu),\M_t(\mu))_{t \geq 0}$ be the couple of the ID flow and the ID solution associated to \eqref{EDS principale Ito}. Let $\check{F} \in \check{\mathscr{A}} \in \C^{2,2} (M \times \P)$. Then,
	
\begin{align*}
	\check{F}(X_t,\mu_t) -\check{F}(x,\mu) &=  \int_0^t \SN \bar{\mathcal{L}}_{\bar{Z}_i}\check{F}(X_s,\cdot)(\mu_s) + \mathcal{L}_{Z_i}\check{F}(\cdot,\mu_s)(X_s) \d W^i_s + \int_0^t \bar{\mathrm{L}} \check{F}(X_s,\cdot)(\mu_s) \d s \\
	& + \int_0^t \mathcal{L}_{Z_0}\check{F}(\cdot,\mu_s)(X_s) + \frac{1}{2} \SN \mathrm{H}\check{F}(\cdot,\mu_s)(X_s)(Z_i,Z_i) \d s \\
	& + \frac{1}{2}\int_0^t \SN \bar{\mathcal{L}}_{\bar{Z}_i(\mu_s)} \mathcal{L}_{Z_i(X_s,\cdot)} F(X_s,\cdot)(\mu_s) + \mathcal{L}_{Z_i(X_s,\mu_s)} \bar{\mathcal{L}}_{\bar{Z}_i(\mu_s)}F(\cdot,\mu_s)(X_s) \d s 
	\end{align*}
\end{prop}

	We define in the same way the generator $\check{\mathrm{L}}$ associated to $(\bar{\mathrm{P}}_t)_{t \geq 0}$. We also have $\check{\mathscr{A}} \subset \mathrm{D}(\check{\mathrm{L}})$ and for any $\check{F} \in \C^{2,2}(M \times \P)$ and any $(x,\mu) \in M \times \P$ it holds:

	\begin{align} \label{Generateur MP}
		\check{\mathrm{L}} \check{F}(x,\mu) &= \bar{\mathrm{L}} \check{F}(x,\cdot)(\mu) + \mathcal{L}_{Z_i(x,\mu)} \check{F}(\cdot,\mu)(x) + \frac{1}{2} \SN \he \check{F}(\cdot,\mu)(Z_i(x,\mu),Z_i(x,\mu)) \\
		& + \frac{1}{2} \SN \bar{Z}_i \left( Z_i \check{F}(\cdot,\mu)(x) \right) + Z_i \left( \bar{Z}_i \check{F}(x,\cdot)(\mu) \right). \nonumber
	\end{align}

	\noindent \textbf{Carré du champ.} The \textit{Carré du champ} $\bar{\Gamma}$ is defined to be:

\begin{equation}
\bar{\Gamma}F = \frac{1}{2} \left( \bar{\mathrm{L}}(F^2) - 2 F \bar{\mathrm{L}}F \right), \quad \forall F \in \C^2(\P).
\end{equation}
For $F \in \C^2(\P)$, the Carré du champ also writes

$$
\left(\bar{\Gamma} F \right) (\mu) = \SN \left( \bar{\mathcal{L}}_{\bar{Z}_i} F(\mu) \right)^2 = \SN \left\langle \grad F(\mu), \bar{Z}_i(\mu) \right\rangle_\mu^2 .
$$
In particular,

$$
|\bar{\Gamma}F|_\infty \leq |\bar{Z}|^2_\infty, |\grad F|^2_\infty,
$$
where $|\bar{Z}|^2_\infty = \SN |\bar{Z}_i|^2_\infty$.

\subsection{Continuity of the ID-flow} \label{S 3.2 la vraie}

	Let us assume that $M$ is a closed manifold isometrically embedded in $\mathbb{R}^d$ (This is always possible by the Nash embedding theorem, see \cite{Nash1956TheIP}). In this case, the following holds:  

	\begin{lem} \label{L ED}
		There exists a constant $\kappa >0$ such that 

		$$
		|x-y|^2 \leq d^2_M(x,y) \leq \kappa |x-y|^2, \quad \forall x,y \in M.
		$$
		Consequently,

		$$
		W_{2,\mathbb{R}^d}^2(\mu,\nu) \leq W_{2,M}^2(\mu,\nu) \leq \kappa W_{2,\mathbb{R}^d}^2(\mu,\nu), \quad \forall \mu,\nu \in \P(M).
		$$
	\end{lem}

 	Let $A_0, A_1, \dots, A_N \in \C^{\infty,2}(M \times \P(M), TM)$ be a family of vector fields indexed by $\P(M)$ and consider the following Image-dependent SDE on $M$:

	\begin{align}\label{E 13}
		\begin{cases}
	        & d^\nabla\varphi_t(x,\mu) = \SN A_i(\varphi_t(x,\mu), \M_t(\mu)) \d W^i_t + A_0(\varphi_t(x,\mu), \M_t(\mu)) \d t, \\
	        & \varphi_0(x,\mu) =x \in M, \, \M_t(\mu) = (\varphi_t(\cdot,\mu))_* \mu, \, \mu \in \P. 
	    \end{cases}
	\end{align}
	In \cite[Subsection 4.3]{Martin2025}, it is shown that, by letting $M_\varepsilon$ be a small enough tubular neighborhood of $M$ such that the projection $p_\varepsilon : M_\varepsilon \rightarrow M$ is smooth, one can extend the driving vector fields in a smooth way by letting 

	$$
	\tilde{A}_i(x) = \chi\left(d^2(x,M)\right) A_i(p_{\varepsilon}(x)), \quad x \in \mathbb{R}^d, \, i=1,2,\dots,N,
	$$
	and 

	$$
	\tilde{A}_0(x) = \chi\left(d^2(x,M)\right) \left( A_0(p_{\varepsilon}(x)) - \frac{1}{2} \SN \II\left(Tp_{\varepsilon}(\tilde{A}_i(x)),Tp_{\varepsilon}(\tilde{A}_i(x))\right) \right),
	$$
	where $\chi$ is a smooth bump function with $\chi(0)=1$ and support $[0; \varepsilon/2]$. In this case, it can be easily seen that the coefficients $\tilde{A}_0,\tilde{A}_1, \dots,\tilde{A}_N$ are continuous in $(x,\mu) \in \mathbb{R}^d \times \P(M)$ and that there exists $K > 0$ such that for any $t \geq 0$,

		\begin{equation} \label{A 11}
			\sum_{i=0}^N |\tilde{A}_i(x,\mu)|^2 \leq K, \, (x,\mu) \in \mathbb{R}^d \times \P(M) 
		\end{equation}
		\begin{align} \label{A 12}
			& \quad \, \, 2 \langle \tilde{A}_0(x,\mu) - \tilde{A}_0(y,\nu), x-y \rangle + \SN |\tilde{A}_i(x,\mu)- \tilde{A}_i(y,\nu)|^2 \\
			& \leq K (|x-y|^2 + W_2^2(\mu,\nu)), \quad   (x,\mu),(y,\nu) \in \mathbb{R}^d \times \P(M). \nonumber
		\end{align}
	As a consequence, from \cite[Theorem 2.1]{wang2020} and \cite[Theorem 4.11]{Martin2025}, there exists a unique solution to:
	
	\begin{align}
		\begin{cases}
	        & \d \tilde{\varphi}_t(x,\mu) = \SN \tilde{A}_i(\varphi_t(x,\mu), \M_t(\mu)) \d W^i_t + \tilde{A}_0(\varphi_t(x,\mu), \M_t(\mu)) \d t, \\
	        & \tilde{\varphi}_0(x,\mu) =x \in \mathbb{R}^d, \, \M_t(\mu) = (\tilde{\varphi}_t(\cdot,\mu))_* \mu, \, \mu \in \P, \label{E 14}
	    \end{cases}
	\end{align}
	and $(\tilde{\varphi}_t(\cdot,\mu))_{t \geq 0}|_M$ is a stochastic flow on $M$. Let $(a,x) \in [0,1] \times M \mapsto U_a(x)$ be a smooth path of diffeomorphisms with $U_0 = \mathrm{id}$. Then, $V_a = \partial_a U_a \circ U_a^{-1}$ defines a smooth vector field on $M$. Let also $\tilde{V}_a$ be a compactly supported smooth extension of $V_a$ and $\tilde{U}_a$ be the unique solution to $\partial_a \tilde{U}_a= V_a(\tilde{U}_a), U_0 = \mathrm{id}$. Let also $\mu_a = (U_a)_* \mu$. The goal of this subsection is to show that the field $(\omega,t,a,x) \mapsto \varphi_t^a(x,\omega) := \varphi_t(U_a(x),\mu_a)(\omega)$ is an $(\mathscr{F}_t)_{t \geq 0}$-adapted $\C^k$ stochastic field of diffeomorphisms on $M$ in the following sense:

	\begin{defi}
		A random field $(\omega,t,a,x) \in \Omega \times \mathbb{R}_+ \times M \times [0,1] \mapsto \varphi_t^a(x,\omega) \in M$ is said to be an \textit{$(\mathscr{F}_t)_{t\geq 0}$-adapted $\C^k$ stochastic field of diffeomorphisms} if: \begin{enumerate}
			\item There exists $\Omega_1 \subset \Omega$ with $\mathbb{P}(\Omega_1)=1$ such that for all $\omega \in \Omega_1$, $\varphi_t^a(x,\omega)$ belong to $\C^{0,k,\infty}$ and for all $(t,a) \in \mathbb{R}_+ \times [0,1]$, $\varphi_t^a(\cdot,\omega)$ is a smooth diffeomorphism. 
			\item For all $0 \leq a \leq 1$ and all $x \in M$, $(\varphi_t^a(x))_{t \geq 0}$ is an $(\mathscr{F}_t)_{t\geq 0}$-adapted process.
			\item For all $t \geq 0$ and all $x \in M$, $(\varphi_t^a(x))_{0 \leq a \leq 1}$ is adapted with respect to the constant filtration $(\mathscr{F}_t)_{0 \leq a \leq 1}$.
		\end{enumerate}
	\end{defi}
	
	In what follows, we shall omit the tilde notation and denote the vector fields by $A_i$ and the solution flow by $(\varphi_t)_{t \geq 0}$.

	\begin{prop} \label{P estim lip sur M}
		For all $p \geq 1$ there exists an increasing function $C^2_p :\mathbb{R}_+ \rightarrow \mathbb{R}_+$ such that

	\begin{equation} \label{E flot sup lipschitz}
		\mathbb{E}[|\varphi_t(x,\mu)- \varphi_s(y,\nu)|^{2p}] \leq C^2_p(T) \left( |x-y|^{2p} + W_2^{2p}(\mu,\nu) + |t-s|^{p} \right),
	\end{equation}
	for all $T >0$, $(x,y) \in M^2$, $(\mu,\nu) \in \P^2$, $(s,t) \in [0,T]^2$.

	\end{prop}

	\begin{proof}
		By \eqref{A 11} and \eqref{Estimee lip sup flot}:

		\begin{align*}
			&\quad \mathbb{E}\left[|\varphi_t(x,\mu)- \varphi_s(y,\nu)|^{2p}\right]\\
			& \leq 2^{2p-1} \left(\mathbb{E}\left[|\varphi_t(x,\mu)- \varphi_t(y,\nu)|^{2p}\right] + \mathbb{E} \left[|\varphi_t(y,\nu)- \varphi_s(y,\nu)|^{2p}\right] \right) \\
			& \leq 2^{2p-1} \left( C_p^1(T) \left( |x-y|^{2p} + W_2^{2p}(\mu,\nu) \right) + \mathbb{E} \left[ \left| \SN \int_s^t A_i(\varphi_r(y,\nu), \nu_r) \d W^i_r + \int_s^t A_0(\varphi_r(y,\nu),\nu_r) \d r \right|^{2p} \right] \right) \\
			& \leq 2^{2p-1} \left( C_p^1(T) \left( |x-y|^{2p} + W_2^{2p}(\mu,\nu) \right) + \mathbb{E} \left[ \left( \SN \int_s^t \left|A_i(\varphi_r(y,\nu), \nu_r)\right|^{2} \d r \right)^p \right] + \mathbb{E} \left[ \left| \int_s^t A_0(\varphi_r(y,\nu), \nu_r) \d r \right|^{2p} \right] \right) \\
			& \leq C_p^2(T) \left( |x-y|^{2p} + W_2^{2p}(\mu,\nu) + |t-s|^p \right),
		\end{align*}
		where $C^1_p$ is the function in \eqref{Estimee lip sup flot}. The conclusion follows.
	\end{proof}

	\begin{prop} \label{Prop continuite proc}
		Let $T> 0$, $p \geq 1$. There exists $C_p(T)>0$ such that for all $s,t \in [0,T]$, $a,b \in [0,1]$ and $x,y \in \mathbb{R}^d$,

		\begin{equation} 
			\mathbb{E}[|\varphi_t(\tilde{U}_a(x),\mu_a) - \varphi_s(\tilde{U}_{b}(y),\mu^b)|^{2p}] \leq C_p(T)\left( |x-y|^{2p} + |a-b|^{2p} + |t-s|^{p} \right).
		\end{equation}
		As a consequence, there exists a continuous version of $(t,x,a) \in \mathbb{R}_+ \times \mathbb{R}^d \times [0,1] \mapsto \varphi_t(U_a(x),\mu_a)$.		
	\end{prop}

	\begin{proof}
		Since $(a,x) \mapsto \tilde{U}_a(x)$ is smooth and the identity outside $M_{\varepsilon/2}$, there exists $c>0$ such that

		$$
		|U_a(x)-U_b(y)| \leq c (|x-y| + |a-b|).
		$$
		Using the same argument, we obtain

		$$
		W_2^2(\mu_a,\mu^b) \leq \int_M |U_a(x) - U_b(x)|^2 \d \mu(x) \leq c^2 |a-b|^2.
		$$
		The conclusion follows from Proposition \ref{P estim lip sur M} and the Kolmogorov-$\check{\text{C}}$encov Theorem.
	\end{proof}

	\begin{prop} \label{T C0}
 Let $k \geq 0$ and assume that $A_0, A_1, \dots, A_N \in \C^{\infty,1}(M \times \P(M),TM)$. Let also $(\varphi_t)_{t \geq 0}$ be the unique ID solution flow associated to \eqref{EDS principale Ito} and $(x,a) \in M \times [0,1] \mapsto U_a(x)$ be a smooth path of diffeomorphisms, $\mu_a = (U_a)_*\mu$. Then, there exists a $\C^{0,k,\infty}$ modification of

 \begin{align*}
 	\Omega \times \mathbb{R}_+ \times [0,1] \times M  & \longrightarrow M \\
 	(\omega,t,a,x) &\longmapsto \varphi_t^a(x,\omega)= \varphi_t(U_a(x),\mu_a)(\omega).
 \end{align*}
 Moreover, the field $(\omega,t,a,x) \mapsto \varphi_t^a(x,\omega)$ is an $(\mathscr{F}_t)_{t \geq 0}$-adapted $\C^0$ stochastic field of diffeomorphisms on $M$.
 \end{prop}

 \begin{proof}
	The first part of the theorem follows by iterating the arguments developed above for the spatial derivatives of the field. The second part can be obtained by showing that the Jacobian matrix of $\varphi_t(x,a)$ is nonsingular and by applying the inverse function theorem, see \cite[Lemma 4.4.3, Theorem 4.6.5]{kunita1997stochastic}.
 \end{proof}

	\subsection{Higher regularity of the flow}

	We consider the space $\P_c(M \times \mathbb{R}^d)$ endowed with the product Wasserstein-$2$ distance, also denoted by $W_2$, namely,
	
	$$
	W_2^2(\zeta,\eta) = \inf_{\pi \in \C(\zeta,\eta)} \int_{M \times \mathbb{R}^d} d^2(x,y) + |u-v|^2 \d \pi ((x,u),(y,v)), \quad \zeta,\eta \in \P_c(M \times \mathbb{R}^d).
	$$
	We denote by $\mathrm{pr}_1$ (resp. $\mathrm{pr}_2$) the projection from $M \times \mathbb{R}^d$ onto $M$ (resp. $\mathbb{R}^d$). Let $V \in \Gamma(TM)$ and consider a smooth extension to $\mathbb{R}^d$. In the sequel, we will also denote by $V$ such an extension. For a measure $\zeta \in \P_c(M \times \mathbb{R}^d)$, we denote $|\zeta|^2 = \int_{M \times \mathbb{R}^d} |u|^2_z \d \zeta(z,u)$.

	Let $(a,x) \in [0,1] \times M \mapsto U_a(x)$ be a smooth path of diffeomorphisms with $U_0 = \mathrm{id}$, $\dot{U}_0 =V \in \Gamma(TM)$. Let us explain heuristically the purpose of this subsection. Assume that $(a,x) \mapsto \varphi_t(U_a(x),\mu_a)$ is regular (in some sense) almost everywhere. Then, for any $f \in \C^\infty(M)$, since $\M_t(\mu_a) = \varphi_t(U_a(\cdot),\mu_a)_*\mu = (\varphi_t(\cdot,\mu_a) \circ U_a)_*\mu = (\varphi_t(\cdot,\mu_a))_*\mu_a$,

	\begin{align*}
			\partial_a F_f(\M_t(\mu_a)) &= \partial_a \int_M f(\varphi_t(x,\mu_a)) \d \mu_a \\
			& = \partial_a \int_M f(\varphi_t(U_a(x),\mu_a)) \d \mu \\
			& = \int_M \langle \nabla f(\varphi_t(U_a(x),\mu_a)), \partial_a \varphi_t(U_a(x),\mu_a) \rangle \d \mu \\
			& = \int_M \langle \nabla f(x), (\partial_a \varphi_t(U_a(\cdot),\mu_a)) \circ (\varphi_t(U_a(\cdot),\mu_a))^{-1} \rangle \d \M_t(\mu_a) \\
			& =\int_M \langle \nabla f(x), \Pi_{\mu^a_t} \left[ (\partial_a \varphi_t(U_a(\cdot),\mu_a)) \circ (\varphi_t(U_a(\cdot),\mu_a))^{-1} \right] \rangle \d \M_t(\mu_a).
	\end{align*}
	Thus, we want to say that

	$$
	\partial_a \M_t(\mu_a) = \Pi_{\mu^a_t} \left[ (\partial_a \varphi_t(U_a(\cdot),\mu_a)) \circ (\varphi_t(U_a(\cdot),\mu_a))^{-1} \right].
	$$
	The goal of this subsection is to make rigourous this calculation by showing that $(a,x) \mapsto X_t^a(x)$ is $\C^{1,0}([0,1] \times M,M)$ is regular enough.
	In the following we denote $X_t^a(x,\omega):= \varphi_t(U_a(x),\mu_a)(\omega)$ and $\mu_t^a = \M_t(\mu_a) $ for simplicity.

	\begin{lem}
		There exists a $\mathscr{F}_T$-measurable map $\bar{r} : \Omega \to (0,1]$ and a family of $\mathscr{F}_T$-measurable smooth vector fields $(V_t^r)_{0 \leq r \leq \bar{r}} = (\log_{X_t}(X_t^r))_{0 \leq r \leq \bar{r}}$ such that for all $r \in [0,\bar{r}]$, the map defined by $(\psi_t^{r,a})_{0 \leq a \leq 1} = (\exp(a V_t^r) \circ X_t)_{0 \leq a \leq 1}$ is a smooth path of diffeomorphisms with $\psi_t^{r,1}= X_t^r$.
	\end{lem}

	\begin{proof}
		By Proposition \ref{T C0} and compactness of $[0,1] \times [0,T] \times M$, for almost all $\omega \in \Omega$, $(a,t,x) \in [0,1] \times [0,T] \times M \rightarrow X_t^a(x,\omega)$ is uniformly continuous, and thus, $a \in [0,1] \mapsto \sup_{t \in [0,T]}\sup_{x \in M} d(X_t(x,\omega),X_t^a(x,\omega))$ is also continuous. In particular, the result follows by taking

		$$
		\bar{r}(\omega) = \sup \left\{ r \in [0,1], \, \sup_{t \in [0,T]}\sup_{x \in M} d(X_t(x,\omega),X_t^r(x,\omega)) \leq \frac{\mathrm{Inj}(M)}{2} \right\}
		$$
		and $V_t^r = \log_{X_t} (X_t^r)$. Moreover, $\bar{r}>0$ by continuity of $a \in [0,1] \mapsto \sup_{t \in [0,T]}\sup_{x \in M} d(X_t(x,\omega),X_t^a(x,\omega))$.
	\end{proof}

	Let $\varepsilon >0$ and $V_t^{\varepsilon \wedge \bar{r}} \in \Gamma(TM)$ with $(\psi_t^a)_{0 \leq a \leq 1} = (\exp(a V_t^{\varepsilon \wedge \bar{r}}) \circ X_t)_{0 \leq a \leq 1}$ a smooth path of diffeomorphisms and $\psi_t^{1} = X_t^{\varepsilon \wedge \bar{r}}$.

	\begin{lem}
		Let $p \geq 1$. For any $F \in \C^{1,1}(M \times \P)$, there exists a constant $C_p>0$ depending on $F$ such that the process

		$$
		\Theta_t^{\varepsilon\wedge \bar{r}}(F)(x) := \frac{F(X_t^{\varepsilon\wedge \bar{r}},\mu_t^{\varepsilon\wedge \bar{r}})-F(X_t,\mu_t)}{\varepsilon\wedge \bar{r}} - \left\langle \nabla F(\cdot,\mu_t)(X_t(x)),\frac{V_t^{\varepsilon\wedge \bar{r}} \circ X_t(x)}{\varepsilon\wedge \bar{r}} \right\rangle - \left\langle \grad F(X_t,\cdot)(\mu_t),\frac{ V_t^{\varepsilon\wedge \bar{r}} \circ X_t}{\varepsilon\wedge \bar{r}} \right\rangle_{\mu}
		$$
		satisfies

		\begin{equation} \label{Equation 1 THETA}
			|\Theta_t^{\varepsilon\wedge \bar{r}}(F)(x)|^{2p} \leq C \left( \left|\frac{V_t^{\varepsilon\wedge \bar{r}} \circ X_t(x)}{\varepsilon\wedge \bar{r}}\right|^{2p}_{X_t(x)} + \left|\frac{V_t^{\varepsilon\wedge \bar{r}} \circ X_t}{\varepsilon\wedge \bar{r}}\right|_{\mu}^{2p} \right), \quad \forall t \in [0,T].
		\end{equation}
		Moreover, for all $p \geq 1$, there exists an increasing function $C_p : \mathbb{R}_+ \to \mathbb{R}_+$ such that 
		
		\begin{equation}\label{Equation 2 THETA}
			\mathbb{E} \left[ \sup_{t \in [0,T]}\left( \left|\frac{V_t^{\varepsilon\wedge \bar{r}} \circ X_t(x)}{\varepsilon\wedge \bar{r}}\right|^{2p}_{X_t(x)} + \left|\frac{V_t^{\varepsilon\wedge \bar{r}} \circ X_t}{\varepsilon\wedge \bar{r}}\right|_{\mu}^{2p} \right) \right] \leq C_p(T) \left( |V(x)|^{2p} + \int_M |V(y)|^{2p} \d \mu(y)\right) + o_{\varepsilon \to 0}(1)  
		\end{equation}
		and

		\begin{equation} \label{Equation 3 THETA}
			\lim_{\varepsilon \to 0} \mathbb{E} \left[|\Theta_t^{\varepsilon\wedge \bar{r}}(F)(x)|^{2p}\right] = 0, \quad \forall t \in [0,T].
		\end{equation}
	\end{lem}

	\begin{proof}
		 Let $\nu_t^a = (\psi_t^a)_* \mu$ for $a \in [0,1]$. By Lemma \ref{Lemme fondamental L}:

		\begin{align*}
			F(X_t^{\varepsilon\wedge \bar{r}},\mu_t^{\varepsilon\wedge \bar{r}})-F(X_t,\mu_t) &= \int_0^{1} \partial_a F(\psi_t^a,\nu_t^a) \d a \\
			& = \int_0^1 \langle \grad F(\psi_t^a,\cdot)(\nu_t^a)(\psi_t^a), \partial_a \psi_t^a \rangle_{\mu} \d a + \int_0^1 \langle \nabla F(\cdot,\nu_t^a)(\psi_t^a(x)), \partial_a \psi_t^a(x) \d a\rangle  
		\end{align*}
		Thus, since for all $x \in M$ and all $a \in [0,1]$, $|\partial_a \psi_t^a(x)|_{\psi_t^a(x)} = |V_t^{\varepsilon\wedge \bar{r}}(X_t(x))|_{X_t(x)}$,

		\begin{align} \label{Equation 4 THETA}
			& \quad |\Theta_t^{\varepsilon\wedge \bar{r}}(F)(x)|^2 \\
			&\leq 2 \left| \frac{1}{\varepsilon\wedge \bar{r}} \int_0^1 \langle \nabla F(\cdot,\nu_t^a)(\psi_t^a(x)), \partial_a \psi_t^a \rangle - \langle \nabla F(\cdot,\mu_t)(X_t), V_t^{\varepsilon\wedge \bar{r}} \circ X_t(x) \rangle \d a \right|^2 \nonumber \\
			& + 2\left|  \frac{1}{\varepsilon\wedge \bar{r}} \int_0^1 \langle \grad F (\psi_t^a(x),\cdot)(\nu_t^a)(\psi_t^a), \partial_a \psi_t^a \rangle + \langle \grad F(X_t,\cdot)(\mu_t), V_t^{\varepsilon\wedge \bar{r}} \circ X_t \rangle_{\mu} \d a \right|^2 \nonumber \\
			& \leq 8 \left[ \sup_{(x',\mu') \in M \times \P} |\nabla F(x',\mu')|_x^2 \cdot \left| \frac{V_t^{\varepsilon\wedge \bar{r}}(X_t(x))}{\varepsilon\wedge \bar{r}} \right|^2 +\sup_{(x',y,\mu') \in M \times \P} |\grad F(x',\mu')(y)|^2 \cdot \left|\frac{V_t^{\varepsilon\wedge \bar{r}} \circ X_t}{\varepsilon\wedge \bar{r}}\right|_\mu^2 \right] \nonumber
		\end{align}
		and \eqref{Equation 1 THETA} follows from the convexity of $r \mapsto r^p$ for $p \geq 1$. To prove \eqref{Equation 2 THETA}, apply \eqref{Estimee lip sup flot} and use the fact $|V_t^{\varepsilon\wedge \bar{r}}(X_t(x))| = d(X_t(x),X_t^{\varepsilon\wedge \bar{r}}(x))$ to obtain:

		\begin{align*}
			& \quad \mathbb{E} \left[ \sup_{t \in [0,T]}\left( \left|\frac{V_t^{\varepsilon\wedge \bar{r}} \circ X_t(x)}{\varepsilon\wedge \bar{r}}\right|^{2p}_{X_t(x)} + \left|\frac{V_t^{\varepsilon\wedge \bar{r}} \circ X_t}{\varepsilon\wedge \bar{r}}\right|_{\mu}^{2p} \right) \right] \\
			&\leq\mathbb{E} \left[  \sup_{t \in [0,T]} \left( \frac{d^{2p}(X_t(x),X_t^{\varepsilon\wedge \bar{r}}(x))}{(\varepsilon\wedge \bar{r})^{2p}} + \int_M \frac{d^{2p}(X_t(y),X_t^{\varepsilon\wedge \bar{r}}(y))}{(\varepsilon\wedge \bar{r})^{2p}} \d \mu(y) \right) \right] \\
			& \leq C_p^1(T)  \left( \frac{d^{2p}(U_{\varepsilon\wedge \bar{r}}(x),U(x))}{(\varepsilon\wedge \bar{r})^{2p}} + \int_M \frac{d^{2p}(U_{\varepsilon\wedge \bar{r}}(y),U(y))}{(\varepsilon\wedge \bar{r})^{2p}} \d \mu(y) \right) \\
			& = C_p^1(T) \left( |\dot{U}_0(x)|^{2p} + \int_M |\dot{U}_0(y)|^{2p} \d \mu(y)\right) + o_{\varepsilon \to 0}(1)  ,
		\end{align*}
		where the last equality comes from $d^{2p}(U_{\varepsilon}(x),U(x))  =|\dot{U}_0(x)|^{2p} \varepsilon^{2p}+ o(\varepsilon^{2p})$ (this fact can be checked by taking an exponential chart centered in $U(x)$).
		To prove \eqref{Equation 3 THETA}, consider smooth extensions of the considered vector fields. The first term on the right-hand side of the equality in \eqref{Equation 4 THETA} satisfies

		\begin{align*}
			& \mathbb{E} \left[ \left| \frac{1}{\varepsilon\wedge \bar{r}} \int_0^1 \langle \nabla F(\cdot,\nu_t^a)(\psi_t^a(x)), \partial_a \psi_t^a \rangle - \langle \nabla F(\cdot,\mu_t)(X_t), V_t^{\varepsilon\wedge \bar{r}} \circ X_t(x) \rangle \d a \right|^2 \right] \\
			& \leq 2 \mathbb{E} \left[ \left| \frac{1}{\varepsilon\wedge \bar{r}} \int_0^1 \langle \nabla F(\cdot,\nu_t^a)(\psi_t^a(x)), \partial_a \psi_t^a(x) - V_t^{\varepsilon\wedge \bar{r}}(x) \rangle \d a\right|^2 \right] \\
			& + 2\mathbb{E} \left[ \left| \frac{1}{\varepsilon\wedge \bar{r}} \int_0^1  \langle \nabla F(\cdot,\nu_t^a)(\psi_t^a(x))-  \nabla F(\cdot,\mu_t)(X_t), V_t^{\varepsilon\wedge \bar{r}} \circ X_t(x) \rangle \d a \right|^2 \right].
		\end{align*}
		We conclude by using dominated convergence Theorem together with the fact that $\sup_{(x,\mu)}|\nabla F(x,\mu)| <\infty$, $\sup_{(x,y,\mu)}|\grad F(x,\mu)(y)|<\infty$, $\sup_{x \in M}d(\psi_t^a(x), X_t(x)) \to_{a \to 0} 0$ and $W_2(\nu_t^a,\mu_t) \to_{a \to 0} 0$. The second term on the right-hand side of the equality in \eqref{Equation 4 THETA} can be treated in the exact same way.
	\end{proof}

	\begin{lem} \label{Lemme difference vecteurs}
		Let $M$ be a closed Riemannian manifold isometrically embedded in $\mathbb{R}^d$. Let $x,y \in M$ and $v \in T_xM$ such that $\exp_x(v)=y$. Then, 

		$$
		|v-(y-x)| \leq \frac{|\II|_\infty}{2} |v|^2 = \frac{|\II|_\infty}{2} d^2(x,y).
		$$
	\end{lem}
	\begin{proof}
		Let $(\gamma_t)_{0 \leq  t\leq 1} = (\exp_x(tv))_{0 \leq t \leq 1}$. Then, $\ddot{\gamma}_t = \II(\dot{\gamma}_t,\dot{\gamma}_t)$. Thus, $\dot{\gamma}_t = v + \int_0^t \II(\dot{\gamma}_s,\dot{\gamma}_s) \d s$. Finally, since $|\dot{\gamma}_s|=|v|$ for all $s \in [0,1]$,

		$$
		|x-y -v| = \left| \int_0^1 \dot{\gamma}_t - v \d t \right| = \left| \int_0^1 \int_0^t \II(\dot{\gamma}_s,\dot{\gamma}_s) \d s \d t \right| \leq |\II|_\infty |v|^2 \int_0^1 \int_0^t \d s \d t = \frac{|\II|_\infty }{2}|v|^2.
		$$
	\end{proof}

	\begin{prop} \label{L inegalite croissance}
		Let $V \in \Gamma(TM)$. Then, the process

		$$
		w_t(x,\mu,V):= T_x\varphi_t(\cdot,\mu)(V(x)) + \T_\mu \varphi_t(x,\cdot)(V), \quad x \in M, \, \mu \in \P(M)
		$$
	exists in $L^2(\Omega, \C([0,T], \mathbb{R}^d);\mathbb{P})$ for all $T>0$. Moreover, the couple $(\varphi_t(x,\mu), w_t(x,\mu,V))$ is solution to the following Image-dependent SDE:

	\begin{equation} \label{MVC TM}
		\begin{cases}
	        & \d \varphi_t(x,\mu) = \SN A_i(\varphi_t(x,\mu), \mu_t) \d W^i_t + A_0(\varphi_t(x,\mu), \mu_t) \d t, \\
	        & \d w_t(x,\mu) = \SN B_i(\varphi_t(x,\mu), w_t(x,\mu)(V), \zeta_t) \d W^i_t + B_0(\varphi_t(x,\mu), w_t(x,\mu)(V), \zeta_t) \d t, \\
			& \varphi_0(x,\mu) = x, \,  w_0(x,\mu)(V)= V(x), \, \zeta_t= (\varphi_t(\cdot,\mu), w_t(\cdot,\mu)(V))_*\mu, \, \mu_t= (\mathrm{pr}_1)_*(\zeta_t), 
	    \end{cases}
	\end{equation}
		where, for all $i=0,1,\dots,N$,

		$$
		B_i(x,v,\zeta) = \underbrace{T_x A_i(\cdot,\mu)(v)}_{B_i^1(x,v,\mu)} + \underbrace{\int_{M \times \mathbb{R}^{d}} \left \langle \grad \tilde{A}_i(x,\mu)(y), u \right \rangle \d \zeta(y,u)}_{B_i^2(x,\zeta)},
		$$
		for all $(x,v) \in \mathbb{R}^{2d}$, $\zeta \in \P_c(M \times \mathbb{R}^d)$ and $\mu= (\mathrm{pr}_1)_* \zeta$.
		Furthermore, there exists an increasing function $C^3_p : \mathbb{R}_+ \rightarrow \mathbb{R}_+$ such that: 

			\begin{equation} \label{Continuite de la diff}
				\mathbb{E}\left[\sup_{t \in [0,T]} \left( |w_t(x,\mu,V)|^{2p} + |\zeta_t|^{2p} \right) \right] \leq C^3_p(T) \left( 1+ |V(x)|^{2p} + |V_*\mu|^{2p} \right) \leq C_p^3(T) \left( 1+ 2 |V|_\infty^{2p}\right),
			\end{equation}
			for all $T >0$.
	\end{prop}
	\begin{proof}
		Let $\tilde{V}_t^{\varepsilon\wedge \bar{r}} = X_t^{\varepsilon\wedge \bar{r}} - X_t$ for $\varepsilon >0$. We show that 

		$$
		\lim_{\varepsilon, \delta \to 0} \mathbb{E} \left[ \sup_{t \in [0,T]} \left|\frac{\tilde{V}_t^{\varepsilon\wedge \bar{r}}(x)}{\varepsilon\wedge \bar{r}}  - \frac{\tilde{V}_t^{\delta\wedge \bar{r}}(x)}{\delta\wedge \bar{r}} \right|^2 \right] =0.
		$$
		Note that, by Lemma \ref{Lemme difference vecteurs} and \eqref{Estimee lip sup flot},

		\begin{equation} \label{difference sympa}
			\mathbb{E} \left[\sup_{t \in [0,T]} \left| \frac{V_t^{\varepsilon\wedge \bar{r}}(x)}{\varepsilon\wedge \bar{r}}-\frac{\tilde{V}_t^{\varepsilon\wedge \bar{r}}(x)}{\varepsilon\wedge \bar{r}} \right|^2 \right] \leq \frac{|\II|_\infty}{2(\varepsilon\wedge \bar{r})^2} \mathbb{E}\left[ \sup_{t \in [0,T]}  d^4(X_t(x),X_t^{\varepsilon\wedge \bar{r}}(x))\right] = o(\varepsilon^2).
		\end{equation}
		We have 

		\begin{align} \label{EDS approchee}
			\frac{\tilde{V}_t^{\varepsilon\wedge \bar{r}}(x)}{\varepsilon\wedge \bar{r}} & = \frac{1}{\varepsilon\wedge \bar{r}} \SN \int_0^t A_i(X_s^{\varepsilon\wedge \bar{r}}(x),\mu_s^{\varepsilon\wedge \bar{r}}) -A_i(X_s(x),\mu_s) \d W^i_s + \frac{1}{\varepsilon\wedge \bar{r}}\int_0^t A_0(X_s^{\varepsilon\wedge \bar{r}}(x),\mu_s^{\varepsilon\wedge \bar{r}}) -A_0(X_s(x),\mu_s) \d s \\
			& = \int_0^t \Theta_s^{\varepsilon\wedge \bar{r}}(A_i)(x) +\left\langle \nabla A_i(\cdot,\mu_s)(X_s(x)),\frac{V_s^{\varepsilon\wedge \bar{r}} \circ X_s(x)}{\varepsilon\wedge \bar{r}} \right\rangle + \left\langle \grad A_i(X_s,\cdot)(\mu_s),\frac{ V_s^{\varepsilon\wedge \bar{r}} \circ X_s}{\varepsilon\wedge \bar{r}} \right\rangle_{\mu} \d W^i_s \nonumber \\
			& + \int_0^t \Theta_s^{\varepsilon\wedge \bar{r}}(A_0)(x) +\left\langle \nabla A_0(\cdot,\mu_s)(X_s(x)),\frac{V_s^{\varepsilon\wedge \bar{r}} \circ X_s(x)}{\varepsilon\wedge \bar{r}} \right\rangle + \left\langle \grad A_0(X_s,\cdot)(\mu_s),\frac{ V_s^{\varepsilon\wedge \bar{r}} \circ X_s}{\varepsilon\wedge \bar{r}} \right\rangle_{\mu} \d s\nonumber.
		\end{align}
		By Burkholder-Davis-Gundy inequality, \eqref{difference sympa} and \eqref{Equation 3 THETA}, there exists a constant $C>0$ such that:

		\begin{align}
			& \quad \mathbb{E} \left[ \sup_{t \in [0,T]} \left|\frac{\tilde{V}_t^{\varepsilon\wedge \bar{r}}}{\varepsilon}  - \frac{\tilde{V}_t^{\delta\wedge \bar{r}}}{\delta} \right|^2 \right] \nonumber\\
			& \leq C \mathbb{E} \left[  \int_0^T \sum_{i=0}^N|\Theta_t^{\varepsilon\wedge \bar{r}}(A_i)(x) - \Theta_t^{\delta\wedge \bar{r}}(A_i)(x)|^2 + \left|\frac{V_t^{\varepsilon\wedge \bar{r}} \circ X_t(x)}{\varepsilon \wedge \bar{r}} -\frac{V_t^{\delta\wedge \bar{r}} \circ X_t(x)}{\delta\wedge \bar{r}}\right|^2 + \left|\frac{V_t^{\varepsilon\wedge \bar{r}}\circ X_t}{\varepsilon \wedge \bar{r}}-\frac{V_t^{\delta\wedge \bar{r}}\circ X_t}{\delta\wedge \bar{r}} \right|^2_\mu \d t  \right]. \nonumber\\
		\end{align}
		Combining \eqref{Equation 1 THETA} and \eqref{Equation 2 THETA}, we obtain that

		$$
		\sum_{i=0}^N\int_0^T \mathbb{E} \left[ |\Theta_t^{\varepsilon\wedge \bar{r}}(A_i)(x) - \Theta_t^{\delta\wedge \bar{r}}(A_i)(x)|^2 \right] \d t= o_{\varepsilon,\delta \to 0}(1).
		$$
		Thus,

		\begin{equation}\label{Inequality}
			 \quad \mathbb{E} \left[ \sup_{t \in [0,T]} \left|\frac{\tilde{V}_t^{\varepsilon\wedge \bar{r}}}{\varepsilon}  - \frac{\tilde{V}_t^{\delta\wedge \bar{r}}}{\delta} \right|^2 \right] \leq C \mathbb{E} \left[ \int_0^T \sup_{s \in [0,t]} \left|\frac{\tilde{V}_s^{\varepsilon\wedge \bar{r}}(x)}{\varepsilon\wedge \bar{r}} -\frac{\tilde{V}_s^{\delta\wedge \bar{r}}(x)}{\delta\wedge \bar{r}}\right|^2 + \sup_{s \in [0,t]} \left|\frac{\tilde{V}_s^{\varepsilon\wedge \bar{r}}}{\varepsilon\wedge \bar{r}}-\frac{\tilde{V}_s^{\delta\wedge \bar{r}}}{\delta\wedge \bar{r}} \right|^2_\mu   \d t \right] +  o_{\varepsilon, \delta \to 0}(1).
		\end{equation}
		Integrating with respect to $\mu$ and applying Gronwall inequality yields: 

		$$
		\lim_{\varepsilon, \delta \to 0} \mathbb{E} \left[ \sup_{t \in [0,T]}\left|\frac{\tilde{V}_t^{\varepsilon\wedge \bar{r}}(x)}{\varepsilon\wedge \bar{r}}-\frac{\tilde{V}_t^{\delta\wedge \bar{r}}(x)}{\delta\wedge \bar{r}} \right|^2_\mu \right] =0.
		$$
		Injecting this result in \eqref{Inequality} and applying once again Gronwall inequality yields:

		$$
		\lim_{\varepsilon, \delta \to 0} \mathbb{E} \left[ \sup_{t \in [0,T]}\left|\frac{\tilde{V}_t^{\varepsilon\wedge \bar{r}}(x)}{\varepsilon\wedge \bar{r}}-\frac{\tilde{V}_t^{\delta\wedge \bar{r}}(x)}{\delta\wedge \bar{r}} \right|^2 \right] =0.
		$$
		Thus, for all $x \in M$, there exists $(\tilde{V}^*_t(x))_{t \in [0,T]}$ such that $\lim_{\varepsilon \to 0} (\tilde{V}_t^{\varepsilon\wedge \bar{r}}(x))_{t \in [0,T]} =(\tilde{V}^*_t(x))_{t \in [0,T]}$ in $L^2(\Omega,\C([0,T],\mathbb{R}^d);\mathbb{P})$. Since for almost all $\omega \in \Omega$, if $\varepsilon < \bar{r}(\omega)/2$ then $\varepsilon \wedge \bar{r}(\omega) = \varepsilon$, it follows that $\sup_{t \in [0,T]}|\tilde{V}_t^{\varepsilon \wedge\bar{r}}(x)/(\varepsilon \wedge\bar{r}) - \tilde{V}_t^\varepsilon(x)/\varepsilon| \overset{a.s.}{\to}_{\varepsilon \to 0} 0$. Then, by dominated convergence theorem together with \eqref{Equation 2 THETA} and \eqref{difference sympa}, it follows that 
		
		$$
		\lim_{\varepsilon \to 0}\mathbb{E} \left[ \sup_{t \in [0,T]} \left|\frac{\tilde{V}_t^{\varepsilon \wedge\bar{r}}(x)}{\varepsilon \wedge\bar{r}}-\frac{\tilde{V}_t^\varepsilon(x)}{\varepsilon}\right|^2 \right] = 0.
		$$
		Thus,
		
		$$
		\lim_{\varepsilon \to 0}\mathbb{E} \left[ \sup_{t \in [0,T]}\left| \frac{\tilde{V}_t^\varepsilon(x)}{\varepsilon}- \tilde{V}^*_t(x)\right|^2 \right] =0.
		$$
		The fact that this process is solution to \eqref{MVC TM} follows by taking the limit when $\varepsilon \to 0$ in \eqref{EDS approchee}. Finally, \eqref{Continuite de la diff} follows by taking the limit when $\varepsilon \to 0$ in \eqref{Equation 2 THETA}.

	\end{proof}

\begin{coro}
		For all $x \in M$, $\mu \in \P$, $V \in T_\mu \P$, $(X_t(x), w_t(x,\mu,V))$ belongs to $TM$ and is solution to the following covariant SDE:
		
		\begin{align*}
			\begin{array}{ll}
			D_t^\nabla w_t(x,\mu,V) &= \SN (\nabla A_i(\cdot,\mu_t))(w_t(x,\mu,V)) + \mathcal{\bar{L}}_{w_t(X_t^{-1},\mu,V)} A_i(X_t(x),\cdot)(\mu_t) \d W^i_t \\
			& +(\nabla A_0(\cdot,\mu_t))(w_t(x,\mu,V)) + \mathcal{\bar{L}}_{w_t(X_t^{-1},\mu,V)} A_0(X_t(x),\cdot)(\mu_t) \d t \\ 
			&- \SN \frac{1}{2} R(w_t(x,\mu,V),A_i(X_t(x),\mu_t))A_i(X_t(x),\mu_t) \d t.
			\end{array}
		\end{align*}
	\end{coro}

	\begin{proof}
		This results follows from Proposition \ref{E I-S Euclid} and by using the same arguments as in Corollary \ref{Coro Pro Tan}. 
	\end{proof}

		The proof of the following lemma can be found in Appendix \ref{Appendix 3}.

		\begin{lem} \label{L estim lip}
			There exists a constant $K>0$ such that for any $ (x,v,\zeta),(y,u,\eta) \in \mathbb{R}^{2d} \times \P_c(M \times \mathbb{R}^d)$ with $(\mathrm{pr}_1)_* \zeta = \mu, \, (\mathrm{pr}_1)_* \eta = \nu$, the following holds:

\begin{align}
	&\sum_{i=0}^N |B_i(x,v,\zeta)|^2 \leq K(1 + |v|^2 + |\zeta|^2) ;\\
	& \begin{array}{l}
		\langle B_0^1(x,u,\mu) - B_0^1(y,v,\nu), u-v \rangle + \SN |B_i^1(x,u,\mu)- B_i^1(y,v,\nu)|^2 \\
	 \leq K \left[ \min(|u|^2,|v|^2) \cdot ( |x-y|^2 + W_2^2(\mu,\nu))  + |u-v|^2 \right] ;
	\end{array} \\
	&\langle B_0^2(x,\zeta) - B_0^2(y,\eta), u-v \rangle + \SN |B_i^2(x,\zeta)- B_i^2(y,\eta)|^2 \leq K \min(|\zeta|^2,|\eta|^2) \left(|x-y|^2  + W_2^2(\zeta,\eta)\right)+ K |u-v|^2.
\end{align}
\end{lem}

		\begin{prop} \label{P Estim lip t cons}
			Let $V,V' \in \Gamma(TM)$. Then, for all $p \geq 1$, there exists an increasing function $C^4_p : \mathbb{R}_+ \rightarrow \mathbb{R}_+$ depending on $\min(|V|_\infty, |V'|_\infty)$ such that

			\begin{align}
			& \quad \,\, \mathbb{E}[|w_t(x,\mu,V)-w_t(y,\nu,V')|^{2p}] \label{E 133}\\
			&\leq C^4_p(T) \left( |x-y|^{2p} + |V(x) - V'(y)|^{2p} + W_2^{2p}\left[(\mathrm{id},V)_*\mu,(\mathrm{id},V')_*\nu\right] \right), \nonumber
		\end{align}
		for all $T>0$.
		\end{prop}

		\begin{proof}
			Let $X_t = \varphi_t(x,\mu)$, $Y_t = \varphi_t(y,\nu)$, $V_t = w_t(x,\mu,V)$, $V'_t = w_t(y,\nu,V')$ and $k_t = \min\left( |V_t|^2, |V_t'|^2 \right)$. Let also 

			$$
				V_t^1 = V(x) + \SN \int_0^t B_i^1(X_s, V_s, \mu_s) \d W^i_s + \int_0^t B_0^1(X_s, V_s, \mu_s) \d s, \quad \mu_s= \varphi_s(\cdot,\mu)_*\mu
			$$
			and

			$$
				V_t^{'1} = V'(y) + \SN \int_0^t B_i^1(Y_s, V_s', \nu_s) \d W^i_s + \int_0^t B_0^1(Y_s, V_s', \nu_s) \d s, \quad \nu_s= \varphi_s(\cdot,\nu)_*\nu.
			$$
			By Itô formula and Lemma \ref{L estim lip}, we obtain that

			\begin{equation} \label{E 156}
				\d |V_t^1 - V_t^{'1}|^2 \leq 2 \d m_t^1 + K \underbrace{\left( \left[ |X_t - Y_t|^2  + W_2^{2}(\mu_t,\nu_t) \right]\cdot k_t + |V_t-V_t'|^2 \right) \d t}_{\d l^1_t}, 
			\end{equation}
			with initial condition $|V-V'|^2$ and where

			$$
			\d m^1_t = \SN \langle B_i^1(X_t,V_t,\mu_t) - B_i^1(Y_t,V_t',\nu_t), V_t - V_t' \rangle \d W^i_t, \quad m^1_0 =0.
			$$
			By Burkholder-Davis-Gundy inequality, Cauchy-Schwarz inequality applied several times, Hölder inequality, the fact that $4ab \leq a^2+4b^2$ and Lemma \ref{L estim lip}, there exists an increasing function $\kappa_p : \mathbb{R}_+ \rightarrow \mathbb{R}_+$ such that

			\begin{align*}
				\mathbb{E}\left[\sup_{t \in [0,T]}|m^1_t|^p\right] &\lesssim \mathbb{E} \left[ \left( \SN \int_0^T \langle B_i^1(X_t,V_t,\mu_t) - B_i^1(Y_t,V_t',\nu_t), V_t - V_t' \rangle^2 \d t \right)^{p/2} \right] \\
				& \leq \mathbb{E} \left[ \left( \sup_{t \in [0,T]} |V_t-V_t'|^{2} \right)^{p/2} \cdot \left( \SN \int_0^T  |B_i^1(X_t,V_t,\mu_t) - B_i^1(Y_t,V_t',\nu_t)|^{2} \d t \right)^{p/2} \right] \\
				& \leq \lambda_{p/2}(N) \SN \mathbb{E} \left[ \sup_{t \in [0,T]} |V_t -V_t'|^{p} \cdot \left( \int_0^T  |B_i^1(X_t,V_t,\mu_t) - B_i^1(Y_t,V_t',\nu_t)|^{2} \d t \right)^{p/2} \right] \\
				& \leq \lambda_{p/2}(N) \mathbb{E} \left[ \sup_{t \in [0,T]} |V_t -V_t'|^{2p} \right]^{1/2}  \SN \mathbb{E} \left[ \left( \int_0^T  |B_i^1(X_t,V_t,\mu_t) - B_i^1(Y_t,V_t',\nu_t)|^{2} \d t \right)^{p} \right]^{1/2} \\
				& \leq  \mathbb{E} \left[ \sup_{t \in [0,T]} |V_t-V_t'|^{2p} \right]^{1/2} \SN \mathbb{E} \left[ \lambda_{p/2}(N)^2  \kappa_p(T)  \left(  \int_0^T  |B_i^1(X_t,V_t,\mu_t) - B_i^1(Y_t,V_t',\nu_t)|^{2p} \d t \right) \right]^{1/2} \\
				& \leq \frac{1}{4} \mathbb{E} \left[ \sup_{t \in [0,T]} |V_t-V_t'|^{2p} \right] + \lambda_{p/2}(N)^2 \kappa_p(T) \left( \SN \mathbb{E} \left[ \int_0^T |B_i^1(X_t,V_t,\mu_t) - B_i^1(Y_t,V_t',\nu_t)|^{2p} \d t \right]^{1/2} \right)^2 \\
				& \leq \frac{1}{4} \mathbb{E} \left[ \sup_{t \in [0,T]} |V_t-V_t'|^{2p} \right] + N\lambda_{p/2}(N)^2 \kappa_p(T)  \SN \mathbb{E} \left[ \int_0^T |B_i^1(X_t,V_t,\mu_t) - B_i^1(Y_t,V_t',\nu_t)|^{2p} \d t \right] \\
				& \leq \frac{1}{4} \mathbb{E} \left[ \sup_{t \in [0,T]} |V_t-V_t'|^{2p} \right] \\
				& + KN  \lambda_{p/2}(N)^2 \kappa_p(T)  \int_0^T \mathbb{E} \left[ |V_t - V_t'|^{2p} \right] + \mathbb{E}\left[  \left(|X_t-Y_t|^{2p} + W_2^{2p}(\mu_t,\nu_t) \right) \cdot k_t^p \right] \d t \\
				& \leq \frac{1}{4} \mathbb{E} \left[ \sup_{t \in [0,T]} |V_t-V_t'|^{2p} \right] +KN  \lambda_{p/2}(N)^2 \kappa_p(T) \int_0^T \mathbb{E} \left[ |V_t - V_t'|^{2p} \right] \d t \\
				&+ KN  \lambda_{p/2}(N)^2 \kappa_p(T)\int_0^T \mathbb{E}\left[ \left( |X_t -Y_t|^{2p} + W_2^{2p}(\mu_t,\nu_t) \right)^{2} \right]^{1/2}  \mathbb{E}\left[ |k_t|^{2p} \right]^{1/2} \d t \\
				& \leq \frac{1}{4} \mathbb{E} \left[ \sup_{t \in [0,T]} |V_t-V_t'|^{2p} \right] +KN  \lambda_{p/2}(N)^2 \kappa_p(T) \int_0^T \mathbb{E} \left[ |V_t - V_t'|^{2p} \right] \d t \\
				&+ 2KN  \lambda_{p/2}(N)^2 \kappa_p(T)\int_0^T \mathbb{E}\left[ |X_t -Y_t|^{4p} + W_2^{4p}(\mu_t,\nu_t)  \right]^{1/2}  \mathbb{E}\left[ |k_t|^{2p} \right]^{1/2} \d t,
				\end{align*}
				where, for $a, b \in \mathbb{R}_+$,

				$$
				\lambda_a(b) = \begin{cases}
					b^{a-1}, \quad \text{if } a \geq 1, \\
					1, \quad \text{else}.
				\end{cases}
				$$
			By Proposition \ref{L inegalite croissance} there exists an increasing function $\tilde{C}^3_p: \mathbb{R}_+\rightarrow \mathbb{R}_+$ depending on $k_0$ such that $\mathbb{E}[|k_t|^{2p}] \leq \tilde{C}^3_p(T)$ for all $0 \leq t \leq T$. Moreover, by \eqref{Estimee lip sup flot}, 

			$$
			\int_0^T \mathbb{E}\left[ |X_t -Y_t|^{4p} + W_2^{4p}(\mu_t,\nu_t)  \right]^{1/2} \d t \leq T \sqrt{C_{2p}^1(T)} \left( |x-y|^{2p} + W_2^2(\mu,\nu) \right).
			$$
			Hence, there exists an increasing function $\tilde{K}^1_p(T) : \mathbb{R}_+ \rightarrow \mathbb{R}_+$ such that:

			\begin{align} \label{E 157}
				\mathbb{E}\left[\sup_{t \in [0,T]}|m^1_t|^p\right] &\leq \frac{1}{4} \mathbb{E} \left[ \sup_{t \in [0,T]} |V_t-V_t'|^{2p} \right] \\
			& + \tilde{K}_1(T) \left( \int_0^T \mathbb{E} \left[ |V_t - V_t'|^{2p} \right] \d t + \mathbb{E} \left[ |x-y|^{2p} + W_2^{2p}(\mu,\nu)  \right] \right) \nonumber.
			\end{align}
			On the other hand, by Hölder inequality and Cauchy-Schwarz inequality there exists an increasing function $\kappa'_p : \mathbb{R}_+ \rightarrow \mathbb{R}_+$ such that:

			\begin{align*}
				& \quad \mathbb{E}\left[ \sup_{t \in [0,T]}|l_t|^{p} \right] \\
				& \leq 2^{p-1} \mathbb{E} \left[ \left( \int_0^T \left[ |X_t-Y_t|^2 +W_2^2(\mu_t^1,\nu_t^1)\right]\cdot k_t \d t \right)^p + \left(\int_0^T |V_t-V_t'|^2 \d t \right)^p \right] \\ 
				& \leq 2^{p-1}\kappa'_p(T) \left( \int_0^T \mathbb{E}\left[ \left( |X_t -Y_t|^{2} + W_2^{2}(\mu_t^1,\nu_t^1) \right)^p \cdot k_t^p \right] \d t  + \mathbb{E} \left[ \int_0^T |V_t -V_t'|^{2p} \d t \right] \right) \\
				& \leq 2^{p-1} \kappa'_p(T) \left( \int_0^T \mathbb{E} \left[ \left(|X_t-Y_t|^{2}+W_2^2(\mu_t^1,\nu_t^1)\right)^{2p} \right]^{1/2} \cdot \mathbb{E} \left[ |k_t|^{2p} \right]^{1/2} \d t + \mathbb{E} \left[ \int_0^T |V_t -V_t'|^{2p} \d t \right] \right).
			\end{align*}
			Since

			$$
			 \mathbb{E} \left[ \left(|X_t-Y_t|^{2}+W_2^2(\mu_t^1,\nu_t^1)\right)^{2p} \right]^{1/2} \leq \mathbb{E} \left[ 2^{2p-1} \left(|X_t-Y_t|^{4p}+W_2^{4p}(\mu_t^1,\nu_t^1) \right) \right]^{1/2} = 2^{p-1/2}\mathbb{E} \left[ |X_t-Y_t|^{4p}+W_2^{4p}(\mu_t^1,\nu_t^1) \right]^{1/2}.
			$$
			and by using the same arguments as above, namely \eqref{Estimee lip sup flot} and the fact that there exists a constant $K(T)$ depending on $k_0$, $p$ and $T$ such that $\mathbb{E}[|k_t|^{2p}] \leq K(T)$, we obtain the existence of an increasing function $\tilde{K}_p^2 : \mathbb{R}_+ \rightarrow \mathbb{R}_+$ such that

			\begin{equation} \label{E 158}
				\mathbb{E}\left[ \sup_{t \in [0,T]}|l_t|^{p} \right] \leq \tilde{K}_p^2(T)\left( \mathbb{E} \left[ |x-y|^{2p} + W_2^{2p}(\mu,\nu)  \right] + \int_0^T \mathbb{E} \left[ |V_t - V_t'|^{2p} \right] \d t \right).
			\end{equation}
			
			Let 

			$$
			V_t^2 = \SN \int_0^t B_i^2(X_s, \zeta_s) \d W^i_s + \int_0^t B_0^2(X_s, \zeta_s) \d s, \quad \zeta_s=  (\varphi_s(\cdot,\mu), w_s(\cdot,\mu)(V))_*\mu
			$$
			and

			$$
				V_t^{'2} = \SN \int_0^t B_i^2(Y_s, \eta_s) \d W^i_s + \int_0^t B_0^2(Y_s, \eta_s) \d s, \quad  \eta_s = (\varphi_s(\cdot,\nu), w_s(\cdot,\nu)(V'))_*\nu.
			$$
			By using the exact same arguments as in the first part of the proof in order to control the $\min(|\zeta_t|^2,|\eta_t|^2)$ term, we obtain that there exists an increasing function $\tilde{K}^3_p : \mathbb{R}_+ \rightarrow \mathbb{R}_+$ such that

			\begin{equation}\label{E 159}
				\begin{array}{l} 
				\quad \mathbb{E}\left[ \sup_{t \in [0,T]}|V_t^2 -V_t^{'2}|^{2p} \right] \\
				 \leq \frac{1}{4} \mathbb{E} \left[ \sup_{t \in [0,T]} |V_t-V_t'|^{2p} \right] + \sup_{t \in [0,T]} \tilde{K}_3(T) \left( |x-y|^{2p} + W_2^{2p}(\zeta_0,\eta_0) + \int_0^T \mathbb{E}\left[ |V_t-V_t'|^{2p} \right]\d t \right).
			\end{array}
			\end{equation}
			
			Finally, by \eqref{E 157}, \eqref{E 158}, \eqref{E 159}, we obtain the existence of an increasing function $\tilde{K}_p : \mathbb{R}_+ \rightarrow \mathbb{R}_+$ depending on $p$ and $k_0$ such that:

			\begin{align*}
				\mathbb{E}\left[\sup_{t \in [0,T]}|V_t - V_t'|^{2p}\right] &\leq \tilde{K}_p(T) \left( |x-y|^{2p} + |V(x)-V'(y)|^{2p} + W_2^{2p}(\zeta_0,\eta_0) + \int_0^T \mathbb{E} [ |V_t - V_t'|^{2p}] \d t \right) \\
				& \leq \tilde{K}_p(T) \left( |x-y|^{2p} + |V(x)-V'(y)|^{2p} + W_2^{2p}(\zeta_0,\eta_0) \right)  + \left( \int_0^T \tilde{K}_p(T) \mathbb{E} \left[ \sup_{s \in [0,t]} |V_s - V_s'|^{2p} \right] \d t \right)
			\end{align*}
			Thus, since $T \mapsto \tilde{K}_p(T)$ is increasing, we can apply Gronwall's lemma to obtain the existence of  an increasing function $\tilde{K}'_p : \mathbb{R}_+ \rightarrow \mathbb{R}_+$ depending on $p$, $k_0$ such that:

			$$
			\mathbb{E}\left[\sup_{t \in [0,T]}|V_t - V_t'|^{2p}\right] \leq \tilde{K}'_p(T) \left( |x-y|^{2p} + |V(x)-V'(y)|^{2p} + W_2^{2p}(\zeta_0,\eta_0) \right).
			$$
		\end{proof}
	
\begin{prop}
		There exists a continuous version of
		
		\begin{align*}
			\mathbb{R}_+ \times M \times [0,1] &\longrightarrow M \times TM \\
			(t,x,a)   &\longmapsto (\varphi_t(U_a(x),\mu_a), w_t(U_a(x),\mu_a,V_a)).
		\end{align*}
	\end{prop}

	\begin{proof}
		Let $X^a_t(x) = \varphi_t(U_a(x),\mu_a)$, $V_t^a(x) = w_t(U_a(x),\mu_a,V_a)$. Since $|(V_a)_a|_\infty=\sup_{(x,a)\in M \times [0,1]}|V_a(x)|<\infty$, the function $C^4_p(T)$ in Proposition \ref{P Estim lip t cons} does not depend on $a,b \in [0,1]$. Thus,

		\begin{align*}
			& \quad \mathbb{E} \left[ |V_t^a(x) -V_s^b(y)|^{2p} \right] \\
			&\lesssim \mathbb{E}\left[|V_t^a(x)-V_t^b(y)|^{2p} \right] + \mathbb{E} \left[ |V_t^b(y) -V_s^b(y)|^{2p} \right] \\
			& \leq C^4_p(T) \left( |U_a(x)-U_b(y)|^{2p} + |V_a(x) - V_b(y)|^{2p} + W_2^{2p}((V_a)_*\mu,(V_b)_*\nu) \right) + \mathbb{E} \left[ |V_t^b(y) -V_s^b(y)|^{2p} \right] \\
			& \lesssim C^4_p(T) \left( |x-y|^{2p} + |a-b|^{2p} + W_2^{2p}(\mu,\nu) \right) + \mathbb{E} \left[ |V_t^b(y) -V_s^b(y)|^{2p} \right].
		\end{align*}
		Moreover, by Itô's isometry and Lemma \ref{L estim lip} it follows

		\begin{align*}
			\mathbb{E} \left[ |V_t^b(y) -V_s^b(y)|^{2p} \right] &= \mathbb{E} \left[ \left| \int_s^t \SN B_i(X^b_r(y), V_r^b(y), \zeta_r^b) \d W^i_r + \int_s^t B_0 (X^b_r(y), V_r^b(y), \zeta_r^b) \d r \right|^{2p} \right] \\
			& \lesssim \mathbb{E} \left[ \left( \int_s^t \SN \left| B_i(X^b_r(y), V_r^b(y), \zeta_r^b) \right|^2 \d r \right)^p \right] + \mathbb{E} \left[ \left| \int_s^t B_0 (X^b_r(y), V_r^b(y), \zeta_r^b) \d r \right|^{2p} \right] \\
			& \lesssim \mathbb{E} \left[ \left( \int_s^t C(1 + |V_r^b(y)|^2 + |\zeta_r^b|^2) \d r\right)^p \right] + \mathbb{E} \left[ \left( \int_s^t \sqrt{C(1 + |V_r^b(y)|^2 + |\zeta_r^b|^2 )} \d r \right)^{2p} \right].
		\end{align*}
		By taking the supremum, and applying Proposition \ref{L inegalite croissance} we finally obtain the existence of an increasing function $C_p^5 : \mathbb{R}_+ \to \mathbb{R}_+$ such that:
		
		\begin{align*}
			& \mathbb{E} \left[ \left( \int_s^t C(1 + |V_r^b(y)|^2 + |\zeta_r^b|^2) \d r\right)^p \right] + \mathbb{E} \left[ \left( \int_s^t \sqrt{C(1 + |V_r^b(y)|^2 + |\zeta_r^b|^2 )} \d r \right)^{2p} \right] \\
			& \leq C_p^5(T) \left(1 + |V_b(y)|^{2p} + |(V_b)_*\mu^b|^{2p} \right) |t-s|^p \\
			& \leq C_p^5(T) \left( 1 + 2|(V_a)_a|_\infty^{2p} \right) |t-s|^p
		\end{align*}
		Thus, there exists an increasing function $C_p^6 : \mathbb{R}_+ \to \mathbb{R}_+$ depending only on $T$ and on $|(V_a)_a|_\infty$ such that 

		$$
		\mathbb{E} \left[ |V_t^a(x) -V_s^b(y)|^{2p} \right] \leq C_p^6(T) \left( |x-y|^{2p} + |a-b|^{2p} + W_2^{2p}(\mu,\nu) + |t-s|^p \right).
		$$
		The conclusion follows from Proposition \ref{Prop continuite proc} and the Kolmogorov-$\check{\text{C}}$encov Theorem.
	\end{proof}

	 \begin{thm} \label{T C1}
 Let $k \geq 0$ and assume that $A_0, A_1, \dots, A_N \in \C^{\infty,k+1}(M \times \P(M),TM)$. Let also $(\varphi_t)_{t \geq 0}$ be the unique ID solution flow associated to \eqref{EDS principale Ito} and $(x,a) \in M \times [0,1] \mapsto U_a(x)$ be a smooth path of diffeomorphisms, $\mu_a = (U_a)_*\mu$. Then, there exists a $\C^{0,k,\infty}$ modification of

 \begin{align*}
 	\Omega \times \mathbb{R}_+ \times [0,1] \times M  & \longrightarrow M \\
 	(\omega,t,a,x) &\longmapsto \varphi_t^a(x,\omega)= \varphi_t(U_a(x),\mu_a)(\omega).
 \end{align*}
 Moreover, the field $(\omega,t,a,x) \mapsto \varphi_t^a(x,\omega)$ is an $(\mathscr{F}_t)_{t \geq 0}$-adapted $\C^k$ stochastic field of diffeomorphisms on $M$.
 \end{thm}

 \begin{proof}
	The first part of the theorem follows by iterating the arguments developed above for the spatial derivatives of the field. The second part can be obtained by showing that the Jacobian matrix of $\varphi_t(x,a)$ is nonsingular and by applying the inverse function theorem, see \cite[Lemma 4.4.3, Theorem 4.6.5]{kunita1997stochastic}.
 \end{proof}

 We can now make rigourous the heuristic reasoning presented at the beginning of the subsection:

	\begin{coro} \label{Coro dérivée processus en a}
		The path $(\M_t(\mu_a))_{0 \leq a \leq 1}$ belong to $\mathscr{R}_s(\M_t(\mu^0),\M_t(\mu^1))$ and its derivative is given almost surely by:

		$$
		\partial_a \M_t(\mu_a) = \Pi_{\mu^a_t} \left[ (\partial_a \varphi_t^a) \circ (\varphi_t^a)^{-1} \right].
		$$
	\end{coro}
	\begin{proof}
		The first point is a straightforward consequence of the fact that $(\varphi_t^a)$ is an $(\mathscr{F}_t)_{t \geq 0}$-adapted $\C^1$ stochastic field of diffeomorphisms. Let $f \in \C^\infty(M)$, we have 

		\begin{align*}
			\partial_a F_f(\M_t(\mu_a)) &= \partial_a \int_M f(\varphi_t^a(x)) \d \mu \\
			& = \int_M \langle \nabla f(\varphi_t^a(x)), \partial_a \varphi_t^a(x) \rangle \d \mu \\
			& = \int_M \langle \nabla f(x), (\partial_a \varphi_t^a(x)) \circ (\varphi_t^a)^{-1} \rangle \d \M_t(\mu_a) =\int_M \langle \nabla f(x), \Pi_{\mu^a_t} \left[ (\partial_a \varphi_t^a) \circ (\varphi_t^a)^{-1} \right] \rangle \d \M_t(\mu_a).
				\end{align*}
				The conclusion follows.
	\end{proof}

Throughout this subsection, the path $(U_a)_{0 \leq a \leq 1}$ has been used simultaneously as the spatial path and, in a certain sense, as the measure path. While this setting simplifies the presentation of the results, the following more general statement can still be established following the exact same steps:

 \begin{prop} \label{P general flow lisse}
		Let $k \geq 0$ and assume that $A_0, A_1, \dots, A_N \in \C^{\infty,k+1}(M \times \P(M),TM)$. Let $(U_a)_{0 \leq a \leq 1}$ be as above, $(U'_a)_{0 \leq a \leq 1}$ be a smooth path of $\C^\infty$ maps and let $\mu_a =(U_a)_* \mu$. Then, there exists a $\C^{0,k,\infty}$ modification of

 \begin{align*}
 	\Omega \times \mathbb{R}_+ \times [0,1] \times M  & \longrightarrow M \\
 	(\omega,t,a,x) &\longmapsto \mathcal{X}_t^a(x,\omega):= \varphi_t(U'_a(x),\mu_a)(\omega).
 \end{align*}
 Moreover, the field $(t,a,x) \mapsto \mathcal{X}_t^a(x)$ is an $(\mathscr{F}_t)_{t \geq 0}$-adapted $\C^k$ stochastic field of diffeomorphisms on $M$. Furthermore, for all $x \in M$, $\mathcal{V}_t = \partial_a |_{a=0}\mathcal{X}_t^a(x)$ belongs to $TM$ and is solution to:
		
		\begin{align*}
			\begin{array}{ll}
			D_t^\nabla \mathcal{V}_t &= \SN (\nabla A_i(\cdot,\mu_t))(\mathcal{V}_t) + \mathcal{\bar{L}}_{w_t(\mathcal{X}_t^{-1},\mu,V)} A_i(\mathcal{X}_t(x),\cdot)(\mu_t) \d W^i_t \\
			& +(\nabla A_0(\cdot,\mu_t))(\mathcal{V}_t) + \mathcal{\bar{L}}_{w_t(\mathcal{X}_t^{-1},\mu,V)} A_0(\mathcal{X}_t(x),\cdot)(\mu_t) \d t \\&
			 - \SN \frac{1}{2} R(\mathcal{V}_t^a,A_i(\mathcal{X}_t(x),\mu_t))A_i(\mathcal{X}_t(x),\mu_t) \d t,
			\end{array}
		\end{align*}
		where $\mathcal{X}_t=\mathcal{X}_t^0$. 
	\end{prop}

\subsection{Regularity of the semigroup}

In the following we will always consider a smooth extension of $\varphi_t$ to $\mathbb{R}^d \times \P$ and smooth extensions of the $A_is$ to $\mathbb{R}^d \times \P \rightarrow \mathbb{R}^d$. Moreover, $\nabla$ will denote here the gradient in $\mathbb{R}^d$. Depending on the context, $|\cdot|$ will denote the euclidean norm on $\mathbb{R}^d$ or the associated canonical norm on $\mathbb{R}^d \otimes \mathbb{R}^d$. Also, note that for $F \in \C^1(\P)$, $\grad F \in \C^0(M \times \P, TM)$. In the remainder of the subsection, we will always consider an associated compactly supported smooth extension $\grad F \in \C^0_c(\mathbb{R}^d \times \P, \mathbb{R}^d)$. The following two propositions can be proved by adapting the arguments in \cite[Subsections 3.1 and 3.2]{wang2020}. The only differences concern the Lipschitz estimates, which follow from the Lipschitz continuity of the driving vector fields, together with standard arguments based on Itô's formula, as carried out above.

	\begin{prop} \label{P gradient spat flow}
		Let $\mu \in \P$, $t \in \mathbb{R}_+$. The map $x \in \mathbb{R}^d \mapsto \varphi_t(x,\mu) \in L^2(\Omega, \mathbb{R}^d;\mathbb{P})$ is differentiable and its differential $\nabla \varphi_t(x,\mu) \in \mathbb{R}^d \otimes \mathbb{R}^d$ is solution to the following:

		$$
		\begin{cases}
			\d \nabla \varphi_t(x,\mu) = \SN \nabla A_i(\varphi_t(x,\mu),\M_t(\mu)) \cdot \nabla \varphi_t(x,\mu) \d W^i_t + \nabla A_0(\varphi_t(x,\mu),\M_t(\mu)) \cdot \nabla \varphi_t(x,\mu) \d t \\
			\nabla \varphi_0(x,\mu) = \mathrm{Id}, \, \M_t(\mu) = (\varphi_t(\cdot,\mu))_* \mu.
		\end{cases}
		$$
		Moreover, let $T>0$, $p \geq 1$. Then, for any $x,y \in \mathbb{R}^d$, $\mu,\nu \in \P$, and $s,t \in [0,T]$, there exists $C_p(T) >0$ such that

		$$
		\mathbb{E} \left[ |\nabla \varphi_t(x,\mu)- \nabla \varphi_s(y,\nu) |^{2p} \right] \leq C_p(T) \left( |t-s|^p + |x-y|^{2p} + W_2^{2p}(\mu,\nu) \right).
		$$
	\end{prop}

	\begin{prop}
		Let $x \in \mathbb{R}^d$, $t \in \mathbb{R}_+$. The map $\mu \in \P \mapsto \varphi_t(x,\mu) \in L^2(\Omega,\mathbb{R}^d;\mathbb{P})$ is intrinsically differentiable and for all $y \in \mathbb{R}^d$, its differential $\grad \varphi_t(x,\mu) \in L^2(\mathbb{R}^d,\mathbb{R}^d \otimes \mathbb{R}^d; \mu)$ satisfies the following:

		\begin{align*}
			& \quad \d \grad \varphi_t(x,\mu)(y) \\
			&= \SN \left( \grad \varphi_t(x,\mu)(y) \right)^\top \cdot \nabla A_i(\varphi_t(x,\mu), \M_t(\mu)) + \left( \nabla \varphi_t(x,\mu) \right)^\top \cdot \left[\grad A_i(\varphi_t(x,\mu), \M_t(\mu))\right](\varphi_t(y,\mu)) \d W^i_t \\
			& +\SN \int_{M} \left( \grad \varphi_t(z,\mu)(y) \right)^\top \cdot \left[ \grad A_i (\varphi_t(x,\mu),\M_t(\mu))\right](\varphi_t(z,\mu)) \d \mu(z) \d W^i_t \\
			& + \left( \grad \varphi_t(x,\mu)(y) \right)^\top \cdot \nabla A_0(\varphi_t(x,\mu), \M_t(\mu)) + \left( \nabla \varphi_t(x,\mu) \right)^\top \cdot \left[ \grad A_0(\varphi_t(x,\mu), \M_t(\mu)) \right](\varphi_t(y,\mu)) \d t \\
			& + \int_{M} \left( \grad \varphi_t(z,\mu)(y) \right)^\top \cdot\left[ \grad A_0 (\varphi_t(x,\mu),\M_t(\mu))\right](\varphi_t(z,\mu)) \d \mu(z) \d t, \\
			& \grad \varphi_0(x,\mu)(y) = 0, \, \M_t(\mu)=(\varphi_t(\cdot,\mu))_*\mu,
		\end{align*}
		Moreover, let $T>0$, $p \geq 1$. There exists $C_p(T)$ such that 

		$$
		\mathbb{E} \left[ |\T_\mu \varphi_t(x,\cdot) (V)|^{2p} \right] \leq C_p(T) |V|_\mu^{2p}, \quad \forall V \in L^2(\mathbb{R}^d,\mathbb{R}^d;\mu),
		$$
		consequently, $\T_\mu \varphi_t(x,\cdot) \in \mathcal{L}(L^2(\mathbb{R}^d,\mathbb{R}^d;\mu), L^2(\Omega,\mathbb{R}^d;\mathbb{P}))$.
		
		Furthermore, for any $x,y \in \mathbb{R}^d$, $\mu,\nu \in \P$, and $s,t \in [0,T]$, there exists $C_p'(T) > 0$ such that

		$$
		\mathbb{E} \left[ |\grad \varphi_t(x_1,\mu)(y_1)- \grad \varphi_s(x_2,\nu)(y_2) |^{2p} \right] \leq C_p'(T) \left( |t-s|^p + |x_1-x_2|^{2p} + |y_1-y_2|^{2p} + W_2^{2p}(\mu,\nu) \right).
		$$
		
	\end{prop}

	\begin{prop} \label{Stab semigroup}
		Let $t>0$, $F \in \C^k(\P)$, then, $\mu \in \P \mapsto \se_t F(\mu)$ belongs to $\C^{k}(\P)$, for $k=0,1,2$. Consequently, $\se_t F \in \mathrm{D}(\bar{\mathrm{L}})$. 
	\end{prop} 

	\begin{proof}
		First, it is clear that $(t,\mu) \in [0,T] \times \P \mapsto \se_t F(\mu)$ belongs to $\C^{0,0}([0,T]\times \P)$. By Corollary \ref{Lemme fondamental L2}, $\mu \mapsto \se_t F(\mu)$ is $L$-differentiable and we have for $x,y \in M$ and $\mu,\nu \in \P$:

		$$
		\grad \se_t F(\mu)(x) = \mathbb{E} \left[ \left( \nabla \varphi_t(x,\mu) \right)^\top \cdot \grad F (\M_t(\mu))(\varphi_t(x,\mu)) + \int_M \left(\grad \varphi_t(z,\mu)(y)\right)^\top \cdot \grad F(\M_t(\mu))(\varphi_t(z,\mu))  \d \mu(z)\right].
		$$
		Let us denote $X_t = \varphi_t(x,\mu)$, $Y_t=\varphi_t(y,\nu)$, $\mu_t = \M_t(\mu)$, $\nu_t = \M_t(\nu)$. Then,

		\begin{align*}
			& \quad \left| \mathbb{E} \left[ \left( \nabla X_t \right)^\top \cdot \grad F (\mu_t)(X_t) - \left( \nabla Y_t \right)^\top \cdot \grad F (\nu_t)(Y_t) \right] \right|^2 \\
			& \leq  \mathbb{E} \left[ \left| \left( \nabla X_t \right)^\top \cdot \grad F (\mu_t)(X_t) - \left( \nabla Y_t \right)^\top \cdot \grad F (\nu_t)(Y_t) \right| \right] ^2 \\
			& \leq \mathbb{E} \left[ |\left( \nabla X_t \right)^\top \cdot \left( \grad F (\mu_t)(X_t) - \grad F(\nu_t)(Y_t) \right) | + | \left( \nabla X_t - \nabla Y_t \right)^\top \cdot \grad F (\nu_t)(Y_t)  \right]^2 \\
			& \lesssim \mathbb{E} \left[ |\nabla X_t|^2 \right] \mathbb{E} \left[ |\grad F (\mu_t)(X_t) - \grad F(\nu_t)(Y_t)|^2 \right] + \mathbb{E} \left[ |\grad F (\nu_t)(Y_t)|^2 \right] \mathbb{E} \left[ |\nabla X_t - \nabla Y_t|^2 \right] \\
			& \leq C \left[ W_2^2(\mu,\nu) + |x-y|^2 \right],
		\end{align*}
		where the last equality follows from the fact that $M \times \P$ is compact, $\mathbb{E}[|\nabla \varphi_t|^2] \in \C^0_b(M \times \P)$, $F \in \C^2(\P)$, $\grad F \in \C^0(M \times \P, \mathbb{R}^d)$ and Proposition \ref{P gradient spat flow}. The second term can be treated in the same way.
		
		One can proceed as in \cite[Subsection 3.3]{wang2020} to show existence and continuity of the second order derivatives of the stochastic flow. Following the same principle as in \cite[Subsection 3.4]{wang2020}, one can show that $\se_{t} F \in \C^2(\P)$.
	\end{proof}

	\section{Itô-Wasserstein curvature condition and consequences} \label{S 4} 

	\subsection{Covariant SDEs} \label{S 3.3}

\noindent \textbf{Covariant SDEs on Hilbert vector bundles.} Let $\mathscr{M}$ be a Hilbert manifold, $\pi_1 : \mathscr{F}^1 \rightarrow \mathscr{M}$ and $\pi_2 : \mathscr{F}^2 \rightarrow \mathscr{M}$ be Hilbert vector bundles over $\mathscr{M}$. Let also $\mathscr{K}$ be a Hilbert space. 
    
    \begin{defi}
        Let $\Psi : \mathscr{M} \rightarrow \mathscr{M}$. A map $J : \mathscr{F}^1 \rightarrow \mathscr{F}^2$ is said to belong to $\mathrm{Hom}_\Psi(\mathscr{F}^1,\mathscr{F}^2)$ if the diagram
        \begin{center}
\begin{tikzpicture}[>=stealth]
  \node (E1) at (0,1) {$\mathscr{F}^1$};
  \node (E2) at (4,1) {$\mathscr{F}^2$};
  \node (E3) at (0,0)   {$\mathscr{M}$};
  \node (E4) at (4,0)   {$\mathscr{M}$};

  \draw[->, >=latex] (E1) -- (E2) node[midway,above]{$J$}; 
  \draw[->, >=latex] (E1) -- (E3) node[midway,left]{$\pi_1$};
  \draw[->, >=latex] (E3) -- (E4) node[midway,below]{$\Psi$}; 
  \draw[->, >=latex] (E2) -- (E4) node[midway,right]{$\pi_2$};
\end{tikzpicture}
\end{center}
is commutative and for all $x \in \mathscr{M}$, the induced map $J_x : \mathscr{F}^1_x \rightarrow \mathscr{F}^2_{\Psi(x)}$ is a continuous linear map. Let $k \in \mathbb{N}$, $J$ is said to belong to $\mathrm{Hom}_{\Psi}^k(\mathscr{F}^1,\mathscr{F}^2)$ if $\Psi \in \C^k(\mathscr{M},\mathscr{M})$ and for all $x \in \mathscr{M}$ there exists $U,U'$ two neighborhoods of $x,\Psi(x)$ and trivializing maps

$$
\theta : \pi_1^{-1}(U) \longrightarrow U \times E
$$
and

$$
\theta' : \pi_2^{-1}(U') \longrightarrow U' \times E'
$$
at $x$ and $\Psi(x)$ respectively, with $E, E'$ some Hilbert spaces, such that $\Psi(U) \subset U'$, and such that the map of $U$ into $\mathcal{L}_c(E,E')$ given by

$$
x \longmapsto \theta' \circ J_x \circ \theta^{-1}
$$
is a $\C^k$ map.
    \end{defi}
    
    Let $(X_t)_{t \geq 0}$ be solution to 
	
	\begin{equation*} 
	    \circ \d X_t = \SN A_i(X_t) \circ \d W^i_t + A_0(X_t) \d t, \quad X_0 = x \in \mathscr{M},
	\end{equation*}
	with $A_0, A_1, \dots,  A_N \in \Gamma^2(T\mathscr{M})$. Let $\pi : \mathscr{F} \rightarrow \mathscr{M}$ be a smooth vector bundle endowed with a symmetric bilinear form $\langle \bullet,\bullet \rangle : \mathscr{F} \times \mathscr{F} \rightarrow \mathbb{R}$ and a connection $ \nabla$ satisfying for $A,B \in \Gamma(\mathscr{F})$ and $V \in TM$:
	
	$$
	V \langle A,B \rangle = \langle \nabla_V A, B \rangle + \langle A, \nabla_{V} B \rangle. 
	$$
	Assume that there exists a unique stochastic parallel transport $(\tau_{0,t})_{t \geq 0}$ along $(X_t)_{t \geq 0}$, namely, a stochastic process $(\tau_{0,t})_{t}$ in $O^x(M)$ which is invertible for all $t \geq 0$ and such that for all $V \in T_x M$ and all $B \in \Gamma(\mathscr{F})$:

	$$
    \circ \d \langle \tau_{0,t}(V), B(X_t) \rangle_{X_t} = \SN \langle \tau_{0,t}(V), \nabla_{A_i(X_t)} B(X_t) \rangle_{X_t} \circ \d W^i_t + \langle \tau_{0,t}(V), \nabla_{A_0(X_t)}B(X_t) \rangle_{X_t} \d t,
    $$ 

	Let $J_0, J_1, \dots, J_N \in \mathrm{Hom}_{\mathrm{id}}^0(\mathscr{F},\mathscr{F})$ be continuous vector bundle morphisms. An adapted $\mathscr{F}$-valued stochastic process $(U_t)_{t \geq 0}$ along $(X_t)_{t \geq 0}$ is said to be solution to the following covariant SDE:
	
	\begin{equation} \label{EDS covariante générale}
	    \circ D_t U_t = \sum_{i=1}^{N} J_i(U_t) \circ \d W^i_t + J_0(U_t) \d t, \quad U_0 = U,
	\end{equation}
	if $(\tau_{0,t}^{-1} U_t)_{t \geq 0}$ is solution to the following SDE on $\mathscr{F}_\varphi$:
	
	$$
	\circ \d \tau_{0,t}^{-1} U_t = \sum_{i=1}^{N} \tau_{0,t}^{-1} J_i(U_t) \circ \d W^i_t + \tau_{0,t}^{-1} J_0 (U_t) \d t, \quad U_0 = U.
	$$
 	In the same way, we say that $(U_t)_{t \geq 0}$ is solution to

	\begin{equation*}
	    D_t^\nabla U_t = \sum_{i=1}^{N} J_i(U_t) \d W^i_t + J_0(U_t) \d t, \quad U_0 = U,
	\end{equation*}
	if $(\tau_{0,t}^{-1}(U_t))_{t \geq 0}$ is solution to

	$$
	\d \tau_{0,t}^{-1} U_t = \sum_{i=1}^{N} \tau_{0,t}^{-1} J_i(U_t) \d W^i_t + \tau_{0,t}^{-1} J_0 (U_t) \d t, \quad U_0 = U.
	$$

\noindent \textbf{Covariant SDEs on $T\mathscr{D}$.} Let $s > m/2$, $j(s) \geq s$. One can easily check that the injection $\iota : \mathscr{D}^{j(s)} \rightarrow \mathscr{D}^s$ is a smooth map. We denote by $\mathscr{F}^s = \iota^\ast T\mathscr{D}^s = \sqcup_{\varphi \in \mathscr{D}^{j(s)}}T_\varphi \mathscr{D}^s$ the pullback tangent bundle over $\mathscr{D}^{j(s)}$. We denote by $\bm \nabla$ the pullback connection $\iota^*(\bm \nabla)$ on $\mathscr{F}^s$ and for any continuous semimartingale $(\varphi_t)_{t \geq 0}$ on $\mathscr{D}^{j(s)}$, $\tau_{0,t}$ is the unique stochastic parallel transport on $\mathscr{F}^s$. Furthermore, for $\bm A \in \Gamma^k(T\mathscr{D}^{j(s)})$, we will also denote by $\bm A$ the $\C^k$ section $\iota_*(\bm A) \in \Gamma^k(\mathscr{F}^s)$. Let
    	
    	\begin{equation*}
    	    \mathscr{J} := \left\{ j : \mathbb{N}^* \rightarrow \mathbb{N}^* \text{ increasing function} \, : \, j(s) \geq s, \, \forall s >m/2 \right\}.
    	\end{equation*}
    	Let $j \in \mathscr{J}$, then $T\mathscr{D}$ is endowed with the ILH manifold structure given by:
    	
    	$$
    	T\mathscr{D} = \bigcap_{s>m/2} \mathscr{F}^s.
    	$$
    	Note that for any $j \in \mathscr{J}$, and any $s \geq s' >m/2$, the injection $\mathscr{F}^s \rightarrow \mathscr{F}^{s'}$ is a smooth map. Indeed, in a neighborhood of $\varphi \in \mathscr{D}^s$, this injection is just the natural injection:
    	
    	$$
    	T_\varphi \mathscr{D}^{j(s)} \times T_\varphi \mathscr{D}^s \longrightarrow T_\varphi \mathscr{D}^{j(s')} \times T_\varphi \mathscr{D}^{s'}.
    	$$
    	In what follows, even when $j(s) \in \mathscr{J}$ is not mentioned explicitly, $\mathscr{F}^s$ is understood to be a vector bundle over $\mathscr{D}^{j(s)}$ with its natural projection.
	
	\begin{defi} \label{Defi vector bundle morphism}
	    Let $\Psi \in \C^k_S(\mathscr{D},\mathscr{D})$. A map $\bm J : T\mathscr{D} \rightarrow T\mathscr{D}$ is said to belong to $\mathrm{Hom}^k_{\Psi}(T\mathscr{D}, T\mathscr{D})$ if there exists $s_0 \in \mathbb{N}$ and $j\in \mathscr{J}$ such that for all $s>s_0$ there exists a $\C^k$ extension $\bm J^s \in \mathrm{Hom}_{\Psi}^k(\mathscr{F}^s, \mathscr{F}^s)$. 
	\end{defi}

	Let $(\varphi_t)_{t \geq 0}$ be the unique solution to the following SDE:
	
	\begin{equation*} 
	    \circ \dd \varphi_t = \SN \bm A_i(\varphi_t) \circ \d W^i_t + \bm A_0(\varphi_t) \d t, \quad \varphi_0 = \varphi,
	\end{equation*}
	with $\bm A_0, \bm A_1, \dots,  \bm A_N \in \Gamma^2_S(T\mathscr{D})$. Let $\bm J_0, \bm J_1, \dots, \bm J_N \in \mathrm{Hom}_{\Psi}^0(T\mathscr{D},T\mathscr{D})$. A $T\mathscr{D}$-valued stochastic process $(\bm A_t)_{t \geq 0}$ along $(\varphi_t)_{t \geq 0}$ is said to be solution to the following covariant SDE:
	
	\begin{equation} \label{EDS covariante générale ILH}
	    \circ \bm D_t \bm A_t = \SN \bm J_i(\bm A_t) \circ \d W^i_t + \bm J_0(\bm A_t) \d t, \quad \bm A_0 = A \in T_\varphi \mathscr{D},
	\end{equation}
	if there exists $s_0 \in \mathbb{N}$ such that for all $s > s_0$, $(\bm A_t)_{t \geq 0}$ is solution to
	
	\begin{equation*}
	    \circ \bm D_t \bm A_t = \SN \bm J^s_i(\bm A_t) \circ \d W^i_t + \bm J^s_0(\bm A_t) \d t, \quad \bm A_0 = A \in T_\varphi \mathscr{D},
	\end{equation*}
	with $\bm J_i^s \in \mathrm{Hom}_{\mathrm{id}}^0(\mathscr{F}^s, \mathscr{F}^s)$ a continuous extension of $\bm J_i$ as in Definition \ref{Defi vector bundle morphism} for all $i=0,1\dots,N$. Analogously, $(\bm A_t)_{t \geq 0}$ is said to be solution to

	\begin{equation*} 
	    \bm D_t^{\bm \nabla} \bm A_t = \SN \bm J_i(\bm A_t) \d W^i_t + \bm J_0(\bm A_t) \d t, \quad \bm A_0 = A \in T_\varphi \mathscr{D},
	\end{equation*}
	if there exists $s_0 \in \mathbb{N}$ such that for all $s>s_0$, $(\bm A_t)_{t \geq 0}$ is solution to

	\begin{equation*}
	    \bm D_t^{\bm \nabla} \bm A_t = \SN \bm J^s_i(\bm A_t) \d W^i_t + \bm J^s_0(\bm A_t) \d t, \quad \bm A_0 = A \in T_\varphi \mathscr{D},
	\end{equation*}
	with $\bm J_i^s \in \mathrm{Hom}_{\mathrm{id}}^0(\mathscr{F}^s, \mathscr{F}^s)$ a continuous extension of $\bm J_i$ as in Definition \ref{Defi vector bundle morphism} for all $i=0,1\dots,N$. 
	
	One can consider $\mathrm{P} : \mathrm{GL}^{\mathrm{id}}(TM) \to M$ to be the vector bundle with fiber $\mathrm{GL}^{\mathrm{id}}_x(TM) := \left\{ u \in \mathrm{GL}(T_xM,T_xM)\right\}$ and consider the Hilbert Lie group given by $\mathscr{G} = H^s(M,\mathrm{GL}^{\mathrm{id}}(TM))$. Then, one can consider $\mathscr{G}$ as a structure group acting on $T\mathscr{D}^s$ and consider the associated principal bundle $\mathcal{P} : \mathscr{G}(T\mathscr{D}^s) \to \mathscr{D}^s$ where the fibers $\mathscr{G}_\varphi(T\mathscr{D}^s)$ are given by the set of linear isomorphisms $u : T_{\mathrm{id}}\mathscr{D}^s \rightarrow T_\varphi \mathscr{D}^s$, see \cite[Appendix A, (G)]{Elworthy1982} for more details. Note that this vector bundle is thus a Hilbert manifold. Then, one can show that there exists $\C^2$ horizontal vector fields $\tilde{\mathcal{A}}_0, \tilde{\mathcal{A}}_1,\dots,\tilde{\mathcal{A}}_N$ on $\mathscr{G}(T\mathscr{D}^s)$ such that 

    $$
    \circ \d \tau_{0,t} = \SN \tilde{\mathcal{A}}_i(\tau_{0,t}) \d W^i_t + \tilde{\mathcal{A}}_0(\tau_{0,t}) \d t, \quad \tau_{0,0}= \mathrm{Id}_{T_\varphi \mathscr{D}^s},
    $$
    with $\bm T  \mathcal{P}(\mathcal{A}_i(u))=A_i( \mathcal{P}(u))$ and $\mathcal{A}_i(R'_g(u)) = \bm T R_g' \mathcal{A}_i(u)$ for all $u \in \mathscr{G}(T\mathscr{D}^s)$ and all $g \in \mathscr{G}$, see \cite[Chapter VII, \textsection 13, (C)]{Elworthy1982}. In particular we can apply all the results in \cite[Chapter VII, \textsection 13]{Elworthy1982}. 
	
	\begin{prop} \label{P Existence unicité ILH covariante ILH}
	    Assume that $\bm J_i\in \mathrm{Hom}^1_{\mathrm{id}}(T\mathscr{D},T\mathscr{D})$ for all $i=0,1,\dots,N$. Then, there exists a unique continuous solution to \eqref{EDS covariante générale ILH} defined for all $t \in \mathbb{R}_+$.
	\end{prop}

	\begin{prop} \label{Ito Strato covariant}
		Let $\bm J_0, \bm J_1, \dots, \bm J_N \in \mathrm{Hom}^1_{\mathrm{id}}(T\mathscr{D},T\mathscr{D})$. Then, the Stratonovich SDE

		\begin{equation*}
	    \circ \bm D_t \bm A_t = \SN \bm J_i(\bm A_t) \circ \d W^i_t + \bm J_0(\bm A_t) \d t, \quad \bm A_0 = A \in T_\varphi \mathscr{D},
	\end{equation*}
	is equivalent to the following Itô's one:

	\begin{equation*} 
	    \bm D_t^{\bm \nabla} \bm A_t = \SN \bm J_i(\bm A_t) \d W^i_t + \bm J_0'(\bm A_t) \d t, \quad \bm A_0 = A \in T_\varphi \mathscr{D},
	\end{equation*}
	with

	$$
\bm J'_0(\cdot) = \bm J_0(\cdot) + \frac{1}{2} \SN \left[ \bm \nabla \bm J_i (\bm A_i)(\cdot) +\bm J_i (\bm J_i(\cdot))\right].
	$$

	\end{prop}
	
	\begin{prop} \label{Vector bundle chain rule ILH}
	    Let $(\bm A_t)_{t \geq 0}$ be solution to \eqref{EDS covariante générale ILH}. Then, for all $\C^2$ linear form $\Phi : T\mathscr{D} \rightarrow \mathscr{K}$ with $\mathscr{K}$ a Hilbert space, the following holds:
	
	\begin{equation*}
	    \circ \dd \Phi(\bm A_t) = (\bm \nabla \Phi(\circ \dd \varphi_t))(\bm A_t) + \Phi(\circ \bm D_t \bm A_t).
	\end{equation*}
	As a consequence, for $\bm L\in \mathrm{Hom}^2_{\Psi}(T\mathscr{D},T\mathscr{D})$ with $\Psi \in \C^2_S(\mathscr{D},\mathscr{D})$, $(\bm L_{\varphi_t}(\bm A_t))_{t \geq 0}$ is a semimartingale solution to:
	    
	    \begin{equation*}
	        \circ \bm D_t \bm L(\bm A_t) = (\bm \nabla\bm L_{\varphi_t}(\circ \dd \varphi_t))(\bm A_t) + \bm L_{\varphi_t}(\circ \bm D_t \bm A_t).
	    \end{equation*}
	\end{prop} 
	
	\vspace{10 pt}

	\noindent \textbf{Covariant SDEs on $T\P_\infty$.} Let $\mu \in \P_\infty$ and $(\mu_t)_{t \geq 0}$ be solution to
    
	\begin{equation} \label{E 517}
	    \circ \d\mu_t = \SN \bar{Z}_i(\mu_t) \circ \d W^i_t + \bar{Z}_0(\mu_t) \d t, \quad \mu_0 =\mu,
	\end{equation}   
	where for all $i=0,1, \dots, N$, $\bar{Z}_i \in \Gamma^2(T\P)$. Let $(\Phi_t)_{t \geq 0}$ be the horizontal lift of $(\mu_t)_{t \geq 0}$, namely the solution to
	
	\begin{equation} \label{EDS du lift horizontal}
	\circ \dd\Phi_t(\varphi) = \SN \bm Z_i(\Phi_t(\varphi)) \circ \d W^i_t + \bm Z_0(\Phi_t(\varphi)) \d t, \quad \Phi_0(\varphi) = \varphi \in \mathscr{D}_\mu,
	\end{equation}
	where $\bm Z_i \in \Gamma^3_S(T\mathscr{D})$ is the horizontal lift of $\bar{Z}_i$ for all $i=0,1, \dots, N$.
	
	\begin{defi}
	    A map $\bar{J} : T \P_\infty \rightarrow T \P_\infty$ such that for all $\mu \in \P_\infty$, $\bar{J}_\mu \in \mathcal{L}_c(T_\mu \P_\infty, T_\mu \P_\infty)$ is said to belong to $\mathrm{Hom}_{\mathrm{id}}^k(T\P_\infty, T \P_\infty)$ if there exists $\bm J \in \mathrm{Hom}_{\mathrm{id}}^k(T\mathscr{D},T\mathscr{D})$ such that the diagram:
	    
	    \begin{center}
\begin{tikzpicture}[>=stealth]
  \node (E1) at (0,1) {$T\mathscr{D}$};
  \node (E2) at (4,1) {$T\mathscr{D}$};
  \node (E3) at (0,0)   {$T\P_\infty$};
  \node (E4) at (4,0)   {$T\P_\infty$};

  \draw[->, >=latex] (E1) -- (E2) node[midway,above]{$\bm J$}; 
  \draw[->, >=latex] (E1) -- (E3) node[midway,left]{$Tp$};
  \draw[->, >=latex] (E3) -- (E4) node[midway,below]{$\bar{J}$}; 
  \draw[->, >=latex] (E2) -- (E4) node[midway,right]{$Tp$};
\end{tikzpicture}
\end{center}
commutes.
	\end{defi}
	We will forget about the subscript $\mu$ when there is no ambiguity about the point considered.
	
	Let $\bar{J}_0, \bar{J}_1, \dots, \bar{J}_N \in \mathrm{Hom}_{\mathrm{id}}^0(T\P_\infty,T\P_\infty)$. A $T\P_\infty$-valued stochastic process $(\bar{A}_t)_{t \geq 0}$ along $(\mu_t)_{t \geq 0}$ is said to be solution to 
	
	$$
	\circ \bar{D}_t \bar{A}_t = \SN \bar{J}_i(\bar{A}_t) \circ \d W^i_t + \bar{J}_0(\bar{A}_t) \d t, \quad \bar{A}_0 = A \in \T_\mu \P_\infty,
	$$
	if for all $\varphi \in \mathscr{D}_\mu$, $(\bm A_t)_{t \geq 0} = (\mathfrak{h}_{\Phi_t(\varphi)}(\bar{A}_t))_{t \geq 0}$ is solution to
	
	\begin{equation*}
	    \circ \bm D_t \bm A_t = \SN \bm J_i(\bm A_t) \circ \d W^i_t + \bm J_0(\bm A_t) \d t, \quad \bm A_0 = \mathfrak{h}_{\varphi}(A) \in T_\varphi \mathscr{D},
	\end{equation*}
	with $\bm J_i = \tilde{\bm J}_i + \mathcal{O}_{\bm Z_i}$, where
	
	$$
	\tilde{\bm J}_{i,\varphi} = \mathfrak{h}_{\varphi} \circ \bar{J}_i \circ \bm T_\varphi p, \quad \forall i=1, \dots,N
	$$
	and with $\bm J_0$ such that
	
	$$
	\bm T p \circ \bm J_0 = \bar{J}_0 \circ \bm T p.
	$$
	In the same way, $(\bar{A}_t)_{t \geq 0}$ is said to be solution to 
	
	\begin{equation} \label{EDS covariante TP}
		\bar{D}^{\grad}_t \bar{A}_t = \SN \bar{J}_i(\bar{A}_t) \d W^i_t + \bar{J}_0(\bar{A}_t) \d t, \quad \bar{A}_0 = A \in \T_\mu \P_\infty,
	\end{equation}
	if for all $\varphi \in \mathscr{D}_\mu$, $(\bm A_t)_{t \geq 0} = (\mathfrak{h}_{\Phi_t(\varphi)}(\bar{A}_t))_{t \geq 0}$ is solution to
	
	\begin{equation*}
	    \bm D^{\grad}_t \bm A_t = \SN \bm J_i(\bm A_t) \d W^i_t + \bm J_0(\bm A_t) - \frac{1}{2} \SN \mathcal{O}^*_{\bm Z_i} \left( P_{\V} (\bm J_i(\bm A_t)) \right) \d t, \quad \bm A_0 = \mathfrak{h}_{\varphi}(A) \in T_\varphi \mathscr{D},
	\end{equation*}
	with $\bm J_i$ defined as above for $i=0,1,\dots,N$.

	\begin{rem}
		On a finite-dimensional Riemannian manifold $M$, a stochastic process $(A_t)_t$ on $TM$ is said to be a solution to 

		$$
		D_t^\nabla A_t = G(A_t) \d N_t,
		$$
		where $G : TM \rightarrow TM \otimes \mathbb{R}^k$ and $(N_t)_t$ is a $ \mathbb{R}^k$ valued semimartingale if the following holds:
		
		$$
		\d \langle A_t, V_t \rangle = \langle D^\nabla_t A_t, V_t \rangle + \langle A_t, D_t^\nabla V_t \rangle + \langle D_t^\nabla A_t, D_t^\nabla V_t \rangle,
		$$
		for any semimartingale $(V_t)_t$ along $(\pi(A_t))_t$. One can actually check that our definition of solutions to covariant SDEs on $T\P_\infty$ is consistent with the finite-dimensional one. In particular, for a solution to \eqref{EDS covariante TP}, it holds

		$$
		\d |\bar{A}_t|^2 = 2 \langle \bar{D}_t^{\grad} A_t, A_t \rangle + |\bar{D}_t^{\grad} A_t|^2.
		$$
	\end{rem}

	\subsection{Tangent process}

	From now on and until the end of the article, we will consider solutions to \eqref{EDS principale Ito} (or equivalenty to \eqref{EDS principale Strato}) with initial conditions in $\P_\infty$. The two following results can be interpreted as stochastic infinite-dimensional counterparts of the Schwarz's theorem.

	Let $(\M_t)$ be the unique solution flow to \eqref{EDS principale Ito} and $(\Phi_t)_{t \geq 0}$ be the unique solution flow to the horizontal lift equation, namely,

	$$
		\dd^{\bm \nabla} \Phi_t(\varphi) = \SN \bm Z_i(\Phi_t(\varphi)) \d W^i_t + \bm Z_0(\Phi_t(\varphi)) \d t, \quad \Phi(\varphi) = \varphi,
	$$
	with $\bm Z_i = \mathfrak{h}(\bar{Z}_i)$.

	\begin{prop} \label{Strato schwarz lemma D}
		There exists a unique solution to the following covariant SDE:

		\begin{equation} \label{E 321}
			\circ \bm D_t \bm V_t = \SN \bm \nabla_{\bm V_t} \bm Z_i(\Phi_t(\varphi)) \circ \d W^i_t + \bm \nabla_{\bm V_t}\bm Z_0(\Phi_t(\varphi)) - \frac{1}{2} \SN \bm \nabla_{\bm V_t} \bm \nabla_{\bm Z_i} \bm Z_i(\Phi_t(\varphi)) \d t, \quad \bm V_0 = A \circ \varphi \in T_\varphi \di,
		\end{equation}
		which can be equivalently rewritten in the concise form

		\begin{equation}
			\circ \bm D_t \bm V_t = \bm \nabla_{\bm V_t} \circ \dd \Phi_t(\varphi).
		\end{equation}
	\end{prop}

	\begin{proof}
		The existence and uniqueness of the solution to this covariant SDE is guaranteed by Proposition \ref{P Existence unicité ILH covariante ILH}, take $A_t = \bm V_t$ and $\bm J_i = \bm \nabla \bm Z_i$.
	\end{proof}

	\begin{rem}
		Let $(\varphi_a)_a$ be a smooth path of diffeomorphisms with $\partial_a|_{a=0} \varphi_a = V \circ \varphi$. The process $(\bm V_t)_{t \geq 0}$ has to be understood as $\bm V_t =\bm T_{\varphi} \Phi_t(V \circ \varphi) = \partial_a|_{a=0} \Phi_t(\varphi^a)$, but at this point we don't know if $a \mapsto \Phi_t(\varphi^a)$ is differentiable.
	\end{rem}
	An application of Proposition \ref{Ito Strato covariant} yields:

	\begin{prop} \label{P Equation Ito D}
		The process $(\bm V_t)_{t \geq 0}$ is the unique solution to the following Itô covariant SDE:

		\begin{equation} \label{Proc Tang}
			 \bm D_t^{\bm \nabla} \bm V_t = \SN \bm \nabla_{\bm V_t} \bm Z_i(\Phi_t(\varphi)) \d W^i_t + \bm \nabla_{\bm V_t}\bm Z_0(\Phi_t(\varphi)) - \frac{1}{2} \SN \bm R (\bm V_t, \bm Z_i ) \bm Z_i  \d t, \quad \bm V_0 = A \circ \varphi \in T_\varphi \di.
		\end{equation}
		As a consequence, for any $x_0 \in M$, $\bm T \bm E_{x_0}(\bm V_t) \sim w_t(\varphi(x_0),p(\varphi), A)$, where $\bm E_{x_0} : \varphi \in \mathscr{D} \mapsto \varphi(x_0)\in M$ is the evaluation map which is smooth (see \cite[Appendix A.3]{Martin2025}), $w_t$ is the process in \eqref{MVC TM} and $\sim$ means indstinguishable.

	\end{prop}

	\begin{proof}
		By Proposition \ref{Ito Strato covariant},

		\begin{equation*} 
	    \bm D_t^{\bm \nabla} \bm V_t = \SN \bm \nabla_{\bm V_t} \bm Z_i(\Phi_t(\varphi)) \d W^i_t +  \bm J_0'(\bm V_t) \d t,
	\end{equation*}
	with

	$$
\bm J'_0(\cdot) = \bm \nabla \bm Z_0(\cdot) - \frac{1}{2} \SN \bm \nabla \left(\bm \nabla_{\bm Z_i}\bm Z_i \right)(\cdot) + \frac{1}{2} \SN \left[ \bm \nabla (\bm \nabla \bm Z_i) (\bm Z_i)(\cdot) +\bm \nabla \bm Z_i (\bm \nabla \bm Z_i(\cdot))\right].
	$$
	The conclusion follows from Proposition \ref{P derivee covariante pullback} since

	$$
	\bm \nabla (\bm \nabla \bm Z_i) (\bm Z_i)(\cdot) = \bm \nabla\left(\bm \nabla_{\bm Z_i}\bm Z_i \right)(\cdot)- \bm \nabla \bm Z_i(\bm \nabla \bm Z_i(\cdot)) - \bm R(\cdot,\bm Z_i)\bm Z_i.
	$$
	\end{proof} 

	Assume now that $A \circ \varphi \in \H_\varphi \mathscr{D}$ and let $\bar{V}_t = T_{\M_t(\mu)}p(\bm V_t) = \Pi_{\M_t(\mu)}(\bm V_t \circ \varphi_t^{-1})$. Note that, by construction, $\mathfrak{h}_{\Phi_t(\varphi)}(\bar{V}_t) = P_{\H_{\Phi_t(\varphi)}}(\bm V_t)$. Since for all $x \in M$, $\bm T \bm E_{x}(\bm V_t) \sim w_t(\varphi(x),p(\varphi), A)$, this implies that for all $t \geq 0$ and all $x \in M$, $\mathbb{P}(\bm T \bm E_{x}(\bm V_t) = w_t(\varphi(x),p(\varphi), A))=1$. In particular,
	
	$$
	\mathbb{E}\left[ \int_M |\bm V_t(x)- w_t(\varphi(x),p(\varphi), A)|^2 \d \vol \right] =\int_M  \mathbb{E}\left[ | \bm V_t(x)- w_t(\varphi(x),p(\varphi), A)|^2  \right]\d \vol =0.
	$$ 
	Hence, $\bm V_t = w_t (\varphi(\cdot), p(\varphi), A)$ in $L^2(\Omega,L^2(M,TM;\vol))$ for all $t \geq 0$. Thus, for all $t \geq 0$,

	\begin{equation} \label{E version P}
		\bar{V}_t = \Pi_{\M_t(\mu)}(\bm V_t \circ \varphi_t^{-1}) = \Pi_{\M_t(\mu)}(w_t (\varphi(\cdot), p(\varphi), A) \circ \varphi_t^{-1}), \quad \text{in $L^2(\Omega,L^2(M,TM;\M_t(\mu)))$.}
	\end{equation}
	We will derive both Itô and Stratonovich equations for $(\bar{V}_t)_{t \geq 0}$.
	
	\begin{coro}[Stratonovich-Schwarz Theorem on $\P_\infty$]  \label{Stochastic Schwarz's lemma}
		The process $(\bar{V}_t)_{t \geq 0}$ is solution to the following Stratonovich covariant SDE:
		
		\begin{equation}
		\circ \bar{D}_t \bar{V}_t  = \SN \grad_{\bar{V}_t} \bar{Z}_i(\M_t(\mu)) \circ \d W^i_t + \grad_{\bar{V}_t} \bar{Z}_0(\M_t(\mu)) - \frac{1}{2} \SN \grad_{\bar{V}_t} \grad_{\bar{Z}_i} \bar{Z}_i(\M_t(\mu)) \d t.
		\end{equation}
	\end{coro}

	\begin{proof}
		Let $\bm V_t^\H = P_{\H_{\Phi_t(\varphi)}}(\bm V_t)$. Proposition \ref{Vector bundle chain rule ILH}, Proposition \ref{Strato schwarz lemma D} and Theorem \ref{Dérivées covariantes projection horizontale et verticale} yield
		
		\begin{align*}
			\circ \bm D_t \bm V_t^\H &= P_{\H_{\Phi_t(\varphi)}}\left( \circ \bm D_t \bm V_t -  \mathcal{O}^*_{\circ \dd \Phi_t(\varphi)}(\bm V_t) \right) - \mathcal{O}_{\circ\d \Phi_t(\varphi)}(\bm V_t^\H)  \\
			& = P_{\H_{\Phi_t(\varphi)}}  \left( \bm \nabla_{\bm V_t} \circ \dd \Phi_t(\varphi) - \bm \nabla_{P_{\V_{\Phi_t(\varphi)}}(\bm V_t)} \circ \dd \Phi_t(\varphi) \right) - \mathcal{O}_{\circ\d \Phi_t(\varphi)}(\bm V_t^\H) \\
			& = P_{\H_{\Phi_t(\varphi)}} \left( \bm \nabla_{\bm V_t^\H } \circ \dd \Phi_t(\varphi) \right) - \mathcal{O}_{\circ\d \Phi_t(\varphi)}(\bm V_t^\H).
 		\end{align*}
		The conclusion follows.
	\end{proof}

		\begin{prop} \label{Existence unicité processus tangent}
		The process $(\bar{V}_t)_{t \geq 0}$ is the unique solution to the following Itô's covariant SDE:
		
		\begin{equation}\label{EDS covariante processus tangent}
		\D_t^{\grad} \bar{V}_t = \SN \grad_{\bar{V}_t} \bar{Z}_i(\M_t) \d W^i_t + \left( \grad_{\bar{V}_t} \bar{Z}_0(\M_t)  - \frac{1}{2} \SN \bar{R}(\bar{V}_t,\bar{Z}_i(\M_t)) \bar{Z}_i(\M_t) \right)  \d t
		\end{equation}
		
	\end{prop}

\begin{proof}
	By Proposition \ref{Ito Strato covariant}, Corollary \ref{Stochastic Schwarz's lemma}, Theorem \ref{Dérivées covariantes projection horizontale et verticale} and Corollary \ref{Dérivée covariante Qt}, we obtain:
	
	\begin{align*}
		\bm D_t^{\bm \nabla} \bm V_t^{\H} & \overset{(V)}{=} \SN P_{\H_{\varphi_t}} \left( \bm \nabla_{\bm V_t^\H} \bm Z_i(\varphi_t)\right) - P_{\V_{\varphi_t}} \left( \bm \nabla_{\bm V_t^\H} \bm Z_i(\varphi_t)\right) \d W^i_t + P_{\H_{\varphi_t}} \left( \bm \nabla_{\bm V_t^\H} \bm Z_0(\varphi_t) \right) \d t \\
		& + \frac{1}{2} \SN P_{\H_{\varphi_t}}\left( \bm \nabla_{P_{\H_{\varphi_t}}\left(\bm \nabla_{\bm V_t^\H}\bm Z_i\right)}\bm Z_i\right)  - P_{\H_{\varphi_t}}\left( \bm \nabla_{P_{\V_{\varphi_t}}\left(\bm \nabla_{\bm V_t^\H}\bm Z_i\right)}\bm Z_i\right) \\
		&- \frac{1}{2} \SN P_{\H_{\varphi_t}}\left(\bm \nabla_{\bm V_t^\H} \bm \nabla_{\bm Z_i} \bm Z_i(\varphi_t) + 2\mathcal{O}^*_{\bm Z_i}(\bm \nabla_{\bm V_t^\H}\bm Z_i) \right) \d t \\
		& + \frac{1}{2} \SN P_{\H_{\varphi_t}} \left( \bm \nabla_{\bm V_t^\H} \bm \nabla_{\bm Z_i} \bm Z_i - \bm \nabla_{\bm \nabla_{\bm V_t^\H} \bm Z_i}\bm Z_i - \bm R(\bm V_t^\H, \bm Z_i)\bm Z_i \right) \d t \\
		& \overset{(V)}{=} \SN P_{\H_{\varphi_t}} \left( \bm \nabla_{\bm V_t^\H} \bm Z_i(\varphi_t)\right) - P_{\V_{\varphi_t}} \left( \bm \nabla_{\bm V_t^\H} \bm Z_i(\varphi_t)\right) \d W^i_t + P_{\H_{\varphi_t}} \left( \bm \nabla_{\bm V_t^\H} \bm Z_0(\varphi_t) \right) \d t \\
		& - \frac{1}{2}\SN P_{\H_{\varphi_t}} \left( \bm R(\bm V_t^\H, \bm Z_i)\bm Z_i + 4 \mathcal{O}^*_{\bm Z_i}(\bm \nabla_{\bm V_t^\H}\bm Z_i)  \right),
	\end{align*}
	where $\overset{(V)}{=}$ denotes equality up to vertical terms with finite variations. The last equality comes from the fact the sum of the fourth, fifth, seventh and ninth terms is $-2 \SN P_{\H_{\varphi_t}}\left( \mathcal{O}^\star_{\bm Z_i}(\mathcal{O}_{\bm Z_i}(V_t^{\H})) \right)$ and the sixth and the eigth terms cancel out. We used the fact that, when $\bm Z$ is a horizontal vector field,
	
	$$
	\mathcal{O}_{\bm Z} = P_{\V} \circ \bm \nabla \bm Z \quad \text{and} \quad \mathcal{O}^*_{\bm Z} = \bm \nabla \bm Z \circ P_{\V}
	$$
	The conclusion follows from Corollary \ref{Trace tenseur Riemann sur P}.
	\end{proof}

	\subsection{Itô-Wasserstein curvature condition on \texorpdfstring{$\P$}{P}}
		We consider once again \eqref{EDS principale Ito} with $\bar{Z}_i \in \bar{\Gamma}^3_S(T \P)$ for $i=0,1,\dots,N$ and we introduce the associated Itô-Wasserstein curvature tensor:
	
	\begin{equation}
	\mathfrak{C}(U) = \SN \left( \langle \bar{R}(U, \bar{Z}_i) \bar{Z}_i,U \rangle  - \left| \grad_{U} \bar{Z}_i \right|^2 \right) -2  \langle \grad_{U}\bar{Z}_0 ,U \rangle.
	\end{equation}
	Let $K \in \mathbb{R}$. If
	
	\begin{equation} \tag{IWC(\textit{K})} \label{E IWCD}
	\mathfrak{C}(U) \geq K | U|^2, \quad \forall \mu \in \P_\infty, \, \forall U \in T_\mu \P_\infty,
	\end{equation}
	then we say that the IWC(K) condition holds. In particular, we shall show in this section that this condition entails important stability properties for the considered SDE. 
	
	\begin{rem}
		By Corollary \ref{Trace tenseur Riemann sur P}, we get
		
		$$
		\langle \bar{R}(V_f, \bar{Z}_i) \bar{Z}_i, V_f \rangle_\mu \geq \int_M \langle R(\nabla f,Z_i)Z_i, \nabla f \rangle d\mu.
		$$
		In particular, 

		$$
		\mathfrak{C}(V_f) \geq \left(-2 \sup_{\mu \in \P}|\grad \bar{Z}_0| + \SN m_i - \sup_{\mu \in \P}|\grad \bar{Z}_i|^2 \right)   |V_f|^2
		$$
		with

		$$
		m_i = \inf_{((x, v), \mu) \in U(TM) \times \P} \left\langle R(v,Z_i(\mu)(x))Z_i(\mu)(x), v \right\rangle,
		$$
		where $U(TM)$ denotes the unitary tangent bundle, namely:

		$$
		U(TM) := \bigsqcup_{x \in M} \left\{ v \in T_xM \, ; \, |v|=1 \right\}.
		$$
		In particular, there always exist a constant $K$ such that $\mathrm{IWC}(K)$ holds.
	\end{rem}

	We now state some consequences of the $\mathrm{IWC}(K)$ condition. Let $\mu, \nu \in \P_\infty$ and $(a,x) \in [0,1] \times M \mapsto U_a(x) \in M$ be a smooth path of diffeomorphisms with $U_0 = \mathrm{id}$, $(\mu_a)_{0 \leq a \leq 1} = ((U_a)_*\mu)_{0 \leq a \leq 1}$ such that $\mu^0 = \mu$ and $\mu^1=\nu$. Let $(\varphi_t^a)$ be the field in Theorem \ref{T C1}. 

	By corollary \ref{Coro dérivée processus en a}, since $(\varphi_t^a)$ is an $(\mathscr{F}_t)_{t \geq 0}$-adapted $\C^1$ stochastic field of diffeomorphisms, $(\M_t(\mu_a))_{0 \leq a \leq 1}$ belong to $\mathscr{R}_s(\M_t(\mu^0),\M_t(\mu^1))$. Thus:
	
	$$
	\mathbb{E}[W_2^2(\M_t(\mu),\M_t(\nu))] \leq \mathbb{E} \left[ \int_0^1 |\partial_a \mu_t^a|^2_{\mu^a_t} \d a \right] = \int_0^1 \mathbb{E}\left[|\partial_a \mu_t^a|_{\mu^a_t}^2\right] \d a.
	$$
	By \eqref{E version P} and the fact that $p : \mathscr{D} \rightarrow \P_\infty$ is a Riemannian submersion, the following holds: 

	$$
	\mathbb{E}\left[|\partial_a \mu_t^a|_{\mu^a_t}^2\right] = \mathbb{E}\left[ \int_M \left| \Pi_{\mu^a_t}\left( (\partial_a \varphi_t^a) \circ (\varphi_t^a)^{-1} \right) \circ \varphi_t^a \right|^2 \d \vol \right] = \mathbb{E} \left[ |\bar{V}_t^a|_{\mu_t^a}^2 \right].
	$$
	Thus, we obtained that

	\begin{equation} \label{Ineq Wdeux}
		\mathbb{E}[W_2^2(\M_t(\mu),\M_t(\nu))] \leq \int_0^1 \mathbb{E}\left[|\bar{V}_t^a|_{\M_t(\mu_a)}^2\right] \d a.
	\end{equation}

	\begin{thm} \label{Décroissance exp}
		Let $K \in \mathbb{R}$ such that the $\mathrm{IWC}(K)$ condition is satisfied. Then, for $\mu, \nu \in \P$, the following holds:
		
		\begin{equation} \label{cvge expo}
		\bm W_2^2(\delta_\mu \se_t, \delta_\nu \se_t) \leq\mathbb{E}[W_2^2(\M_t(\mu),\M_t(\nu))] \leq e^{-Kt} W_2^2(\mu, \nu).
		\end{equation}
	\end{thm}
	
	\begin{proof}
		By definition:

		$$
		\bm W_2^2(\delta_\mu \se_t, \delta_\nu \se_t) \leq\mathbb{E}[W_2^2(\M_t(\mu),\M_t(\nu))], \quad \forall \mu,\nu \in \P.
		$$
	We first prove the result for $\mu,\nu \in \P_\infty$. Let $\varepsilon>0$ and $(\mu_a^\varepsilon)_{0 \leq a \leq 1}$ be the $\varepsilon$-entropic interpolation between $\mu$ and $\nu$. We denote by $(\bar{V}_t^{a,\varepsilon})_{t \geq 0}$ the solution to \eqref{Proc Tang} with initial condition $\partial_a \mu_a^\varepsilon$. Then, we have:
    		
    		$$
    		\d \|\bar{V}_t^{a,\varepsilon}\|^2 = 2 \langle \D_t^{\grad} \bar{V}_t^{a,\varepsilon}, \bar{V}_t^{a,\varepsilon} \rangle + \|\D_t^{\grad} \bar{V}_t^{a,\varepsilon}\|^2.
    		$$
    		The $\mathrm{IWC}(K)$ condition yields:
    		
    		$$
    		\d \mathbb{E}[\| \bar{V}_t^{a,\varepsilon}\|^2] = \mathbb{E}[\d \| \bar{V}_t^{a,\varepsilon}\|^2] \leq -K \mathbb{E}\left[ \|\bar{V}_t^{a,\varepsilon}\|^2 \right]\d t.
    		$$
    		By Grönwall's lemma we get
    		
    		$$
    		\mathbb{E}[\|\bar{V}_t^{a,\varepsilon}\|^2] \leq e^{-Kt} \mathbb{E} \left[ \|\partial_a \mu_a^\varepsilon \|^2 \right].
    		$$
    		Now, by \eqref{Ineq Wdeux}, it follows:
    		
    		$$
    		\W_2^2(\delta_\mu \se_t, \delta_\nu \se_t) \leq \mathbb{E}\left[ W_2^2(\M_t(\mu),\M_t(\nu)) \right] \leq \mathbb{E} \left[ \int_0^1 \| \bar{V}_t^{a,\varepsilon}\|^2 da \right] \leq e^{-Kt} \mathbb{E} \left[ \int_0^1 \|\partial_a \mu_a^\varepsilon\|^2 da \right].
    		$$
    		Then, the conclusion follows from Proposition \ref{Convergence entropic interpolation} since

			$$
			\lim_{\varepsilon \to 0} \int_0^1 \| \partial_a \mu^\varepsilon_a \|^2_{\mu^\varepsilon_a} da =W_2^2(\mu,\nu).
			$$
			Thus, since $\M_t : \mu \in \P_\infty \mapsto \M_t(\mu) \in L^2(\Omega, \P; \mathbb{P})$ is uniformly continuous on $\P_\infty$, we can extend the result to the case $\mu,\nu \in \P$ by density of $\P_\infty$ in $\P$.
			
	\end{proof}

    \begin{coro} \label{Coro mesu inv}
	Let $K \in \mathbb{R}$ such that $\mathrm{IWC}(K)$ is satisfied. Then, for all $\vartheta_1,\vartheta_2 \in \P(\P)$, the following holds:

	\begin{equation} \label{Decroissance exp forte}
		\bm W_2^2(\vartheta_1 \se_t, \vartheta_2 \se_t) \leq e^{-Kt} \bm W_2^2(\vartheta_1, \vartheta_2).
	\end{equation}
    Consequently, if $K >0$, there exists a unique $\Lambda \in \P(\P)$ such that

    \begin{equation} \label{Conv mesu inv}
		\lim_{t \rightarrow \infty } \W_2^2(\vartheta_1 \se_t, \Lambda) =0, \quad \forall \vartheta_1 \in \P(\P).
	\end{equation}
    Moreover, $\Lambda$ is the unique invariant measure in $\P(\P)$ for the semigroup $(\se_t)_{t \geq 0}$.
    \end{coro}
    
    \begin{proof}
		Let $\Pi$ be an optimal coupling of $\vartheta_1$ and $\vartheta_2$ and $\bar{Q}_t$ denote the transition semigroup associated to the Markov process $(\mathcal{Z}_t)$ defined by $\mathcal{Z}_t(\mu,\nu)=(\M_t(\mu),\M_t(\nu))$. Then, $\Pi \bar{Q}_t$ is a coupling of $\mu \se_t$ and $\nu \se_t$. By Proposition \ref{Décroissance exp},

		\begin{align*}
			\bm W_2^2(\vartheta_1 \se_t, \vartheta_2 \se_t) &\leq \int_{\P} W_2^2(\mu,\nu) \d (\Pi \bar{Q}_t)(\mu,\nu) \\
			& = \int_{\P} \mathbb{E}\left[W_2^2(\M_t(\mu), \M_t(\nu)) \right] \d \Pi(\mu,\nu) \\
			& \leq e^{-Kt} \int_{\P} W_2^2(\mu,\nu) \d \Pi(\mu,\nu) = e^{-Kt} \bm W_2^2(\vartheta_1,\vartheta_2).
		\end{align*}

    		The existence of an ergodic invariant measure follows from the semigroup property $\se_{t+s} = \se_s \se_t$. Indeed, by \eqref{Decroissance exp forte}, we obtain:
    		
    		$$
    		\W_2^2(\vartheta_1 \se_{t+s},\vartheta_1 \se_{t}) \leq e^{-Kt} \bm W_2^2(\vartheta_1 \se_s, \vartheta_1) \leq e^{-Kt} \mathrm{Diam}(M)^2.
    		$$
    	We conclude by using the fact that $\P(\P)$ is a complete space. The fact that $\Lambda$ does not depend on $\vartheta_1$ is a straightforward consequence of \eqref{Decroissance exp forte}. To see that $\Lambda$ is invariant for the semigroup $(\se_t)_{t \geq 0}$, note that for any $\vartheta \in \P(\P)$,

		$$
		\vartheta \se_{t+s} = (\vartheta \se_s) \se_t.
		$$
		By \eqref{Conv mesu inv}, $\lim_{s \to \infty} \vartheta \se_{t+s} = \Lambda$, and by Proposition \ref{Semigroup property}, $\lim_{s \to \infty}(\vartheta \se_s) \se_t = \Lambda \se_t $. The invariance property follows.

    	Assume that $\Lambda' \in \P(\P)$ is another invariant measure. Then, 

		$$
		\W_2^2( \Lambda'\se_t, \Lambda) = \W_2^2(\Lambda', \Lambda).
		$$
		Hence, by \eqref{Conv mesu inv}, $\Lambda'= \Lambda$.
    \end{proof}
    
	\begin{prop} \label{Sub-Commutativity property}
		Let $K \in \mathbb{R}$ such that the $\mathrm{IWC}(K)$ condition is satisfied and let $F \in \C^1(\P)$. Then, $\se_t F \in \C^1(\P)$ and the following sub-commutativity inequality holds:
		
		$$
		|\grad \se_t F(\mu)|_{\mu}^2 \leq e^{-Kt} \se_t \left( |\grad F(\mu)|_\mu^2 \right), \quad \forall \mu \in \P.
		$$
	\end{prop}
	
	\begin{proof}
	The fact that $\se_t F \in \C^1(\P)$ has been proved in Proposition \ref{Stab semigroup}. Let $\mu \in \P_\infty$, $V \in T_\mu \P_\infty$ with $|V|_\mu = 1$. Let $\varepsilon> 0$ such that for $0 \leq a < \varepsilon$, the path $(U_a)_{0 \leq a < \varepsilon} = (\exp(aV))_{0 \leq a < \varepsilon}$ is a $\C^\infty$ path of diffeomorphisms. Then, for $a \in [0, \varepsilon)$, the path $(\mu_a)_{0 \leq a \leq \varepsilon}$ defined by $\mu_a = (U_a)_* \mu$ belongs to $\mathscr{R}(\mu,\nu)$. Let also $\varphi_t^a = \varphi_t(U_a,\mu_a)$ be the version of the ID flow which is a $\C^1$ stochastic field of diffeomorphisms by Theorem \ref{T C1} and $\mu_t^a= (\varphi_t^a)_*\mu$. Applying Corollary \ref{Lemme fondamental L2} yields for $0 \leq a <\varepsilon$:
	
	$$
	\partial_a F(\mu_t^a(\omega)) = \int_M \langle \grad F(\mu_t^a(\omega))(\varphi_t^a(x,\mu_a,\omega)), \partial_a \varphi_t^a(x,\mu_a,\omega) \rangle \d \mu, \quad \text{a.s.}.
	$$
	By \eqref{E version P}, $\partial_a \mu_t^a= \bar{V}^a_t$ in $L^2(\Omega,L^2(M,TM;\mu^a_t))$ and, by following the same reasoning as in the proof of Theorem \ref{Décroissance exp}, one can show that $\mathbb{E} \left[|\partial_a \mu^a_t|_{\mu_t^a}^2 \right] = \mathbb{E}[|\bar{V}^a_t|_{\mu_t^a}] \leq e^{-Kt} |\bar{V}^a|_{\mu_a}^2$, where $V^a = \partial_a U_a \circ U_a^{-1}$. Thus, by Cauchy-Schwarz inequality, we obtain:

	\begin{align*}
		|\se_t F(\mu) - \se_t F(\mu_a)|^2 & \leq \mathbb{E} \left[ \int_0^a \partial_r F(\mu_t^r) \d r \right]^2 \\
		&= \left(\int_0^a \mathbb{E} \left[ \int_M \langle \grad F(\mu_t^r)(\varphi_t^r(x,\mu_r)), \partial_r \varphi_t^r(x,\mu_r) \rangle \d \mu  \right]\d r \right)^2 \\
		&= \left(\int_0^a \mathbb{E} \left[ \int_M \left\langle \grad F(\mu_t^r)(x), \bar{V}_t^r \right\rangle \d \M_t(\mu_r)(x)  \right]\d r \right)^2 \\
		& \leq \left( \int_0^a  \mathbb{E} \left[|\grad F(\mu_t^r)|_{\mu_t^r}^2 \right]^{1/2} \mathbb{E}\left[ |\bar{V}_t^r|_{\mu_t^r}^2 \right]^{1/2} \d r \right)^2 \\
		& \leq e^{-Kt} \left(\int_0^a  \mathbb{E} \left[|\grad F(\mu_t^r)|_{\mu_t^r}^2 \right]^{1/2} \cdot |\bar{V}^r|_{\mu^r} \d r\right)^2.
	\end{align*}
	Dividing by $a^2$ and taking the supremum in $ 0 \leq s \leq a $ yields:
	
	$$
	\frac{|\se_t F(\mu_a) - \se_t F(\mu)|^2}{a^2} \leq e^{-Kt} \sup_{0 \leq s \leq a} \mathbb{E} \left[ |\grad F(\mu^s_t)|_{\mu^s_t}^2 \right] \frac{1}{a^2}\left(\int_0^a  |V^r|_{\mu^r}  \d r\right)^2,
	$$
	and, by letting $a \rightarrow 0$, we finally get:
	
	$$
	|T_\mu \se_t F(V)|^2 \leq e^{-Kt} \mathbb{E} \left[ |\grad F(\mu_t) |_{\mu_t}^2 \right] \cdot |V|_\mu^2 = e^{-Kt} \se_t \left(|\grad F(\mu)|_\mu^2 \right).
	$$
	Thus, we obtained

	$$
		|\grad \se_t F(\mu)|_{\mu}^2 \leq e^{-Kt} \se_t \left( |\grad F(\mu)|_\mu^2 \right), \quad \forall \mu \in \P_\infty.
	$$
	We obtain the result for all $\mu \in \P$ by using both the density of $\P_\infty$ in $\P$ and the fact that $|\grad \se_t F(\cdot)|$ and $\se_t \left( |\grad F(\cdot)| \right)$ are continuous mapping by Proposition \ref{Stab semigroup}.
	\end{proof} 

		\begin{prop}[Local Poincaré Inequality] \label{Local Poincare infini}
		Let us assume that $\mathrm{IWC}(K)$ holds for some constant $K \in \mathbb{R}$ and let $F \in \C^2(\P)$. Then, the following local Poincaré inequality is satisfied:
		
		\begin{equation*}
		\se_t (F^2)(\mu) - (\se_t F)^2(\mu) \leq 2|\bar{Z}|_{2,\infty} C(K,t) \se_t ( | \grad F |^2)(\mu),\quad \forall t \geq 0, \, \mu \in \P,
		\end{equation*}
		where $|\bar{Z}|^2_{2,\infty} = \sup_{\mu \in \P} \SN |\bar{Z}_i(\mu)|_\mu^2$ and $C(K,t)$ is defined in \eqref{ckt}:

		$$
C(K,t) = \begin{cases}
	\frac{1-e^{-Kt}}{K} \text{ if } K \neq 0, \\
	t \text{ if } K =0.
\end{cases}
		$$
	\end{prop}
	
	\begin{proof}
	The proof is just an adaptation of \cite[Proposition 5.4.1]{logsob2000} to our infinite dimensional context. Let $t > 0$ and $\Psi : [0,t] \rightarrow \mathbb{R}$ the function defined as
	
	$$
	\Psi_s = \se_s\left((\se_{t-s} F)^2\right).
	$$
	Note that $\Psi \in \C^{1,2}(\mathbb{R}_+ \times \P)$ by Proposition \ref{Stab semigroup}.
	By applying the chain rule to $\Psi$, \eqref{E Dérivée semi groupe} and the definition of $\bar{\Gamma}$ we obtain: 
	
	$$
	\dot{\Psi}_s = \se_s  \left(  \bar{\mathrm{L}}\left( \se_{t-s} F \right)^2 -  2 (\se_{t-s} F) \bar{\mathrm{L}}(\se_{t-s} F) \right) = 2 \se_s \left( \bar{\Gamma}(\se_{t-s} F) \right).
	$$
	Hence, 
	
$$
	\se_t (F^2) - \left(\se_t F \right)^2 = \Psi_t - \Psi_0 = 2 \int_0^t \se_s \left( \Gamma \left( \se_{t-s} F \right) \right) \d s.
$$
	Now, by applying Proposition \ref{Sub-Commutativity property} to the above inequality, by setting $r=t-s$ and by using the semigroup property (Proposition \ref{Semigroup property}), we get:
	
	$$
	2\int_0^t \se_s \left( \Gamma \left( \se_{t-s} F \right) \right) \d s \leq 2\int_0^t e^{-Kr}\se_t (\Gamma(F)) \d r.
	$$
	Since $\Gamma(F) \leq |\bar{Z}|_{2,\infty}^2 | \grad F|^2$, we obtain finally:
	
	$$
	2\int_0^t e^{-Kr}\se_t (\Gamma(F)) \leq \left(\int_0^t e^{-Kr} \d r \right) 2 |\bar{Z}|_{2,\infty}^2 \se_t (|\grad F|^2) = 2|Z|_{2,\infty}^2 \frac{1- e^{-Kt}}{K} \se_t (|\grad F|^2).
	$$
	Consequently, we obtain
	
	$$
	\se_t (F^2) - (\se_t F)^2 \leq 2|\bar{Z}|_{2,\infty}^2\frac{1- e^{-Kt}}{K} \se_t ( | \grad F |^2).
	$$
	\end{proof}

\subsection{Itô-Wasserstein curvature condition on \texorpdfstring{$M \times \P$}{MP}}

	In the following we consider $\P(M \times \P)$ endowed with the Wasserstein-2 metric, namely,

	$$
	\mathcal{W}_2^2(\check{\vartheta}_1, \check{\vartheta}_2) = \inf \left\{ \int_{(M \times \P)^2} d^2(x,y) + W_2^2(\mu,\nu) \, \d \check{\Pi}((x,\mu),(y,\nu)), \, \check{\Pi} \in \C(\check{\vartheta}_1,\check{\vartheta}_2) \right\}, \quad \check{\vartheta}_1, \check{\vartheta}_2 \in \P(M \times \P).
	$$
	Let $(\varphi_t)_{t \geq 0}$ denote the unique solution flow to the following Image-dependent SDE

	\begin{equation} \label{MCV2}
		\begin{cases}
		\d^\nabla \varphi_t(x,\mu) = \SN Z_i(\varphi_t(x,\mu),\M_t(\mu)) \d W^i_t + Z_0(\varphi_t(x,\mu), \M_t(\mu)) \d t, \\
		\varphi_0(\cdot,\mu) = \mathrm{id}, \, \M_t(\mu)= (\varphi_t(\cdot,\mu))_* \mu.
	\end{cases}
	\end{equation}
	and consider the following operator:

	\begin{align*}
	\mathfrak{C}'((U,V)) &=\SN \left\langle R(U, A_i(x,\mu))A_i(x,\mu), U \right\rangle - \left| \nabla_{U} A_i(x,\mu)  + \bar{\mathcal{L}}_{V} A_i(x,\mu) \right|^2  \\
		&-2\left\langle \nabla_{U} A_0(x,\mu) + \bar{\mathcal{L}}_{V} A_0(x,\mu), U \right\rangle.
	\end{align*}
	Let $K_1 \in \mathbb{R}, \, K_2 \leq 0$. If
	
	\begin{equation} \tag{$\mathrm{IWC}_{M}$($K_1,K_2$)} \label{E IWCM}
	\mathfrak{C}'((U,V)) \geq K_1 |U|_x^2 + K_2|V|_\mu^2, \quad \forall x\in M, \, \mu \in \P_\infty, \, U \in T_xM, \, V \in T_\mu \P,
	\end{equation}
	then we say that the $ \mathrm{IWC}_M(K_1,K_2)$ condition holds. Note that $\mathrm{IWC}_M(K_1,K_2)$ is strictly stronger than $\mathrm{IWC}(K_1+K_2)$. This follows by setting $U=V(x)$, integrating \eqref{E IWCM} with respect to $\mu$ and by noticing that

	$$
	\begin{cases}
		\int_M \langle R(V(x),A_i(x,\mu))A_i(x,\mu), V(x) \rangle \d \mu(x) \leq \langle \bar{R}(V,A_i)A_i, V \rangle_\mu \\
		 |\grad_{V}A_i|^2_\mu \leq \int_M \left| \nabla_{V(x)} A_i(x,\mu)  + \bar{\mathcal{L}}_{V} A_i(x,\mu) \right| \d \mu(x).
	\end{cases}
	$$

	\begin{rem}
		This condition is exactly the infinitesimal form of the generalization to manifolds of the condition introduced by Wang in \cite[Section 4]{wang2020}.
	\end{rem}

	\begin{rem}
		The $ \mathrm{IWC}_M(K_1,K_2)$ condition may look surprising at first since we have a cross term on the left hand-side, $\left\langle \bar{\mathcal{L}}_{V} A_0(x,\mu), U \right\rangle$ and no cross term on the right-hand side. Actually, we could reformulate $ \mathrm{IWC}_M(K_1,K_2)$ by introducing a third constant $K_3 <0$ and asking for

		$$
		\mathfrak{C}'((U,V)) \geq K_1 |U|_x^2 + K_2|V|_\mu^2 + K_3 |U|_x |V|_\mu, \quad \forall x\in M, \, \mu \in \P_\infty, \, U \in T_xM, \, V \in T_\mu \P.
		$$
		Since $K_3 <0$, by Young's inequality, for all $\varepsilon >0$, $ K_3 |U|_x |V|_\mu \geq K_3 \left( \frac{|U|_x^2}{2\varepsilon} + \frac{\varepsilon |V|_\mu^2}{2} \right)$. Obviously, this constant is unnecessary and we will continue with the $ \mathrm{IWC}_M(K_1,K_2)$ condition.
	\end{rem}

	\begin{rem}
		Assume that the $ \mathrm{IWC}_M(K_1,K_2)$ condition holds for $K_1 \in \mathbb{R}$ and $K_2 \geq 0$. For $U =0$ and $0 \neq V \in T_\mu \P$, this means that

		$$
		- \left| \bar{\mathcal{L}}_{V} A_i(x,\mu) \right|^2 \geq K_2 |V|_\mu^2.
		$$
		This inequality can be true if and only if $A_i$ is constant in $\mu$ for $i=1,\dots,N$ and $K_2=0$. Moreover, let $ U \in T_xM$ and $V \in T_\mu \P$, this also means that 

		$$
		\left\langle R(U, A_i(x,\mu))A_i(x,\mu), U \right\rangle -2\left\langle \nabla_{U} A_0(x,\mu), U \right\rangle - K_1 |U|^2 \geq   2 \left\langle \bar{\mathcal{L}}_{V} A_0(x,\mu),U \right\rangle.
		$$
		Since $V \in T_\mu \P \mapsto \left\langle \bar{\mathcal{L}}_{V} A_0(x,\mu),U \right\rangle$ is a linear map, it is constant equal to $0$, in particular, $A_0$ is also constant in $\mu$. In conclusion, $K_2$ has to be nonnegative and, if $K_2=0$, the driving vector fields $A_0,A_1,\dots,A_N$ are constant in $\mu$. In the case $K_2=0$, one should apply the finite-dimensional methods developped in Subsection \ref{S 1.3}.
	\end{rem}

	\begin{thm} \label{T Ergo MP}
		Assume that $\mathrm{IWC}_M(K_1,K_2)$ holds for $K_1 \in \mathbb{R},\, K_2 < 0$. Then: \begin{enumerate}
			\item For all $(x,\mu), (y,\nu) \in \M \times \P$,

		$$
		\mathbb{E}\left[ d^2(\varphi_t(x,\mu), \varphi_t(y,\nu)) \right] \leq e^{-K_1t} d^2(x,y) + e^{-(K_1+K_2)t}W_2^2(\mu,\nu)
		$$
		\item For all $\check{\vartheta}_1,\check{\vartheta}_2 \in \P(M \times \P)$, the following holds:

	\begin{equation*} 
		\mathcal{W}_2^2(\check{\vartheta}_1 \check{\mathrm{P}}_t, \check{\vartheta}_2 \check{\mathrm{P}}_t) \leq 2e^{-(K_1+K_2)t} \mathcal{W}_2^2(\check{\vartheta}_1, \check{\vartheta}_2).
	\end{equation*}
	If $K_1+K_2 >0$, there exists a unique probability measure $\check{\Lambda} \in \P(M \times \P)$ such that

    \begin{equation*}
		\lim_{t \rightarrow \infty } \mathcal{W}_2^2(\check{\vartheta} \check{\mathrm{P}}_t, \check{\Lambda}) =0, \quad \forall \check{\vartheta} \in \P(M \times \P).
	\end{equation*}
    Moreover, $\check{\Lambda}$ is the unique invariant measure in $\P(M \times \P)$ for the semigroup $(\check{\mathrm{P}}_t)_{t \geq 0}$.
		\end{enumerate}
		
	\end{thm}

	\begin{proof}
		Let $x,y \in M$, $\mu,\nu \in \P_\infty$, let us show the first point for $\check{\vartheta}_1 = \delta_{(x,\mu)}$ and $\check{\vartheta_2}=\delta_{(y,\nu)}$. Let $(\gamma_a)_{0 \leq a \leq 1}$ be a geodesic from $x$ to $y$, consider $Z$ a smooth extension of $(\dot{\gamma}_a)_{0 \leq a \leq 1}$ and

		$$
		\partial_a U_a = Z(U_a), \quad U_0 = \mathrm{id}.
		$$
		Let $\varepsilon >0$ and $(\mu_a^\varepsilon)_{0 \leq a \leq 1}$ be the $\varepsilon$-entropic interpolation from $\mu$ to $\nu$. By Proposition \ref{P general flow lisse}, $(a,t,x) \mapsto \varphi_t(U_a(x),\mu_a^\varepsilon)$ is $\C^{1,0,\infty}([0,1] \times \mathbb{R}_+\times M,M)$ and $\mathcal{V}_{t}^{a,\varepsilon}=(\partial_a \varphi_t(U_a(x),\mu^\varepsilon_a))_{t \geq 0}$ is the unique solution to 

		\begin{align*}
			D^\nabla_t \mathcal{V}_{t}^{a,\varepsilon} &= \SN \left( \nabla_{\mathcal{V}_{t}^{a,\varepsilon}} Z_i(U_a(x),\mu^\varepsilon_a) \d W^i_t - \frac{1}{2} R\left(\mathcal{V}_{t}^{a,\varepsilon},Z_i(U_a(x),\mu^\varepsilon_a)\right)Z_i(U_a(x),\mu^\varepsilon_a) \d t\right) + \nabla_{\mathcal{V}_{t}^{a,\varepsilon}} Z_0(U_a(x),\mu^\varepsilon_a) \d t \\
			& + \SN \bar{\mathcal{L}}_{\T_{\mu_a^\varepsilon}\M_t(\partial_a \mu^\varepsilon_a)} Z_i(U_a(x),\mu^\varepsilon_a) \d W^i_t + \bar{\mathcal{L}}_{\T_{\mu_a^\varepsilon}\M_t(\partial_a \mu^\varepsilon_a)} Z_0(U_a(x),\mu^\varepsilon_a) \d t.
		\end{align*}
		Thus, by $\mathrm{IWC}_M(K_1,K_2)$,

		$$
		\d \mathbb{E} \left[ |\mathcal{V}_{t}^{a,\varepsilon}|^2 \right] \leq  -K_1 \mathbb{E}\left[|\mathcal{V}_{t}^{a,\varepsilon}|^2\right] -K_2 \mathbb{E} \left[ |\T_{\mu_a^\varepsilon}\M_t(\partial_a \mu^\varepsilon_a)|^2 \right] \d t.
		$$
		By Grönwall's lemma, the fact that $\mathrm{IWC}_M(K_1,K_2)$ implies $\mathrm{IWC}(K_1+K_2)$ and by applying the same reasoning as in the proof of Theorem \ref{Décroissance exp}, we obtain

		\begin{align*}
			\mathbb{E} \left[ |\mathcal{V}_{t}^{a,\varepsilon}|^2 \right] &\leq e^{-K_1 t} |\dot{\gamma}_a|^2_{\gamma_a} -K_2 \int_0^t e^{-K_1(t-s)} \mathbb{E} \left[ |\T_{\mu_a^\varepsilon}\M_s(\partial_a \mu^\varepsilon_a)|^2 \right] \d s \\
			& \leq e^{-K_1 t} |\dot{\gamma}_a|^2_{\gamma_a} - |\partial_a \mu^\varepsilon_a|^2_{\mu_a} e^{-K_1 t} \int_0^t K_2 e^{-K_2s} \d s \\
			& \leq e^{-K_1 t} |\dot{\gamma_a}|^2_{\gamma_a} + e^{-(K_1+K_2)t} |\partial_a \mu^\varepsilon_a|^2_{\mu_a}.
		\end{align*}
		Thus, we obtained:

		$$
		\mathbb{E}[d^2(\varphi_t(x,\mu),\varphi_t(y,\nu))] \leq \int_0^1 \mathbb{E} \left[ |\mathcal{V}_{t}^{a,\varepsilon}|^2 \right] \d a \leq e^{-K_1 t} d^2(x,y) + e^{-(K_1+K_2)t} \int_0^1 |\partial_a \mu^\varepsilon_a|^2_{\mu_a^\varepsilon} \d a.
		$$
		The result for $\mu, \nu \in \P_\infty$ follows by letting $\varepsilon \to 0$, the general result is obtained by uniform continuity on a dense subspace once again.

		By using the fact that the couple $(\varphi_t(\cdot,\cdot), \M_t(\cdot))$ is a strong Markov process, we can repeat the same arguments as in the proof of Theorem \ref{Décroissance exp} to prove the first point. 
		
		We now prove the second point, let $\check{\Pi}$ be an optimal coupling be an optimal coupling of $\check{\vartheta}_1$ and $\check{\vartheta}_2$. Let also $\check{\mathcal{Z}}_t$ denote the synchronous coupling defined by $\check{\mathcal{Z}}_t((x,\mu),(y,\nu))= \left((\varphi_t(x,\mu),\M_t(\mu)), (\varphi_t(y,\nu),\M_t(\nu))\right)$ and consider the associated strong Markov semigroup $\check{\mathrm{Q}}_t$. Then, $\check{\Pi} \check{\mathcal{Z}}_t$ is a coupling of $\check{\vartheta}_1 \check{\mathrm{P}}_t$ and $\check{\vartheta}_2 \check{\mathrm{P}}_t$. It follows from the first point and from Theorem \ref{Décroissance exp} that:

		\begin{align*}
			\mathcal{W}_2^2(\check{\vartheta}_1 \check{\mathrm{P}}_t,\check{\vartheta}_2 \check{\mathrm{P}}_t) & \leq \int_{(M \times \P)^2} \mathbb{E} \left[ d^2(\varphi_t(x,\mu), \varphi_t(y,\nu)) + W_2^2(\M_t(\mu),\M_t(\nu)) \right] \d \check{\Pi}((x,\mu),(y,\nu)) \\
			& \leq 2e^{-K_3 t} \int_{(M \times \P)^2} d^2(x,y) + W_2^2(\mu,\nu) \d \check{\Pi}((x,\mu),(y,\nu)) \\
			& = 2e^{-K_3t} \mathcal{W}_2^2(\check{\vartheta}_1,\check{\vartheta}_2).
		\end{align*}
		The remainder of the proof can be performed in the exact same way as in the proof of Corollary \ref{Coro mesu inv}.
	\end{proof}

	Let $\langle \bullet, \bullet \rangle$ also denotes the product Riemannian metric on $M \times \P$, that is, 

	$$
	\langle (u,U), (v,V) \rangle_{(x,\mu)} = \langle u, v \rangle_x + \langle U,V \rangle_\mu, \quad \forall (u,U), (v,V) \in T_x M \times \T_\mu \P.
	$$
	For $\check{F} \in \C^{1,1}(M \times \P)$, there exists a unique gradient $\check{\nabla} \check{F}$ such that 

	$$
	\langle \check{\nabla}\check{F}(x,\mu), (u,U) \rangle = T_x F(\cdot,\mu)(u) + \T_\mu F(x,\cdot)(U).
	$$
	Moreover, by Proposition \ref{Stab semigroup}, it is clear that, if $\check{F} \in \C^{k,k}(M \times \P)$, then, $\check{\mathrm{P}}_t \check{F} \in \C^{k,k}(M \times \P)$.

	\begin{prop} Assume that $\mathrm{IWC}_M(K_1,K_2)$ holds for some constants $K_1 \in \mathbb{R}, \, K_2 < 0$. Then,
		\begin{enumerate}
			\item Let $\check{F} \in \C^{1,1}(M \times \P)$. Then, the following sub-commutativity inequality holds:
		
		$$
		|\check{\nabla} \check{\mathrm{P}}_t\check{F}(x,\mu)|_{(x,\mu)}^2 \leq e^{-(K_1+K_2)t} \check{\mathrm{P}}_t \left( |\check{\nabla} \check{F}(x,\mu)|_{(x,\mu)}^2 \right), \quad \forall (x,\mu)\in M \times \P.
		$$
		\item Let $\check{F} \in \C^{2,2}(M \times \P)$. Then, the following local Poincaré inequality is satisfied:
		
		\begin{equation*}
		\check{\mathrm{P}}_t (\check{F}^2)(x,\mu) - (\check{\mathrm{P}}_t \check{F})^2(x,\mu) \leq 2|(Z,\bar{Z})|_\infty C(K,t) \check{\mathrm{P}}_t ( | \check{\nabla} \check{F} |^2)(x,\mu),\quad \forall t \geq 0, \, (x,\mu) \in M \times \P,
		\end{equation*}
		where $|(Z,\bar{Z})|^2_{2,\infty} = \sup_{(x,\mu) \in M \times \P} \SN |Z_i(x,\mu)|^2_x + |\bar{Z}_i(\mu)|_\mu^2$ and $C(K,t)$ is defined in \eqref{ckt}.
		\end{enumerate}
	\end{prop}

	\begin{proof}
		The idea of the proof of the first point is the same as in the one of Proposition \ref{Sub-Commutativity property}, just consider the path $[(\gamma_a,\mu_a)]_{0 \leq a < \varepsilon} = [(\exp_x(av), (\exp_x(aV))_*\mu]_{0 \leq a < \varepsilon}$ for $\varepsilon>0$ sufficiently small instead of the path $(\mu_a)_{0 \leq a < \varepsilon}$. The proof of the second point is exactly the same as the one of Proposition \ref{Local Poincare infini} replacing the overlined objects by the check ones.
	\end{proof}
	\subsection{The Image-dependent Langevin case}

	Let $\phi=(\phi_1,\dots,\phi_N) : M \rightarrow \mathbb{R}^N$ be an isometric embedding, $\phi_0 \in \C^\infty_S(M \times \P)$ and denote by $\mathrm{Cut}$ the cut-locus of $M$ and for any $k >0$, 
	
	$$
	\mathrm{Cut}_k = \left\{ (x,y) \in M^2 ,\, d[(x,y),\mathrm{Cut}] < 1/k \right\}.
	$$
	Let $h_{n,\varepsilon} \in \C^\infty(M \times M)$ such that $0 \leq h_{n,\varepsilon} \leq \sqrt{1-\varepsilon}$ and

	$$
	h_{n,\varepsilon}(x,y) = \begin{cases}
		\sqrt{1-\varepsilon} \, \, \text{if } (x,y) \in (M \times M) \setminus \mathrm{Cut}_n; \\
		0 \, \, \text{if } (x,y) \in \mathrm{Cut}_{2n}
	\end{cases}
	$$
	The goal of this subsection is to construct a coupling between two stochastic processes on $M \times \P$ based on the idea of the parallel coupling as in \cite[Subsection 2.1]{Wang2005}. We consider the following $M \times M$-valued Image-dependent SDE:

	\begin{equation} \label{Coupl para}
		\begin{cases}
			\d^\nabla \mathcal{Z}_t^{n,\varepsilon} [(x,y),\pi] = \SN \mathbb{A}_i^{1,n,\varepsilon}(\mathcal{Z}^{n,\varepsilon}_t[(x,y),\pi]) \d W^i_t + \mathbb{A}_i^{2,n,\varepsilon} (\mathcal{Z}^{n,\varepsilon}_t[(x,y),\pi]) \d \tilde{W}^i_t + \mathbb{A}_0(\mathcal{Z}_t[(x,y),\pi], \pi^{n,\varepsilon}_t) \d t, \\
			\pi^{n,\varepsilon}_t = \mathcal{Z}^{n,\varepsilon}_t[(\cdot,\cdot),\pi]_* \pi, \, \mathcal{Z}^{n,\varepsilon}_0[(x,y),\pi]=(x,y),
		\end{cases}
	\end{equation}
	where $W_t$, $\tilde{W}_t$ are two independent $\mathbb{R}^N$-valued Brownian motions and

	$$
	\begin{cases}
		\mathbb{A}_i^{1,n,\varepsilon}[(x,y)] = \left[ \nabla \phi_i(x), h_{n,\varepsilon} (x,y) P_{x,y} \nabla \phi_i(x) \right]; \\
		\mathbb{A}_i^{2,n,\varepsilon}[(x,y)] = \left[ 0, \sqrt{1 - h_{n,\varepsilon}^2 (x,y)} \nabla \phi_i(y) \right]; \\
		\mathbb{A}_0 ([(x,y),\pi]) = \left[ \nabla \phi_0(x,\pi^1), \nabla \phi_0(y,\pi^2) \right],
	\end{cases}
	$$
	where $\pi^i$ denotes the $i$-th marginal of $\pi$. Note that \eqref{Coupl para} is an Image-dependent SDE driven by smooth vector fields on $M^2 \times \P(M^2)$ and thus admits a unique solution by Theorem \ref{Existence Unicité}. We denote by $\mathcal{Z}_t^{i,n,\varepsilon}[(x,y),\pi]$ the projection on the $i$-th coordinate of $\mathcal{Z}_t^{n,\varepsilon}[(x,y),\pi]$ and by $\pi_t^{i,n,\varepsilon}$ the $i$-th marginal of $\pi_t^{n,\varepsilon}$. The idea of this subsection is to take the limit when $(n,\varepsilon) \to (\infty,0)$ to obtain two stochastic processes driven by a (different) Brownian noise. The difference with the finite-dimensional parallel coupling is that the drifts are not the same because they depend on the measure $\pi_t^{i}$ for $i=1,2$. However, we can still recover a Bakry-Emery condition for this kind of SDE, see \eqref{Bakry emery}. 

	\begin{rem}
		Note that $\mathcal{Z}_t^{1,n,\varepsilon}[(x,y),\pi]$, and thus $\pi_t^{1,n,\varepsilon}$, does not depend on $n,\varepsilon,y$ and only depends on $\pi$ through $\pi^1$.
	\end{rem}
	
	Denote by $\mathcal{K}^{n,\varepsilon}_t$ the transition kernel
	
	\begin{align*}
		\mathcal{K}^{n,\varepsilon}_t : M \times \mathscr{B}(M) &\longrightarrow [0,1] \\
		(x,A) &\longmapsto \left[ \left( \mathcal{Z}_t^{2,n,\varepsilon}[(x,\cdot),\pi] \right)_* \pi_x \right](A)
	\end{align*}
	By measure disintegration it follows that 

	$$
		\pi_t^{2,n,\varepsilon} = \int_M \mathcal{K}_t^{n,\varepsilon}(x,\cdot) \, \d \pi^{1}(x).
	$$
	Consider the space $M^2 \times \P(M^2)$ endowed with the distance 

	$$
	\check{d}\left([(x_1,y_1),\pi_1],[(x_2,y_2),\pi_2] \right) = \left( d^2[(x_1,y_1),(x_2,y_2)] + W_2^2(\pi_1,\pi_2) \right)^{1/2}.
	$$
	This metric space is a compact Polish space. Let $[(x,y),\pi] \in M^2 \times \P(M^2)$, $T>0$,
	
	$$
	\Omega' = \C([0,T], M^2 \times \P(M^2)),
	$$
	$\mathscr{F}' = \mathscr{B}(\Omega')$ and consider $\mathbb{Q}^{n,\varepsilon}_{[(x,y),\pi]}$ the law of $\left(\mathcal{Z}^{n,\varepsilon}_\cdot[(x,y),\pi], \pi_\cdot^{n,\varepsilon}\right)$. Denote by $[(X_t,\tilde{X}_t),\sigma_t]_{0 \leq t \leq T}$ the coordinate process on $\Omega'$ and by $\sigma^i_t$ the $i$-th marginal of $\sigma_t$ for $i=1,2$. We will need the following result:

	\begin{thm}[\cite{kallenberg2002}, Corollary 14.9] \label{KCTC}
		Let $T>0$ and $(E,d)$ be a Polish space. Let $\left\{ (X_t^n)_{0 \leq t \leq T} \right\}_{n \in \mathbb{N}}$ in $\C([0,T],E)^\mathbb{N}$ and denote by $\left\{\mathbb{Q}_n\right\}_{n \in \mathbb{N}}$ the associated sequence of probability laws. Assume that $(X_0^n)_{n \in \mathbb{N}}$ is tight in $E$ and that there exists $C,b >0$, $a> 1+b$ such that 

		$$
		\mathbb{E} \left[ d^a(X_t^n, X_s^n) \right] \leq C |t-s|^{1+b}.
		$$
		Then, $\left\{\mathbb{Q}_n\right\}_{n \in \mathbb{N}}$ is tight.
	\end{thm}

	\begin{lem} 
		For all $p \geq 1$, there exists $C_p(T)>0$ such that for all $s,t \in [0,T]$ and all $(n,\varepsilon) \in \mathbb{N} \times \mathbb{R}_+^*$:

		\begin{equation} \label{KCT}
			\mathbb{E} \left[ d^{2p} \left( \mathcal{Z}^{n,\varepsilon}_t[(x,y),\pi],\mathcal{Z}^{n,\varepsilon}_s[(x,y),\pi] \right) \right] \leq C_p(T) |t-s|^p , \quad \forall [(x,y),\pi] \in M^2 \times \P(M^2).
		\end{equation}
		As a consequence, $\left\{\mathbb{Q}^{n,\varepsilon}_{[(x,y),\pi]}\right\}_{(n,\varepsilon) \in \mathbb{N}\times \mathbb{R}_+^*}$ is a tight family of measures on $\Omega'$. 
	\end{lem}
	
	\begin{proof}
		Inequality \eqref{KCT} follows from Itô's formula applied to $|\mathcal{Z}^{n,\varepsilon}_t[(x,y),\pi] - \mathcal{Z}^{n,\varepsilon}_s[(x,y),\pi]|^2$ seen as a $\mathbb{R}^{2N}$-valued process, the fact that $h_{n,\varepsilon}^2 \leq 1$ and compactness of the support of the vector fields. Thus

		\begin{align*}
			\mathbb{E} \left[ W_2^{2p}(\pi_t^{n,\varepsilon}, \pi_s^{n,\varepsilon}) \right] &\leq \mathbb{E} \left[ W_{2p}^{2p}(\pi_t^{n,\varepsilon}, \pi_s^{n,\varepsilon}) \right] \\
			& \leq \mathbb{E} \left[ \int_{M^2} d^{2p}\left( \mathcal{Z}^{n,\varepsilon}_t[(x,y),\pi],\mathcal{Z}^{n,\varepsilon}_s[(x,y),\pi] \right) \d \pi(x,y)  \right] \\
			& \leq C_p(T) |t-s|^{p}.
		\end{align*}
		The conclusion follows from Theorem \ref{KCTC}.
	\end{proof}

	 A straightforward consequence of the tightness of $\left\{\mathbb{Q}^{n,\varepsilon}_{[(x,y),\pi]}\right\}_{(n,\varepsilon) \in \mathbb{N}\times \mathbb{R}_+^*}$ is that there exists two sequences $\{n_k\}_{k \in \mathbb{N}}$, $\{\varepsilon_l\}_{l \in \mathbb{N}}$ with $\lim_{l \to \infty} \varepsilon_l=0$ such that $\mathbb{Q}^{n_k,\varepsilon_l}_{[(x,y),\pi]} \to_{k \to \infty} \mathbb{Q}^{\varepsilon_l}_{[(x,y),\pi]}$ and $\mathbb{Q}^{\varepsilon_l}_{[(x,y),\pi]} \to_{l \to \infty} \mathbb{Q}_{[(x,y),\pi]}$. Moreover, we denote by $\mathbb{E}^{n_k,\varepsilon_l}_{[(x,y),\pi]}$ (resp. $\mathbb{E}_{[(x,y),\pi]}$) the expectation on $\Omega'$ with respect to $\mathbb{Q}^{n,\varepsilon}_{[(x,y),\pi]}$ (resp. $\mathbb{Q}_{[(x,y),\pi]}$).
	 Note that, since the map

	$$
	[(x,y),\pi] \in M^2 \times \P(M^2) \mapsto \mathbb{E} \left[ d^2(\mathcal{Z}_t^{1,n,\varepsilon}[(x,y),\pi],\mathcal{Z}_t^{2,n,\varepsilon}[(x,y),\pi]) \right]
	$$
	is measurable, the limit map $[(x,y),\pi] \in M^2 \times \P(M^2) \mapsto \mathbb{E}_{[(x,y),\pi]}\left[ d^2(X_t,\tilde{X}_t) \right]$ is also measurable. If there exists $K_1 \in \mathbb{R}, K_2 \leq 0$ such that

	\begin{equation} \tag{$\mathrm{IWC}_{M}'$($K_1,K_2$)} \label{Bakry emery}
		\langle \mathrm{Ric}(U), U \rangle - 2 \langle \nabla_{U} \nabla \phi_0(x,\mu) + \nabla \bar{\mathcal{L}}_{V} \phi_0(x,\mu), U \rangle \geq K_1 |U|_x^2 + K_2|V|_\mu^2,
	\end{equation}
	we say that $\mathrm{IWC}'_M(K_1,K_2)$ holds.

	\begin{rem}
		As in the previous subsection, if we assume that $K_2 \geq 0$, by taking $U= 0$ and $0 \neq V \in T_\mu \P$, we obtain

		$$
		0 \geq K_2 |V|^2_\mu.
		$$
		Thus, $K_2=0$. In particular, for $0 \neq U \in T_xM$, we obtain 

		$$
		\langle \mathrm{Ric}(U), U \rangle - 2 \langle \nabla_{U} \nabla \phi_0(x,\mu), U \rangle - K_1 |U|^2_x \geq 2 \langle \bar{\mathcal{L}}_{V} \phi_0(x,\mu), U \rangle.
		$$
		Since the right-hand side term is linear and bounded by a constant independent of $V$, $\bar{\mathcal{L}}_{V} \phi_0(x,\mu) =0$ for all $V$, thus, $\phi_0$ is constant in $\mu$. As in the previous subsection, in the case $K_2=0$, one should apply the finite-dimensional methods developped in Subsection \ref{S 1.3}.
	\end{rem}

	\begin{thm}
		Assume that $\mathrm{IWC}'_M(K_1,K_2)$ holds for $K_2 <0$. Then, 

		$$
		\mathbb{E}_\pi \left[W_2(\sigma^1_t, \sigma^2_t) \right] \leq  e^{-(K_1+K_2)t} W_2(\pi^1,\pi^2)
		$$
		and 

		$$
		\mathbb{E}_{[(x,y),\pi]} \left[d(X_t,\tilde{X}_t) \right] \leq e^{-K_1t} d^2(x,y) + e^{-(K_1+K_2)t}W_2^2(\pi^1,\pi^2).
		$$
	\end{thm}
	
	\begin{proof}
		First, note that by definition it holds:

	$$
	\mathbb{E}^{n_k,\varepsilon_l}_{[(x,y),\pi]}\left[W_2^2(\sigma^1_t, \sigma^2_t) \right] = \mathbb{E} \left[ W_2^2(\pi_t^{1,n_k,\varepsilon_l}, \pi_t^{2,n_k,\varepsilon_l}) \right] \leq \int_{M^2} \mathbb{E}^{n_k,\varepsilon_l}_{[(x,y),\pi]} \left[ d^2(X_t,\tilde{X}_t) \right] \d \Pi(x,y),
	$$
	for all $\Pi \in \P(M^2)$. Since $\mathbb{Q}^{n_k,\varepsilon_l}_{[(x,y),\pi]} \to_{k,l \to \infty} \mathbb{Q}_{[(x,y),\pi]}$ and since the considered maps are continuous, letting $k,l \to \infty$ yields

	\begin{equation} \label{Ineg Was lim}
		\mathbb{E}_{[(x,y),\pi]}\left[W_2^2(\sigma^1_t, \sigma^2_t) \right] \leq \int_{M^2} \mathbb{E}_{[(x,y),\pi]} \left[ d^2(X_t,\tilde{X}_t) \right] \d \Pi(x,y).
	\end{equation}
	Denote by $I_{\nabla \phi_0}$ the following map

		$$
	I_{\nabla \phi_0}((x,\mu),(y,\nu)) = I(x,y) + \nabla \phi_0(\cdot,\mu) d(\cdot,y)(x) + \nabla \phi_0(\cdot,\nu)d(x,\cdot)(y).
	$$
	Then, by adapting the arguments of Wang in \cite[Subsection 2.1]{Wang2005} (just replace the mirror reflection operator with the parallel transport) or in \cite[Section 3]{ArnaudonThalmaierWang2006}, we obtain that the process $d(X_t,\tilde{X}_t)$ has $\mathbb{Q}_{[(x,y),\pi]}$-almost surely finite variations and by applying $\mathrm{IWC}'_M(K_1,K_2)$, it follows that:

	\begin{align*}
		\d d(X_t, \tilde{X}_t) &\leq  I_{\nabla \phi_0}\left((X_t,\sigma^1_t), (\tilde{X}_t,\sigma^2_t) \right) \d t  \\
		& \leq   -K_1 \left(\int_0^1 |\partial_a \gamma_a|^2 \d a\right)^{1/2}  -K_2  \left(\int_0^1 |\partial_a \mu_a|^2 \d a \right)^{1/2} \d t  \\
		& \leq  -K_1 d(X_t, \tilde{X}_t) -K_2 W_2(\sigma^1_t, \sigma^2_t) \d t,
	\end{align*}
	where $(\gamma_a)_{0 \leq a \leq 1}$ is a geodesic from $X_t$ to $\tilde{X}_t$ and $(\mu_a)_{0 \leq a \leq 1}$ is a geodesic from $\sigma^1_t$ to $\sigma^2_t$. Let $\Pi$ be an optimal coupling of $\pi^1,\pi^2$, then, by considering the squared distance, taking the expectation and integrating with respect to $\Pi$, it follows from Cauchy-Schwarz inequality and \eqref{Ineg Was lim} that:

	\begin{align*}
		& \quad \d \int_{M^2} \mathbb{E}_{[(x,y),\pi]} \left[d^2(X_t,\tilde{X}_t) \right] \d \Pi(x,y) \\
		&\leq \int_{M^2} -2K_1 \mathbb{E}_{[(x,y),\pi]} \left[ d^2(X_t,\tilde{X}_t) \right] -2K_2 \mathbb{E}_{[(x,y),\pi]} \left[d(X_t,\tilde{X}_t) W_2(\sigma^1_t,\sigma^2_t) \right] \d \Pi(x,y) \d t \\
		&\leq \int_{M^2} -2K_1 \mathbb{E}_{[(x,y),\pi]} \left[ d^2(X_t,\tilde{X}_t) \right] -2K_2 \mathbb{E}_{[(x,y),\pi]} \left[d^2(X_t,\tilde{X}_t) \right]^{1/2}  \mathbb{E}_{[(x,y),\pi]}\left[W_2^2(\sigma^1_t,\sigma^2_t) \right]^{1/2} \d \Pi(x,y) \d t \\
		&\leq -2(K_1+K_2) \int_{M^2} \mathbb{E}_{[(x,y),\pi]} \left[ d^2(X_t, \tilde{X}_t) \right] \d \Pi(x,y) \d t.
	\end{align*}
	By Grönwall's inequality, it holds:

	$$
	\mathbb{E}_\pi \left[W_2^2(\sigma^1_t,\sigma^2_t) \right] \leq \int_{M^2} \mathbb{E}_{[(x,y),\pi]} \left[ d^2(X_t,\tilde{X}_t) \right] \d \Pi(x,y) \leq e^{-2(K_1+K_2)t} W_2^2(\pi^1,\pi^2).
	$$
	The remainder of the proof is as in the proof of Theorem \ref{T Ergo MP}.
	\end{proof}

	\section{Examples} \label{S 5}

	\subsection{Lift of equations on $M$} \label{S 5.1}

 Consider the following SDE

\begin{equation} \label{E Lift EDS classique sur M}
	\dw^{\grad} \mu_t = \SN V_{\phi_i}(\mu_t) \d W^i_t + V_{\phi_0}(\mu_t) \d t, \quad \mu_0 = \mu \in \P_\infty,
\end{equation}
where $\phi_i \in \C^\infty(M)$ for all $i=0,1, \dots, N$. The unique ID-solution is given by $(X_t)_* \mu$ where $(X_t)_{t \geq 0}$ is solution to:

\begin{equation} \label{EDS classique M}
\d^\nabla X_t(x) = \SN \nabla \phi_i(X_t(x)) \d W^i_t + \nabla \phi_0(X_t(x))\d t, \quad X_0(x) = x.
\end{equation}

\begin{prop} \label{P implication HCK IWCD}
	Assume that the $\mathrm{HC}(K)$ condition is satisfied for \eqref{EDS classique M}. Then, the $\mathrm{IWC}(K)$ condition is satisfied for \eqref{E Lift EDS classique sur M}. Assume moreover that $K > 0$ and let $\bar{\mu}$ be the unique invariant measure associated to $(X_t)_{t \geq 0}$ obtained by applying Corollary \ref{C Mesure invariante dim finie}. Then, the unique invariant measure $\Lambda$ associated to $(\mu_t)_{t \geq 0}$ is $\Lambda = \int_M \delta_{\delta_x} \d \bar{\mu}(x)$.
\end{prop}

\begin{proof}
	Let $\mu \in \P_\infty$, then,
	
	\begin{align}
	&\SN \langle \bar{R}(U, \bar{Z}_i) \bar{Z}_i, U \rangle  - \langle \grad_{U} \bar{Z}_i, \grad_{U} \bar{Z}_i \rangle -2  \langle \grad_{U}\bar{Z}_0 ,U \rangle \nonumber \\
	& \geq  \int_M  \SN \left(\langle R(U(x), \nabla \phi_i(x)) \nabla \phi_i(x), U(x) \rangle_x -  | \nabla_U \nabla  \phi_i(x) |^2 \right) - 2 \he \phi_0(x)(U,U) \, \d \mu(x) \label{E SU 1}\\
	& \geq K \|U\|^2_\mu, \label{E SU 2}
	\end{align}
	where \eqref{E SU 1} comes from Corollary \ref{Trace tenseur Riemann sur P} and the fact that $\Pi_\mu$ is a contraction and \eqref{E SU 2} from the $\mathrm{HC}(K)$ condition. 
	
	Assume that $\mathrm{HC}(K)$ holds for $K> 0$ and let $\bar{\mu}$ be the associated unique invariant measure and $\Lambda = \int_M \delta_{\delta_x} \d \bar{\mu}(x)$. Let also $F \in \mathscr{B}(\P)$ and $F_1 \in \mathscr{B}(M)$ be $F \circ p^1$, where $p^1 : M \rightarrow \P^1$ is the projection introduced in Appendix \ref{App ito}. Then, since $\M_t(\delta_x) = \delta_{X_t(x)}$, the conclusion follows from:

	\begin{align*}
		\int_{\P} \se_t F(\mu) \d \Lambda(\mu) &=  \int_M\int_{\P}  \mathbb{E} \left[F(\M_t(\mu)) \right] \d \delta_{\delta_x}(\mu) \d \bar{\mu}(x) \\
		&= \int_M \mathbb{E} \left[ F(\delta_{X_t(x)})\right] \d \bar{\mu}(x) \\
		& = \int_M \mathbb{E} \left[ F_1(X_t(x)) \right] \d \bar{\mu}(x) \\
		& = \int_M F_1(X_t(x)) \d \bar{\mu}(x) = \int_{\P} \int_M F(\mu) \d \delta_{\delta_x}(\mu) \d \bar{\mu}(x).
	\end{align*}
\end{proof}

		Let $\eta \in \P_{\mathrm{p.a}}$. One can check that, in this case,
		
		$$
		\T_{\eta} \P = L^2(M,TM,\eta).
		$$
		Note that this is not the case in general since
		
		$$
		L^2(M,TM,\mu) = \T_\mu \P \oplus \T_\mu \P^{\perp}, \quad \forall \mu \in \P_\infty,
		$$
		with
		
		$$
		\T_\mu \P^\perp = \overline{\left\{ Y \in \Gamma(M) \, : \, \mathrm{div}_\mu(Y) = 0 \right\}}^{L^2(\mu)}.
		$$
		Let $x_0 \in M$, $A \in \Gamma(M)$ and $(\mu_n)_{n \geq 0} \in (\P_\infty)^{\mathbb{N}}$ be a sequence of smooth measures converging toward $\delta_{x_0}$. Intuitively, one is inclined to think that 
		
		\begin{equation} \label{E convergence partie verticale}
		\lim_{n \to \infty} \int_M |\Pi_{\mu_n}^\perp(A)(x) |^2 \d \mu_n(x) = 0 
		\end{equation}
		since $\T_{\mu_n} \P^\perp$ tends to "disappear". If that were the case, we would have
		
		\begin{align*}
		\lim_{n \rightarrow \infty} \langle \bar{R}(U,V_{\phi_i})V_{\phi_i}, U \rangle_{\mu_n} &= \lim_{n \to \infty}\int_M \langle R(U(x), \nabla \phi_i(x)) \nabla \phi_i(x), U(x) \rangle \d \mu_n(x) \\
		&+ \lim_{n \to \infty} \int_M | \Pi_{\mu_n}^\perp( \nabla \nabla \phi_i) \cdot U)(x) |^2 \d \mu_n(x) \\
		& = \langle R(U(x_0), \nabla \phi_i(x_0)) \nabla \phi_i(x_0), U(x_0) \rangle.
		\end{align*}
		and
		
		\begin{align*}
		\lim_{n \to \infty} |\grad_{U} V_{\phi_i}|_{\mu_n}^2 &= \lim_{n \to \infty} \left( \int_M |\nabla_{U}\nabla \phi_i(x)|^2 \d \mu_n(x) - \int_M |\Pi_{\mu_n}^\perp( \nabla_U \nabla \phi_i)(x)|^2 \d \mu_n(x) \right) \\
		& = |\nabla_U \nabla \phi_i(x_0)|^2.
		\end{align*}
		In particular, assuming that the $\mathrm{IWC}(K)$ condition holds, taking the limit in \eqref{E IWCD} yields: 
		
		$$
		\langle R(U(x_0), \nabla \phi_i(x_0)) \nabla \phi_i(x_0), U(x_0) \rangle - |\nabla_U \nabla \phi_i(x_0)|^2 - 2 \he \phi_i(U(x_0), U(x_0)) \geq K |U(x_0)|^2.
		$$
		However, since $\Pi_\mu^\perp(A)$ is difficult to determine for $A \in \Gamma(M)$, we cannot guarantee that \eqref{E convergence partie verticale} holds for any base manifold $M$. Nevertheless, the computations in Example \ref{Exemple tore} allows us to prove the following result:
		
		\begin{prop}
			Let $M= \mathbb{S}^1$. The $\mathrm{HC}(K)$ condition for
			
			$$
			\d X_t(x) = \SN \nabla \phi_i(X_t(x)) \d W^i_t + \nabla \phi_0(X_t(x))\d t, \quad X_0(x) = x
			$$
			is equivalent to the $\mathrm{IWC}(K)$ condition for $(\mu_t)_{t \geq 0}$ solution to
			
			$$
			\dw^{\grad} \mu_t = \SN V_{\phi_i}(\mu_t) \d W^i_t + V_{\phi_0}(\mu_t) \d t, \quad \mu_0 = \mu \in \P_\infty.
			$$
		\end{prop}
		
		\begin{proof}
		We showed in Example \ref{Exemple tore} that for any vector field $A \in \Gamma(\mathbb{S}^1)$,
		
		$$
		\int_{\mathbb{S}^1} \left|\Pi_{\mu_n}^\perp(A)(x)\right|^2 \d \mu_n(x) = \left( \int_{\mathbb{S}^1} A(x) \d \vol(x) \right)^2 \left( \frac{1}{\int_{\mathbb{S}^1} \frac{\d \vol(x)}{\rho_n(x)} } \right),
		$$
		with $\rho_n = \d \mu_n / \d \vol$. Then, let $B^c = \mathbb{S}^1 \setminus B(x_0,r)$ for $r>0$ sufficiently small. By Cauchy-Schwarz inequality,

		$$
		\vol(B^c)^2 = \left(\int_{B^c} \d \vol(x)\right)^2 \leq \left(\int_{B^c} \rho_n(x) \d \vol(x) \right) \left( \int_{B^c} \frac{1}{\rho_n(x)} \d \vol(x) \right).
		$$
		Since $B^c$ is closed and $\delta_{x_0}(B^c)=0$, Portmanteau Lemma ensures that $\mu_n(B^c) \to_{n \to \infty}0$. Then, since $\rho_n >0$,

		$$
		\int_{\mathbb{S}^1} \frac{\d \vol(x)}{\rho_n(x)} \geq \frac{\vol(B^c)^2}{\mu_n(B^c)} \underset{n \to \infty}{\rightarrow} \infty.
		$$
		By following the above reasoning we obtain that the $\mathrm{IWC}(K)$ condition implies the $\mathrm{HC}(K)$ condition. The converse comes from Proposition \ref{P implication HCK IWCD}.
		
		\end{proof}

		\subsection{The Gaussian space}

		This section concerns the Gaussian space, which is a subspace of $\P_2(\mathbb{R}^d)$. Even if $\mathbb{R}^d$ is not a closed manifold, the Euclidean case is very well understood, the reader is refered to \cite{wang2020} for results on Image-dependent SDEs and to \cite{BHATIA2019165} for a presentation of the Riemannian structure of the Gaussian space. It can be shown that the Gaussian space $\mathscr{N}^d$ is isometric to $\mathbb{R}^d \times \mathscr{N}_0^d$, where $\mathscr{N}_0^d$ denotes the space of centered Gaussian probability measures:
	
	\begin{equation}
	\mathscr{N}_0^d := \left \{ \mu = \mathscr{N}(0,\Sigma) \, : \, \Sigma \in \mathscr{S}^{++}_d \right \},
\end{equation}
seen as a subspace of $\P_\infty$. Note that the restriction of the $W_2$ distance to $\mathscr{N}^d$ yields for two gaussian measures $N(x,\Sigma_1), N(y,\Sigma_2) \in \mathscr{N}^d$

\begin{equation*}
	W_2^2(N(x,\Sigma_1), N(y,\Sigma_2)) = |x-y|^2 + \mathrm{tr} \Sigma_1 + \mathrm{tr} \Sigma_2 - 2 \mathrm{tr}(\Sigma_1^{1/2}\Sigma_2 \Sigma_1^{1/2})^{1/2}.
	\end{equation*}
	Because of the flatness of $\mathbb{R}^d$, all the non-trivial geometric informations of $\mathscr{N}^d$ are concentrated in $\mathscr{N}_0^d$. Therefore, we will focus only on $\mathscr{N}^d_0$. The tangent bundle is given fiberwise by

	$$
	T_\mu \mathscr{N}^d_0 := \mathscr{S}_d,
	$$
	and the Riemannian metric is given by the restriction of the one of $\P_\infty$ on $\mathscr{N}^d_0$

	$$
	\langle Z_1, Z_2 \rangle_{N(\Sigma)} = \int_M \langle Z_1x, Z_2x \rangle_x d N(\Sigma)(x) = \mathrm{Tr}(Z_1 \Sigma Z_2).
	$$

	\begin{figure}[t] 
     \centering
     \includegraphics[width=0.7\textwidth]{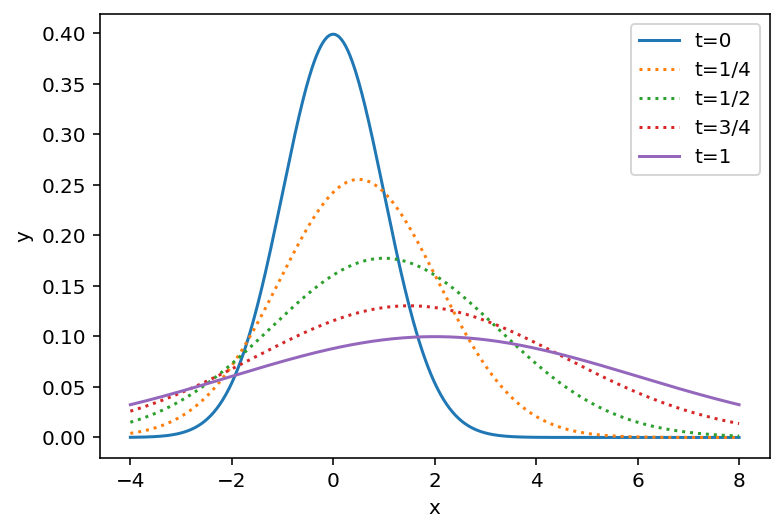}
     \caption{The displacement interpolation between $N(0,1)$ and $N(2,4)$ at times $t\in \left\{1,\frac{1}{4},\frac{1}{2},\frac{3}{4},1\right\}$.}\label{Transport optimal gaussiennes}
 \end{figure}

	Let us show that there exists a smooth Riemannian metric on $\mathscr{S}_d^{++}$ such that the latter is isometric to $\mathscr{N}_{0}^d$ with the above metric. We consider the usual trivialization for the tangent space of $\mathscr{S}_d^{++}$, namely,

	$$
	T\mathscr{S}_d^{++} = \mathscr{S}_d^{++} \times \mathscr{S}_d.
	$$

		 We denote by $\mathscr{S}^{++}_d$ the space of symmetric definite positive matrices on $\mathbb{R}^d$. For any $V \in \mathscr{S}^{++}_d$, we denote by $V^{1/2}$ the unique element of $\mathscr{S}^{++}_d$ such that $(V^{1/2})^2 = V$. 

		Let $\mathrm{GL}_d := \mathrm{GL}_d(\mathbb{R})$ be the space of invertible matrices seen as a smooth manifold and let 

		$$
	T_U \Inv := \left \{ AU \, : \, A \in M_d(\mathbb{R}) \right \},
	$$
	denote its tangent space in $U \in \Inv$. On the one hand, this space can be endowed with a Riemannian manifold structure with the following Riemannian metric:

	\begin{equation}\label{Metrique gld}
		\langle AU, BU \rangle_U = \mathrm{Tr}(AUU^\top B^\top).
	\end{equation}
	On the other hand, since an element $U \in \Inv(\mathbb{R})$ is as a smooth diffeomorphism of $\mathbb{R}^d$, we will also consider $\mathscr{G} = \Inv(\mathbb{R})$ seen as a finite-dimensional Riemannian submanifold of $\mathscr{D}$. With the Riemannian metric inherited from $\mathscr{D}$, one can easily check that the obvious bijection between $\mathscr{G}$ and $\Inv$ is an isometry. 
	Note that, the map
	
	\begin{align*}
		p_0 : \Inv & \longrightarrow \mathscr{S}^{++}_d \\
		U & \longmapsto U U^\top
	\end{align*}
	is a smooth map with differential given by

	\begin{align*}
		T_U p_0(AU) = A U U^\top + U U^\top A^\top 
	\end{align*}
	in these trivializations. One can easily show that the vertical space of $\Inv$ is
	
	\begin{equation}
	\V_M \Inv := \ker(Tp_0) =\left \{ Y (U^\top)^{-1}= Y (U U^\top)^{-1} U \, : \, Y \in \mathscr{A}_d(\mathbb{R}) \right \},
	\end{equation}
	and the horizontal space is
	
	\begin{equation}
	\H_M \Inv = \left \{ ZU \, : \, Z \in \mathscr{S}_d (\mathbb{R}) \right \}.
	\end{equation}

	\begin{lem} \label{L VBI}
		The vector bundle morphism $\Psi : T\mathscr{S}_d^{++} \rightarrow T\mathscr{S}_d^{++}$ defined fiberwise by:

		$$
		\Psi_{\Sigma}(Z) = \Sigma Z + Z \Sigma
		$$
		is a smooth vector bundle isomorphism over the identity.
	\end{lem}

	\begin{proof}
		Clearly, $\Psi \in \mathrm{Hom}^{\infty}_{\mathrm{id}}(T\mathscr{S}^{++}_d, T\mathscr{S}_d^{++})$. Moreover, for all $\Sigma \in \mathscr{S}^{++}_d$, $\Psi_\Sigma : \mathscr{S}_d \rightarrow \mathscr{S}_d$ is a linear injective map (one can check this fact by looking at the problem with $\Sigma = \mathrm{diag}(\lambda_1, \dots, \lambda_d)$ with $\lambda_i >0$ for all $i$). Thus, $\Psi_\Sigma$ is a bijective linear map for all $\Sigma \in \mathscr{S}^{++}_d$, the conclusion follows from the inverse mapping theorem.
	\end{proof}
	From now on we will consider the trivialization of $T \mathscr{S}^{++}_d$ given by $\Psi^{-1}(T\mathscr{S}_d^{++})$ and we endow the tangent bundle with the Riemannian metric given by
	
	$$
	\langle Z_1, Z_2 \rangle_{\Sigma} = \mathrm{Tr}(Z_1 \Sigma Z_2).
	$$
	In this new trivialization, for $U \in \Inv$ and $ZU \in \H_{U}\Inv$, the following holds

	$$
	T_{U}p_0(ZU) = Z \in T_{p_0(M)}\mathscr{S}^{++}_d.
	$$
	In particular, for $Z_1U, Z_2U \in \H_{U} \Inv$,

	$$
	\langle T_Up_0(Z_1U), T_Up_0(Z_2U) \rangle_{p_0(U)} = \mathrm{Tr}(Z_1UU^\top Z_2) = \langle Z_1 U, Z_2 U \rangle_U.
	$$
	We obtained the following:

	\begin{prop}
		The projection $p_0 : \Inv \rightarrow \mathscr{S}^{++}_d$ is a Riemannian submersion.
	\end{prop}
	Note that, as a straightforward corollary of Lemma \ref{L VBI}, the horizontal projection of a vector $AU \in T_U \Inv$ is given by:

	$$
	P_{\H_U}(AU) = \Psi_{p_0(U)}^{-1}(A p_0(U) + p_0(U)A^\top).
	$$

	\begin{rem}
		As already mentioned in Example \ref{Exemple tore}, the normal tensor and more generally the vertical projection of vector fields on $\mathscr{D}$ is hard to compute. Proposition \ref{L VBI} allows us to compute explicitly the horizontal part and, consequently, the vertical part of any vector of $\Inv$.
	\end{rem}

	Finally, if $\tilde{p}_0 : \mathscr{G} \rightarrow \mathscr{N}^d_0$ is the map defined by

	$$
	\tilde{p}_0(\varphi) = (\varphi)_* N(0,\mathrm{Id}),
	$$
	then, the diagram

	 \begin{center}
\begin{tikzpicture}[>=stealth]
  \node (E1) at (0,1) {$\mathscr{G}$};
  \node (E2) at (4,1) {$\Inv$};
  \node (E3) at (0,0)   {$\mathscr{S}_d^{++}$};
  \node (E4) at (4,0)   {$\mathscr{N}_{0}^d$};

  \draw[<->, >=latex] (E1) -- (E2); 
  \draw[->, >=latex] (E1) -- (E3) node[midway,left]{$p_0$};
  \draw[<->, >=latex] (E3) -- (E4); 
  \draw[->, >=latex] (E2) -- (E4) node[midway,right]{$\tilde{p}_0$};
\end{tikzpicture}
\end{center}
is commutative with the obvious isometries between $\mathscr{N}^d_0$ and $\mathscr{S}_d^{++}$ and between $\mathscr{G}$ and $\Inv$. In particular,
	
	\begin{equation} \label{Distance Wasserstein cas gaussien}
	W_2^2(\mu_1,\mu_2) = d^{WB}(\Sigma_1,\Sigma_2)^2,
	\end{equation}
	where
	
	\begin{equation*}
	d^{WB}(\Sigma_1,\Sigma_2) = \left( \mathrm{tr} \Sigma_1 + \mathrm{tr} \Sigma_2 - 2 \mathrm{tr}(\Sigma_1^{1/2}\Sigma_2 \Sigma_1^{1/2})^{1/2} \right)^{1/2}
	\end{equation*}
	is the so-called Wasserstein-Bures distance. We can also derive the following result which has been proven in \cite{oneill1966fundamental} and which is also a straightforward consequence of Corollary \ref{Trace tenseur Riemann sur P}.
	
	\begin{prop}
		The following formula for the Riemann curvature tensor on $\mathscr{N}^d_0$ holds:
		
		$$
		\left \langle R(Z_1,Z_2)Z_2, Z_1 \right \rangle_\Sigma = \frac{3}{4} \left \| \left[ Z_1\Sigma^{1/2},Z_2 \Sigma^{1/2} \right]^\V \right \|^2_{\Sigma^{1/2}},
		$$
		where $(\cdot)^{\mathscr{V}}$ denote the vertical projection. In particular, $\mathscr{N}^d_0$ has non-negative sectional curvature.
	\end{prop}

	\begin{prop}
		Let $\bar{Z}_0, \bar{Z}_1, \dots, \bar{Z}_N \in \Gamma^1(\mathscr{N}_0^d)$. Then, the SDE
		
		\begin{equation} \label{E EDS N}
		\dw^{\grad} \M_t(\mu) = \SN \bar{Z}_i(\M_t(\mu)) \d W^i_t + \bar{Z}_0(\M_t(\mu)) \d t, \quad \M_0(\mu) = \mu \in \mathcal{N}_0^d,
		\end{equation}
	admits a unique ID solution given by $\M_t(\mu) = (\varphi_t)_* \mu$ where $(\varphi_t)_{t \geq 0}$ is solution to the following linear Image-dependent SDE:
	
	\begin{equation} \label{E ID N}
	\d\varphi_t(x) = \SN Z_i(\M_t(\mu)) \cdot \varphi_t(x) \d W^i_t + \SN Z_0(\M_t(\mu))\cdot \varphi_t(x) \d t, \quad \varphi_0 = \mathrm{id}, \, \M_t(\mu) = (\varphi_t)_* \mu.
	\end{equation}
	In particular, $(\varphi_t)_{t \geq 0}$ is a family of linear maps and, consequently, $\M_t(\mu)$ lies in $\mathcal{N}_0^d$.
	\end{prop}
	
	\begin{proof}
	This result is an application of \cite[Theorem 2.1]{wang2020}. 
	\end{proof}
	
	\noindent \textbf{Study of diffusions on $\mathscr{N}_0^d$.} Let $\bar{Z}_0, \bar{Z}_1, \dots, \bar{Z}_N \in \Gamma^1(\mathscr{N}_0^d)$ and $(\M_t(\mu))_{t \geq 0}$ be solution to

			$$
		\d^\nabla \M_t(\mu) = \SN \bar{Z}_i(\M_t(\mu)) \d W^i_t + \bar{Z}_0(\M_t(\mu)) \d t, \quad \M_0(\mu) \in \mathscr{N}_0^d.
		$$
	In order to study the stability of such diffusions we can explore at least two options: study this process as a $\mathcal{N}_0^d$-valued process and use the $\mathrm{IWC}(k)$ condition or look at it as an $\mathcal{S}_{d}^{++}(\mathbb{R})$-valued process with the above Riemannian structure and use the $\mathrm{HC}(K)$ condition. As one may expect, these conditions are equivalents in this context since the finite-dimensional Riemannian structure is exactly the one induced by the Otto's metric.
	
	\begin{prop}
		The $\mathrm{IWC}(K)$ condition on $\mathscr{N}_0^d$ seen as a submanifold of $\P_\infty$ is equivalent to the $\mathrm{HC}(K)$ condition on $\mathscr{S}_d^{++}$ with the Riemannian structure described above.
	\end{prop}
	
	One can also consider a Langevin diffusion on $\mathscr{N}^d_0$ by setting $(X_t)_{t \geq 0}$ a Brownian motion on $\mathcal{N}^d_0$ with drift, namely, a diffusion with generator
	
	$$
	L = \frac{1}{2} \left( \Delta - \langle \nabla F, \cdot \rangle\right),
	$$
	where $F$ is choosen to be a geodesically convex function on $\mathcal{N}_0^d$. For instance, one may take $F= F_f: \mathcal{S}^{++}_d \rightarrow \mathbb{R}$, the smooth function defined to be
	
	$$
	F_f(A) = \int_M f(x) \d\mathcal{N}(A)(x), \quad \forall A \in \mathcal{S}^{++}_d,
	$$
	with $f : \mathbb{R}^d \rightarrow \mathbb{R}$ a convex function. In particular, all results of Subsection \ref{S 1.3} apply.

	\appendix
	
	\section{Appendix} \label{Appendix A}

	\subsection{Covariant equations on the Euclidean space} \label{appendice A1}

	Let $M$ be a Riemannian submanifold of $\mathbb{R}^d$ with the usual flat metric. We denote by $\II : TM \times TM \rightarrow TM^\perp$ the associated second fundamental form and for $v \in T_xM$, we denote by $S_v : T_xM \rightarrow T_xM^\perp$ the shape operator defined by

	$$
	S_v(u) = \II(u,v).
	$$
	We recall that $\II$ is a symmetric tensor. Let $A \in \Gamma(TM)$, the adjoint $S_A^\top$ is given by

	\begin{align*}
		TM^\perp &\longrightarrow TM \\
		N &\longmapsto P_{TM}\left( \left(D A\right)^\top \cdot N \right).
	\end{align*}
	We extend this definition to the whole $TM \oplus TM^\perp$ by setting $S_A^\top(B)=0$ for all $B \in \Gamma(TM)$. Note that for $A,B \in \Gamma(TM)$, $N \in \Gamma(TM^\perp)$, the following holds:

	$$
	\langle B, S^\top_A(N) \rangle = \langle \II(A,B), N \rangle = \langle D B(A), N \rangle = - \langle B, D N(A) \rangle.
	$$
	Consequently, it holds:

	\begin{equation} \label{E 17}
		S_A^\top(N) = P_{TM}\left( \left(D A\right)^\top \cdot N \right) = - P_{TM}\left( D N (A) \right) .
	\end{equation}
	Let $A_0, A_1, \dots, A_N \in \Gamma(TM)$. One can easily check that if $\tilde{A}_0, \tilde{A}_1, \dots, \tilde{A}_N$ are smooth extensions with compact support of the driving vector fields, then the unique solution flow to

	\begin{equation} \label{E 15}
		\d \tilde{X}_t(x) = \SN \tilde{A}_i(\tilde{X}_t(x)) \d W^i_t + \tilde{A}_0(\tilde{X}_t(x)) \d t + \frac{1}{2} \SN \tilde{\II}(\tilde{A}_i(\tilde{X}_t), \tilde{A}_i(\tilde{X}_t)), \quad X_0(x) = x \in M
	\end{equation}
	lie in $M$ and is actually solution to

	\begin{equation} \label{E 16}
		\d^\nabla X_t(x) = \SN A_i(X_t(x)) \d W^i_t + A_0(X_t(x)) \d t, \quad X_0(x) =x \in M.
	\end{equation}
	For any $A \in \Gamma(TM)$, it is shown in \cite[Equation (2.21)]{kunita1981} that:

	$$
	D_t^{\nabla} A(X_t(x)) = \SN \nabla A(A_i(X_t(x))) \d W^i_t + \nabla A(A_0(X_t(x))) \d t + \frac{1}{2} \SN \nabla^2 A (A_i(X_t(x)),A_i(X_t(x))) \d t
	$$

	\begin{prop} \label{E I-S Euclid}
		Let $\tilde{L}^0, \tilde{L}^1, \dots, \tilde{L}^N : \mathbb{R}_+ \rightarrow TM \oplus TM^\perp$ such that for all $i=0,1,\dots,N$:
	
	$$
	\tilde{L}_t^i = L_t^i \oplus (L_t^i)^\perp, \quad t \in \mathbb{R}_+.
	$$
	Assume that $(\tilde{V}_t)_{t \geq 0}$ is a $TM$-valued stochastic process solution to

	\begin{equation}\label{E 19}
		\d \tilde{V}_t = \SN \tilde{L}_t^i \d W^i_t + \tilde{L}_t^0 - \frac{1}{2} \SN S_{A_i}^\top((L_t^i)^\perp) \d t. \quad \tilde{V}_0 = V \in T_xM.
	\end{equation}
	Then, $(\tilde{V}_t)_{t \geq 0}$ is the unique solution to

	$$
	D^\nabla_t V_t = \SN L_t^i \d W^i_t + L_t^0(V_t) \d t, \quad V_0 =V \in T_xM.
	$$
	\end{prop}

	\begin{proof}
		Let $(X_t)_{t \geq 0}$ be solution to \eqref{E 15}, $A \in \Gamma(TM)$ and $\tilde{A} : \mathbb{R}^d \rightarrow \mathbb{R}^d$ a smooth extension. Then, by applying Itô's formula we obtain that $A(X_t)$ is solution to
		
		\begin{align*}
			\d A(X_t) = \SN D A(A_i(X_t)) \d W^i_t + D A(A_0(X_t)) + \frac{1}{2} \SN DA (\II(A_i,A_i)) + D^2 A (A_i,A_i) \d t 
		\end{align*}
		Since $D^2 A (A_i,A_i) = D \left( D A(A_i) \right)(A_i)- DA (D A_i (A_i))$, we obtain that

		\begin{align*}
			\d A(X_t) &= \SN D A(A_i(X_t)) \d W^i_t + D A(A_0(X_t)) + \frac{1}{2} \SN + D(DA(A_i))(A_i) - DA (D A_i (A_i)) + DA( \II(A_i,A_i)) \d t \\ 
			& = \SN D A(A_i(X_t)) \d W^i_t + D A(A_0(X_t)) + \frac{1}{2} \SN + D(DA(A_i))(A_i) - DA (\nabla_{A_i} A_i) \d t
		\end{align*}
		Thus,

		$$
		D(DA(A_i))(A_i) = D(\nabla_{A_i}A)(A_i) + D(\II(A_i,A))(A_i) = D(\nabla_{A_i}A)(A_i) + D (S_{A_i}(A))(A_i).
		$$
		Finally, $A(X_t)$ is solution to

		\begin{align*}
			\d A(X_t) &=\SN D A(A_i(X_t)) \d W^i_t + D A(A_0(X_t)) + \frac{1}{2} \SN D(\nabla^M_{A_i}A)(A_i) +  D(S_{A_i}(A))(A_i) - DA( \nabla^M_{A_i}A_i) \d t.
		\end{align*}
		By \eqref{E 17} if $(\tilde{V}_t)_t$ is solution to \eqref{E 19}, the following holds:

		\begin{align*}
			\d \langle \tilde{V}_t, A(X_t) \rangle_{X_t} &= \langle d \tilde{V}_t, A(X_t) \rangle + \langle \tilde{V}_t, d A(X_t) \rangle + \langle \d \tilde{V}_t, \d A(X_t) \rangle \\
			& = \SN \left\langle L_t^i, A(X_t) \right\rangle \d W^i_t + \left\langle L_t^0, A(X_t) \right\rangle \d t \\
			& + \SN \left\langle \tilde{V}_t, \nabla_{A_i}A(X_t) \right\rangle \d W^i_t +\left\langle \tilde{V}_t, \nabla_{A_0}A(X_t) \right\rangle \d t \\
			& + \SN \left\langle \tilde{L}_t^i, DA(A_i(X_t)) \right\rangle - \frac{1}{2} \left(\left\langle S_{A_i}^\top ((L_t^i)^\perp), A(X_t) \right\rangle + \left\langle \tilde{V}_t, S_{A_i}^\top(S_{A_i}(A)) \right\rangle \right) \d t \\
			& + \SN \left\langle \tilde{V}_t, \nabla^2 A(A_i,A_i) \right\rangle \d t.
		\end{align*}
		On the one hand:

		$$
		\left\langle \left( \d \tilde{V}_t \right)^\perp, \left( \d A(X_t(x)) \right)^\perp \right\rangle =  \SN \left\langle (L_t^i)^\perp, \II(A_i,A) \right\rangle \d t,
		$$
		and, on the other hand,

		$$
		\left\langle \left( \d \tilde{V}_t \right)^\perp, \left( \d A(X_t(x)) \right)^\perp \right\rangle = \left\langle S_{\d X_t(x)} (\tilde{V}_t), S_{\d X_t(x)} (A(X_t(x))) \right\rangle = \SN \langle S_{A_i}(\tilde{V}_t), S_{A_i}(A(X_t(x))) \rangle \d t.
		$$
		Thus, the conclusion follows from the fact that

		\begin{align*}
			&  \quad \SN \left\langle \tilde{L}_t^i, DA(A_i(X_t)) \right\rangle - \frac{1}{2} \left(\left\langle S_{A_i}^\top ((L_t^i)^\perp), A(X_t) \right\rangle + \left\langle \tilde{V}_t, S_{A_i}^\top(S_{A_i}(A)) \right\rangle \right) \\
			& = \SN \left\langle \tilde{L}_t^i, DA(A_i(X_t)) \right\rangle - \left\langle (L_t^i)^\perp, \II(A_i,A) \right\rangle \\
			& = \SN \left\langle L_t^i, \nabla_{A_i}A \right\rangle.
		\end{align*}
	\end{proof}

	\begin{coro} \label{Coro Pro Tan}
		Let $V \in T_xM$, $\tilde{V}_t = D_x \tilde{X}_t(V)$. Then, $\tilde{V}_t = T_x X_t(V)$ and is solution to

	$$
	D^\nabla_t V_t = \SN \nabla A_i(V_t) \d W^i_t + \nabla A_0(V_t) - \frac{1}{2} \SN R(V_t,A_i)A_i \d t, \quad V_0 =V \in T_xM.
	$$

	\end{coro}

	\begin{proof}
	The fact that $\tilde{V}_t = T_x X_t(V)$ is obvious by definition. Then, note that
	
	$$
	D \II(A,B)(U) \overset{(n)}{=} -S^\top_{U}(\II(A,B)),
	$$
	where $\overset{(\mathrm{n})}{=}$ denotes equality up to normal terms. Thus,

	$$
	\d \tilde{V}_t = \SN D A_i(V_t) \d W^i_t + D A_0(V_t) + L'_0(V_t) \d t,
	$$
	with

	$$
	L_0'(U) \overset{(\mathrm{n})}{=} -\frac{1}{2}\SN S^\top_{U}(S_{A_i}(A_i)), \quad U \in T_xM.
	$$
	A straightforward application of the Gauss equation (see \cite[Theorem 8.4]{Lee2018} for instance) yields:

	$$
	L_0'(U) \overset{(\mathrm{n})}{=} -\frac{1}{2}\SN \left(R^M(U,A_i)A_i + S^\top_{A_i}(S_{A_i}(U))\right).
	$$
	The conclusion follows from \ref{E I-S Euclid}. 
	\end{proof}

	\subsection{Proof of Itô formula} \label{App ito}
	Let $M$ be a closed Riemannian manifold and endow $M^n$ with the Cartesian product Riemannian metric. Let also

	$$
	\P^n = \left\{ \mu = \frac{1}{n} \sum_{i=1}^n  \delta_{x_i}, \, \mathbb{X} = \{x_i\}_{i=1}^n \in M^n \right\}
	$$
	endowed with the $W_2$ distance. The map 

	\begin{align*}
		p^n : M^n &\longrightarrow \P^n \\
		\mathbb{X} = \{x_i\}_{i=1}^n &\longmapsto \mu^n_{\mathbb{X}} := \frac{1}{n} \sum_{i=1}^n \delta_{x_i}
	\end{align*}
	is clearly an onto continuous map. For any function $F :  \P^n\to \mathbb{R}$, the $M^n$-lift $F_n : M^n \to \mathbb{R}$ is defined to be $F_n = F \circ p^n$. Note that, in the case $n=1$, $p^1$ is a homeomorphism.
	
	\begin{lem} \label{Lemme A4}
		If  $F \in \C^2(\P)$, then $F_n \in \C^2(M^n)$ and, for $\mathbb{U} = \left\{U_i\right\}_{i=1}^n$ and $\mathbb{V} = \left\{V_i\right\}_{i=1}^n$ in $T_xM^n$,

		$$
		\nabla^n F_n(\mathbb{X}) = \frac{1}{n} \left( \grad F(\mu^n_{\mathbb{X}},x_1), \dots, \grad F(\mu^n_{\mathbb{X}},x_n) \right)
		$$
		and

		$$
		\mathrm{H}^n F_n(\mathbb{X})[\mathbb{U},\mathbb{V}] = \frac{1}{n} \sum_{i=1}^n \left\langle \nabla_{V_i} \grad F(\mu^n_{\mathbb{X}}, x_i), U_i \right\rangle + \frac{1}{n^2} \sum_{i,j=1}^n \left\langle \grad^2 F(\mu^n_{\mathbb{X}},x_i,x_j), U_i \otimes V_j \right\rangle.
		$$
	\end{lem}

	\begin{proof}
		 Let $s \mapsto \mathbb{X}_s =\{\exp_{x_i(0)}(U_i)\}_{i=1}^n$. By definition of the $L$-differentiability,

		$$
		\left. \frac{\d}{\d s} F_n(\mathbb{X}_s) \right|_{s=0}= \left. \frac{\d}{\d s} F(\mu^n_{\mathbb{X}_s}) \right|_{s=0} = \left\langle \grad F(\mu^n_{\mathbb{X}_0}), \partial_s|_{s=0}\mu^n_{\mathbb{X}_s} \right\rangle_{\mu^n_{\mathbb{X}_0}}  = \frac{1}{n} \sum_{i=1}^n \left\langle \grad F (\mu^n_{\mathbb{X}_0}, x_i), U_i \right\rangle.
		$$
		Further derivatives can be computed in the exact same way.
	\end{proof}

	\begin{coro}
		Let $\mu^n = \frac{1}{n} \sum_{i=1}^n \delta_{x_i}$, $\mu_t^n = \M_t(\mu^n)$ and $\mathbb{X}_t^n= \left\{\varphi_t(x_i,\mu^n)\right\}_{i=1}^n $. Then, for $F \in \C^2(\P)$,

		$$
		F(\mu_t^n) = F(\mu^n) + \SN \int_0^t \bar{\mathcal{L}}_{\bar{Z}_i} F(\mu_s^n) \d W^i_s + \int_0^t \underbrace{\left[\bar{\mathcal{L}}_{\bar{Z}_0} F(\mu_s^n) + \frac{1}{2} \SN \bar{\mathrm{H}} F(\mu^n_s)(\bar{Z}_i, \bar{Z}_i) \right]}_{\bar{\mathrm{L}} F(\mu^n_s)} \d s.
		$$
	\end{coro}

	\begin{proof}
		Since each component of $\mathbb{X}_t^n$ is solution to an Image-dependent SDE on $M$ and since $F_n \in \C^2(M^n)$, one can apply Itô formula to compute $F_n(\mathbb{X}_t^n)$. To conclude, it is sufficient to identify the terms with Lemma \ref{Lemme A4}.
	\end{proof}

	\begin{lem} \label{Lemme continuité Lp}
		Let $(E,d)$ be a metric space, $F: E \rightarrow \mathbb{R}$ be a bounded and uniformly continous map. Then, $F : L^p(\Omega,E;\mathbb{P}) \rightarrow L^p(\Omega,\mathbb{R};\mathbb{P})$ is continuous for all $p \geq 1$.
	\end{lem}

	\begin{proof}
		Let $p \geq 1$, $\varepsilon>0$ and $\delta >0$ such that if $d^{p}(x,y)< \delta$, $|F(x)-F(y)|^p< \varepsilon$. Let also $X,Y \in L^p(\Omega,E;\mathbb{P})$, $A = \{\omega \in \Omega, \, d^p(X(\omega),Y(\omega) < \delta) \}$ and assume that 

		$$
		\mathbb{E}\left[d^p(X,Y)\right] < \frac{\delta \varepsilon}{C_p |F|_\infty^p },
		$$
		where $C_p$ is a constant such that $(a+b)^p \leq C_p(a^p + b^p)$. 
		Then, it follows by Markov inequality:

		\begin{align*}
			\mathbb{E} \left[ |F(X)-F(Y)|^p \right] &< \varepsilon + \mathbb{E} \left[ |F(X)-F(Y)|^p \bm 1_{\Omega \setminus A} \right] \\
			& \leq \varepsilon + C_p |F|_\infty^p \mathbb{P}(d^p(X,Y) \geq \delta) \\
			& \leq \varepsilon + \frac{C_p |F|_\infty^p}{\delta} \mathbb{E} \left[ d^p(X,Y) \right] < 2\varepsilon.
		\end{align*}
		The proof is completed.
	\end{proof}

	\textit{Proof of Proposition \ref{Ito formula}.} Since $\P^n$ is dense in $\P$ (see \cite[Proposition 7.1.5]{ambrosioGradientFlowsMetric2008} or \cite[Theorem 6.18]{villani2009optimal} for instance), let $\{\mu^n\}_{n \in \mathbb{N}}$ be a sequence in $\mathbb{\P}^{\mathbb{N}}$ such that for all $n \in \mathbb{N}$, $\mu^n \in \P^n$ and 

	$$
	\lim_{n \to \infty} W_2^2(\mu^n,\mu) = 0.
	$$
	By \eqref{Estimee lip sup flot}, we obtain that 

	\begin{equation} \label{Continuite en tps}
		\lim_{n \to \infty} \mathbb{E} \left[ \sup_{t \in [0,T]} W_2^2(\mu_t^n, \mu_t) \right] = 0.
	\end{equation}
	By continuity of $F$, $F(\mu^n) \to F(\mu)$. Moreover, from \eqref{Continuite en tps}, continuity of $\nu \mapsto \mathcal{L}_{\bar{Z}_i}F(\nu)$, Lemma \ref{Lemme continuité Lp} with $(E,d)= (\P,W_2)$ and by taking the limit under the integral sign, it follows:

	$$
	\lim_{n \to \infty}  \int_0^T \mathbb{E} \left[ \left| \bar{\mathcal{L}}_{\bar{Z}_i}F (\mu^n_s) - \bar{\mathcal{L}}_{\bar{Z}_i}F (\mu_s) \right|^2 \right] \d s  = 0.
	$$
	The Burkholder-Davis-Gundy inequality then implies

	$$
	\lim_{n \to \infty} \mathbb{E} \left[ \sup_{t \in [0,T]} \left| \int_0^t \bar{\mathcal{L}}_{\bar{Z}_i}F (\mu^n_s) - \bar{\mathcal{L}}_{\bar{Z}_i}F (\mu_s) \d W^i_s \right|^2 \right] = 0.
	$$
	Also by \eqref{Continuite en tps}, continuity, Lemma \ref{Lemme continuité Lp} with $(E,d)= (\P,W_2)$ and compactness of $\P$, one can show that 

	$$
	\mathbb{E} \left[ \sup_{t \in [0,T]} \left| \int_0^t \bar{\mathrm{L}}F (\mu^n_s) - \bar{\mathrm{L}}F (\mu_s) \d s \right|^2 \right] \leq T \int_0^T \mathbb{E} \left[ |\bar{\mathrm{L}}F (\mu^n_s) - \bar{\mathrm{L}}F (\mu_s)|^2 \right] \d s \underset{n \to \infty}{\rightarrow} 0.
	$$
	Finally, the conclusion follows from \eqref{Continuite en tps} and

	$$
	\mathbb{E} \left[ \sup_{t \in [0,T]}|F(\mu^n_t)-F(\mu_t)|^2 \right] \leq |\grad F|^2_\infty \mathbb{E}\left[\sup_{t \in [0,T]}W_2^2(\mu^n_t, \mu_t)\right].
	$$

	\subsection{Proof of Lemma \ref{L estim lip}} \label{Appendix 3}

	\begin{lem}[= Lemma \ref{L estim lip}]
			There exists a constant $C>0$ such that for any $ (x,v,\zeta),(y,u,\eta) \in \mathbb{R}^{2d} \times \P_c(M \times \mathbb{R}^d)$ with $\pi_* \zeta = \mu, \, \pi_* \eta = \nu$, the following holds:

\begin{align}
	&\sum_{i=0}^N |B_i(x,v,\zeta)|^2 \leq C(1 + |v|^2 + |\zeta|^2) ;\\
	& \begin{array}{l}
		\langle B_0^1(x,u,\mu) - B_0^1(y,v,\nu), u-v \rangle + \SN |B_i^1(x,u,\mu)- B_i^1(y,v,\nu)|^2 \\
	 \leq K \left[ \min(|u|^2,|v|^2) \cdot ( |x-y|^2 + W_2^2(\mu,\nu))  + |u-v|^2 \right] ;
	\end{array} \\
	&\langle B_0^2(x,\zeta) - B_0^2(y,\eta), u-v \rangle + \SN |B_i^2(x,\zeta)- B_i^2(y,\eta)|^2 \leq K \min(|\zeta|^2,|\eta|^2) \left(|x-y|^2  + W_2^2(\zeta,\eta)\right)+ K |u-v|^2.
\end{align}

		\end{lem}
		\begin{proof}
			The first point follows from

			\begin{align*}
				|B_i(x,v,\zeta)|^2 &\leq 2 \left( |T_x A_i(\cdot,\mu)(v)|^2 + \left| \int_{M \times \mathbb{R}^d} \left\langle \grad A_i(x,\mu)(y),u \right\rangle \d \zeta(y,u) \right|^2 \right) \\
				& \leq 2  \sup_{(x,\mu) \in M \times \P(M)} |T_x A_i(\cdot,\mu)|^2 \cdot |v|^2  \\
				& + 2 \left( \int_{M} |\grad A_i(x,\mu)(y)|^2 \d \mu_1(y) \right) \cdot \left(\int_{M \times \mathbb{R}^d} |u|^2 \d \zeta(y,u) \right) \\
				& \leq 2\left( \sup_{(x,\mu) \in M \times \P(M)} |T_x A_i(\cdot,\mu)|^2 \cdot |v|^2 + \sup_{(x,y,\mu)\in M^2 \times \P(M)} |\grad A_i(x,\mu)(y)|^2 \cdot |\zeta|^2 \right).
			\end{align*}
			The second point follows from the fact that for all $i=0,1,\dots,N$:
			
			\begin{align*}
				& \quad |B_i^1(x,u,\mu)-B_i^1(y,v,\nu)|^2 \\
				& \leq 3 \left( |B_i^1(x,u,\mu)-B_i^1(y,u,\mu)|^2 + |B_i^1(y,u,\mu)-B_i^1(y,v,\mu)|^2 + |B_i^1(y,v,\mu)-B_i^1(y,v,\nu)|^2 \right) \\
				& \leq \sup_{(x,\mu) \in M \times \P}|\nabla^2 A_i(x,\mu)|^2 |x-y|^2 \cdot |u|^2 + \sup_{(x,\mu) \in M \times \P} |\nabla A_i(x,\mu)| |u-v|^2 \\
				& + \sup_{(x,\mu) \in M \times \P}|\grad \nabla A_i(x,\mu)| W_2^2(\mu,\nu) \cdot |u|^2.
			\end{align*}
			Since 

			$$
			\langle B_0^1(x,u,\mu) - B_0^1(y,v,\nu), u-v \rangle \leq \frac{1}{2} \left( |B_0^1(x,u,\mu) - B_0^1(y,v,\nu)|^2 + |u-v|^2 \right),
			$$
			and since the roles of $u$ and $v$ can be exchanged, the result follows.
			
			Let us prove the third point.
			Let $\pi$ be an optimal coupling between $\zeta$ and $\eta$. The Cauchy-Schwarz inequality yields together with the compactness of $M$ and the regularity properties of $A_i$ for all $i=0,1,\dots,N$ yields:

			\begin{align*}
				& \quad |B_i^2(x,\zeta) - B_i^2(y,\eta)|^2 \\
				& \lesssim |B_i^2(x,\zeta)- B_i^2(y,\zeta)|^2 + |B_i^2(y,\zeta) - B_i^2(y,\eta)|^2 \\
				& \lesssim \int_M |x-y|^2 \cdot |u|^2 \d \zeta(z,u) + \left( \int_{M \times \mathbb{R}^d} \langle \grad A_i(y,\mu)(z), u \rangle \d \zeta(z,u) - \langle \grad A_i(y,\nu)(z), u \rangle \d \eta(z,u) \right)^2 \\
				& \leq |x-y|^2 \cdot |\zeta|^2 + \left( \int_{M \times \mathbb{R}^d} \langle \grad A_i(y,\mu)(z_1), u_1 \rangle - \langle \grad A_i(y,\nu)(z_2), u_2 \rangle \d \pi[(z_1,z_2),(u_1,u_2)] \right)^2.
			\end{align*}
			Now, by the Lipschitz property of $A_i$, it follows,

			\begin{align*}
				& \quad  \int_{M \times \mathbb{R}^d} \langle \grad A_i(y,\mu)(z_1), u_1 \rangle - \langle \grad A_i(y,\nu)(z_2), u_2 \rangle \d \pi[(z_1,z_2),(u_1,u_2)] \\
				& \leq  \int_{M \times \mathbb{R}^d} |\langle \grad A_i(y,\mu)(z_1)- \grad A_i(y,\nu)(z_1), u_1 \rangle| \d \pi[(z_1,z_2),(u_1,u_2)] \\
				& + \int_{M \times \mathbb{R}^d} |\langle \grad A_i(y,\nu)(z_1)- \grad A_i(y,\nu)(z_2), u_1 \rangle| \d \pi[(z_1,z_2),(u_1,u_2)] \\
				& + \int_{M \times \mathbb{R}^d} \langle \grad A_i(y,\nu)(z_2), u_1-u_2 \rangle  \d \pi[(z_1,z_2),(u_1,u_2)] \\
				& \leq C_1 \int_{M \times \mathbb{R}^d} W_2(\mu,\nu) |u_1| \d \pi[(z_1,z_2),(u_1,u_2)] +  C_2\int_{M \times \mathbb{R}^d} |u_1| \cdot |z_1-z_2| \d \pi[(z_1,z_2),(u_1,u_2)] \\
				& + C_3  \int_{M \times \mathbb{R}^d} |u_1-u_2| \d \pi[(z_1,z_2),(u_1,u_2)] \\
				& \leq C_1 |\zeta| \cdot W_2(\mu,\nu) + C_2 |\zeta| \cdot \left( \int_{M\times \mathbb{R}^d} |z_1-z_2|^2 \d \pi[(z_1,z_2),(u_1,u_2)] \right)^{1/2} \\
				&+  C_3  \left( \int_{M \times \mathbb{R}^d} |u_1-u_2|^2 \d \pi[(z_1,z_2),(u_1,u_2)] \right)^{1/2}.
			\end{align*}
			Thus, we obtained:

			$$
			\quad |B_i^2(x,\zeta) - B_i^2(y,\eta)|^2 \lesssim |\zeta|^2 \cdot \left( |x-y|^2 + W_2^2(\zeta,\eta) \right).
			$$
			In particular, for $i=0$, and by using Young's inequality we obtain 

			$$
			\langle B_0^2(x,\zeta) - B_0^2(y,\eta), u-v \rangle \lesssim |\zeta|^2 \cdot \left(|x-y|^2 + W_2^2(\zeta,\eta)\right) + |u-v|^2
			$$
			Once again, the results follows from the symmetric roles of $\zeta$ and $\eta$.

		\end{proof} 

		\subsection{Differential formulas on \texorpdfstring{$\mathscr{D}$}{D}} \label{Appendice formule geo diff}

		The proofs of the following can be found in \cite{Martin2025}.

		\begin{prop} \label{P derivee covariante pullback}
    Let $k \geq 1$, $\bm A \in \Gamma^k(T\mathscr{D})$. Then, $\bm \nabla \bm A \in \mathrm{Hom}_{\mathrm{id}}^{k-1}(T\mathscr{D}, T\mathscr{D})$ and the following holds:
    
    $$
    \bm \nabla \left( \bm \nabla\bm A\right)( B \circ \varphi) = \bm \nabla^2 \bm A (B \circ \varphi, \cdot) =  \bm \nabla^2 \bm A (\cdot, B \circ \varphi) - \bm R(\cdot, B \circ \varphi) \bm A , \quad \forall \varphi \in \mathscr{D}, \, \forall B \circ \varphi \in T_\varphi \mathscr{D},
    $$
	where $\bm \nabla^2 A(\bm B,\bm C) = \bm \nabla_{\bm B} \bm \nabla_{\bm C} \bm A - \bm \nabla_{\bm \nabla_{\bm B} \bm C}\bm A$, for $\bm B, \bm C \in \Gamma^1(T\mathscr{D})$.
\end{prop}

\begin{thm}\label{Dérivées covariantes projection horizontale et verticale} 
        The projections $P_\H$ and $P_\V$ are both elements of $\mathrm{Hom}^\infty_{\mathrm{id}}(T\mathscr{D}, T\mathscr{D})$. Moreover, for $\varphi \in \mathscr{D}$ and $ Z\circ \varphi \in \H_\varphi\mathscr{D}$, the following holds:
	    
	    \begin{align}
	        &(\bm \nabla P_{\H_\varphi}) \cdot (Z \circ \varphi) = -\mathcal{O}_{  Z \circ \varphi}\left(P_{\H_{\varphi}}(\cdot)\right) - P_{\H_{\varphi}}\left(\mathcal{O}^\top_{Z \circ \varphi}(\cdot) \right), \label{E 56} \\
	        & (\bm \nabla P_{\V_\varphi}) \cdot (Z \circ \varphi) = \mathcal{O}_{Z \circ \varphi}\left(P_{\H_{\varphi}}(\cdot)\right)  + P_{\H_{\varphi}}\left(\mathcal{O}^\top_{Z \circ \varphi}(\cdot) \right) .\label{E 57}
	    \end{align}
\end{thm}

\begin{coro}\label{Dérivée covariante Qt}
    Let $k \geq 1$, $\bm Z \in \Gamma^k(T \mathscr{D})$ be a $G$-right invariant stratified vector field on $\mathscr{D}$. Then, $\mathcal{O}_{\bm Z}$ and $\mathcal{O}^\top_{\bm Z}$ belong to $\mathrm{Hom}_{\mathrm{id}}^{k-1}(T\mathscr{D}, T\mathscr{D})$. Moreover, 
	     
	     \begin{align*}
	          (\bm\nabla \mathcal{O}_{\bm Z(\varphi)}) \cdot ( Z_2 \circ \varphi) & = \mathcal{O}_{\bm \nabla_{Z_2 \circ \varphi} \bm Z(\varphi)}\left(\cdot  \right)  -  P_{\V_{\varphi}} \left( \bm R( \cdot, Z_2 \circ \varphi ) \bm Z(\varphi) \right)\nonumber  \\
     &+ P_{\H_{\varphi}}\left( \mathcal{O}_{ Z_2\circ \varphi}^\top(\mathcal{O}_{\bm Z(\varphi)}(\cdot) \right) - \mathcal{O}_{ Z_2 \circ \varphi}(  \mathcal{O}_{\bm Z(\varphi)}(\cdot))
	     \end{align*}
	     and
	     
	     \begin{align*}
	          (\bm\nabla \mathcal{O}^\top_{\bm Z(\varphi)}) \cdot ( Z_2 \circ \varphi) & =  \mathcal{O}^\top_{\bm \nabla_{  Z_2 \circ \varphi} \bm Z(\varphi)}\left(\cdot   \right)  -    \bm R( P_{\V_{\varphi}}(\cdot),  Z_2 \circ \varphi ) \bm Z(\varphi) \nonumber  \\
     &+  \mathcal{O}_{Z_2 \circ \varphi}^\top(\mathcal{O}_{\bm Z(\varphi)}(P_{\H_{\varphi}}(\cdot)))  - \mathcal{O}^\top_{Z_2 \circ \varphi}(  \mathcal{O}^\top_{\bm Z(\varphi)}(\cdot)),
	     \end{align*}
	     for any $ Z_2 \circ \varphi \in \H_\varphi \mathscr{D}$.
\end{coro}

%

\section*{Notation table}

\noindent\textbf{Conventions.}
A bar generally denotes an object on the Wasserstein space (for example
$\grad$, $\bar{R}$, and $\se_t$), whereas a check accent is used for
objects on the product space $M\times\P$ and bold characters generally stand for objects on the diffeomorphism group $\mathscr{D}$.

\small
\begin{longtable}{>{\raggedright\arraybackslash}p{0.23\textwidth} >{\raggedright\arraybackslash}p{0.70\textwidth}}
\toprule
\textbf{Symbol} & \textbf{Meaning} \\
\midrule
\endfirsthead

\toprule
\textbf{Symbol} & \textbf{Meaning} \\
\midrule
\endhead

\midrule
\multicolumn{2}{r}{\small Continued on the next page}\\
\endfoot

\bottomrule
\endlastfoot

\multicolumn{2}{l}{\textbf{Base Riemannian geometry}}\\*[-0.2em]
$M$ & Smooth closed Riemannian manifold. In the introductory finite-dimensional discussion it may only be assumed complete. \\
$d_M(x,y)$ & Riemannian distance between $x,y\in M$. \\
$\vol$ & Riemannian volume measure on $M$, normalized by $\vol(M)=1$ in the Wasserstein sections. \\
$T_xM$, $TM$ & Tangent space at $x$ and tangent bundle of $M$. \\
$\langle\cdot,\cdot\rangle_x$, $|\cdot|_x$ & Riemannian inner product and norm on $T_xM$. \\
$\nabla$ & Levi--Civita connection on $M$; also the Riemannian gradient when applied to a scalar function. \\
$\Delta$ & Laplace--Beltrami operator on $M$. \\
$R$ & Riemann curvature operator of $M$. \\
$\Ric$ & Ricci curvature tensor of $M$. \\
$\mathrm{H}f$ & Riemannian Hessian of a smooth function $f$. \\
$\exp_x(v)$ & Riemannian exponential map at $x\in M$. \\
$\operatorname{inj}(M)$ & Injectivity radius of $M$. \\
$\Gamma(TM)$ & Space of smooth vector fields on $M$. \\
$\II$ & Second fundamental form of an isometric embedding $M\hookrightarrow\mathbb{R}^N$. \\
$\mathrm{P}_t$ & Markov semigroup of a finite-dimensional diffusion on $M$. \\

\addlinespace
\multicolumn{2}{l}{\textbf{Probability measures and Wasserstein geometry}}\\*[-0.2em]
$\P=\P(M)$ & Space of Borel probability measures on $M$, equipped with the $2$-Wasserstein distance. \\
$\P_2(M)$ & Probability measures on $M$ with finite second moment. \\
$\P_\infty$ & Smooth Wasserstein space: measures $\mu=\rho\,\vol$ with $\rho\in \C^\infty(M)$ and $\rho>0$. \\
$\P_{\mathrm{p.a.}}$ & Space of purely atomic probability measures: $\left\{ \mu \in \P : \, \mu = \sum_{i=1}^\infty \alpha_i \delta_{x_i}, \, \left\{\alpha_i\right\}_{i \in \mathbb{N}} \in \mathbb{R}_+^{\mathbb{N}}, \, \left\{x_i\right\}_{i \in \mathbb{N}} \in M^{\mathbb{N}} \right\}$. \\
$\C(\mu,\nu)$ & Set of transport plans (couplings) with marginals $\mu$ and $\nu$. \\
$W_2$ & $2$-Wasserstein distance on $\P(M)$, with cost $d_M^2$. \\
$\bm W_2$ & $2$-Wasserstein distance on $\P(\P)$, with cost $W_2$. \\
$\mathcal{W}_2$ & Wasserstein distance on $\P(M\times\P)$ with cost $d_M^2(x,y)+W_2^2(\mu,\nu)$. \\
$T_*\mu$ & Pushforward of $\mu$ by a measurable map $T$. \\
$F_f(\mu)$ & For $f\in \C^\infty(M)$, $F_f(\mu)=\int_M f\,\d\mu$ \\
$\T_\mu\P, \T\P$ & Tangent space in $\mu \in \P$ and the corresponding tangent bundle. \\
$T_\mu\P_\infty, T\P_\infty$ & Tangent space in $\mu \in \P_\infty$ and the corresponding tangent bundle. \\
$V_\varphi$ & Tangent vector represented by the gradient field $\nabla\varphi$. \\
$\Pi_\mu$, $\Pi_\mu^\perp$ & $L^2(\mu)$ orthogonal projections onto $\T_\mu\P$ and its orthogonal complement. \\
$\langle Z_1,Z_2\rangle_\mu$ & Otto metric: $\int_M\langle Z_1(x),Z_2(x)\rangle\,\d\mu(x)$. \\
$\L_ZF$ & Intrinsic Lie/directional derivative of a functional $F$ in direction $Z\in\T_\mu\P$. \\
$\grad F(\mu)$ & Wasserstein (intrinsic/Lions) gradient of $F$ at $\mu$. \\
$\grad^2F(\mu)(x,y)$ & Second derivative of $F$ with respect to the measure variable. \\
$\nabla\grad F(\mu)(x)$ & Spatial covariant derivative of the Wasserstein gradient. \\
$\C^k(\P)$, $\C^{k,\ell}(M \times \P)$ & L-differentiability classes. \\
$\Gamma^k(T\P)$ & $\C^k$ Wasserstein vector fields. \\
$\grad$ & Levi--Civita connection for the Otto metric on $\P$ or Wasserstein gradient. \\
$\bar{\he}F$ & Wasserstein Hessian of a functional $F$. \\
$\R$ & Riemann curvature operator on $\P$. \\
$\mathscr{R}(\mu,\nu)$ & Set of regular paths in $\P_\infty$ joining $\mu$ to $\nu$. \\
$\mathscr{R}_s(\mu,\nu)$ & Set of stochastic regular paths joining random measures $\mu$ and $\nu$. \\
$T_x \varphi(\cdot,\mu)$ or $\nabla \varphi(x,\mu) $ & Differential of $\varphi : M \times \P \rightarrow M$ with respect to the spatial variable. \\
$\T_\mu\varphi(x,\cdot)$ or $\grad\varphi(x,\mu)$ & Lions differential of $\varphi : M \times \P \rightarrow M$ with respect to the measure variable. \\

\addlinespace
\multicolumn{2}{l}{\textbf{Diffeomorphism group and Riemannian submersion}}\\*[-0.2em]
$\di=\mathrm{Diff}(M)$ & ILH Lie group of smooth diffeomorphisms of $M$. \\
$\di_\mu=p^{-1}(\mu)$ & Fiber over $\mu\in\P_\infty$. \\
$T_\varphi\di=\Gamma(\varphi^*TM)$ & Tangent space to $\di$ at $\varphi$. \\
$\Gamma_S^k(T\di)$ & $\C^k$ stratified vector fields on the ILH tangent bundle. \\
$\V_\varphi\di$, $\H_\varphi\di$ & Vertical and horizontal tangent subspaces at $\varphi$. \\
$P_{\H}$, $P_{\V}$ & Horizontal and vertical projections on $T\di$. \\
$\h_\varphi$ & Horizontal lift $T_{p(\varphi)}\P_\infty\to \H_\varphi\di$, with $\h_\varphi(V_\psi)=\nabla\psi\circ\varphi$. \\
$\N_\varphi$ & Normal tensor of the Riemannian submersion $p:\di\to\P_\infty$. \\
$\mathcal{O}_Z$ & Linear operator associated with the normal tensor, equivalently $\mathcal{O}_Z=P_{\V}\circ\nabla Z$ for horizontal $Z$. \\
$\langle\cdot,\cdot\rangle_\varphi$ & Right-invariant $L^2$ metric on $T_\varphi\di$. \\
$\mathrm{Hom}_\Psi^k(F^1,F^2)$ & $\C^k$ bundle morphisms covering a base map $\Psi$. \\

\addlinespace
\multicolumn{2}{l}{\textbf{Stochastic dynamics}}\\*[-0.2em]
$\se_t$ & Transition semigroup of $(\M_t)_{t\ge0}$: $\se_tF(\mu)=\mathbb E[F(\M_t(\mu))]$. \\
$\bar{\mathrm{L}}$ & Infinitesimal generator associated with $\se_t$. \\
$\vartheta\se_t$ & Law obtained by evolving an initial distribution $\vartheta\in\P(\P)$ under the semigroup. \\
$\circ D_t$, $D_t^{\bm \nabla}$ & Stratonovich and It\^o covariant differentials along a process on a vector bundle. \\
$\circ\D_t$, $\D_t^{\grad}$ & Corresponding covariant differentials along a process in $T\P_\infty$. \\
$\circ \bm D_t$, $\bm D_t^{\bm \nabla}$ & Stratonovich and It\^o covariant differentials along a process in $T\di$. \\
$\check{\mathrm{P}}_t$ & Transition semigroup of the joint process $(\varphi_t(x,\mu),\M_t(\mu))$ on $M\times\P$.

\end{longtable}
\normalsize

     \bibliographystyle{alpha}
     \bibliography{sample} 

\end{document}